\documentstyle[amscd,verbatim,leqno,12pt]{amsart}

\numberwithin{equation}{section}

\theoremstyle{plain}
\newtheorem{thm}{Theorem}[section]
\newtheorem{defn}[thm]{Definition}
\newtheorem{prop}[thm]{Proposition}
\newtheorem{lem}[thm]{Lemma}
\newtheorem{cor}[thm]{Corollary}
\newtheorem{conj}[thm]{Conjecture}

\theoremstyle{definition}
\newtheorem{rem}[thm]{Remark}
\newtheorem{exam}[thm]{Example}

\renewcommand{\b}{\bullet}
\newcommand{\beast}{\begin{eqnarray*}}
\newcommand{\east}{\end{eqnarray*}}

\newcommand{\I}{{\Bbb I}}
\newcommand{\II}{\ul{{\Bbb I}}}
\newcommand{\N}{{\Bbb N}}
\newcommand{\Z}{{\Bbb Z}}
\newcommand{\Q}{{\Bbb Q}}

\newcommand{\D}{{\Bbb D}}

\newcommand{\X}{{\Bbb X}}
\newcommand{\Af}{{\Bbb A}}
\newcommand{\fAf}{\wh{{\Bbb A}}}

\newcommand{\reldim}{{\mathrm{rel.dim}}}

\newcommand{\Spec}{{\mathrm{Spec}}\,}

\newcommand{\Spf}{{\mathrm{Spf}}\,}

\newcommand{\Proj}{{\mathrm{Proj}}\,}

\newcommand{\lra}{\longrightarrow}
\newcommand{\ra}{\rightarrow}
\newcommand{\hra}{\hookrightarrow}

\newcommand{\lla}{\longleftarrow}

\newcommand{\isom}{\overset{\sim}{=}}
\newcommand{\ti}[1]{\widetilde{#1}}

\newcommand{\ul}[1]{\underline{#1}}
\newcommand{\ol}[1]{\overline{#1}}
\newcommand{\os}{\overset}

\newcommand{\DR}{{\mathrm{DR}}}

\newcommand{\conv}{{\mathrm{conv}}}

\newcommand{\Tor}{{\mathrm{Tor}}}

\newcommand{\Ker}{{\mathrm{Ker}}}

\newcommand{\Coker}{{\mathrm{Coker}}}
\newcommand{\tor}{{\mathrm{tor}}}

\newcommand{\triv}{{\mathrm{triv}}}

\newcommand{\Mod}{{\mathrm{Mod}}}

\newcommand{\sat}{{\mathrm{sat}}}

\newcommand{\id}{{\mathrm{id}}}

\newcommand{\End}{{\mathrm{End}}}
\newcommand{\gp}{{\mathrm{gp}}}

\newcommand{\res}{{\mathrm{res}}}

\newcommand{\red}{{\mathrm{red}}}

\newcommand{\pr}{{\mathrm{pr}}}
\newcommand{\dlog}{{\mathrm{dlog}}\,}

\newcommand{\an}{{\mathrm{an}}}

\newcommand{\cB}{{\cal B}}
\newcommand{\cC}{{\cal C}}
\newcommand{\cD}{{\cal D}}
\newcommand{\cE}{{\cal E}}
\newcommand{\cF}{{\cal F}}
\newcommand{\cG}{{\cal G}}

\newcommand{\cI}{{\cal I}}

\newcommand{\cO}{{\cal O}}

\newcommand{\cS}{{\cal S}}

\newcommand{\cU}{{\cal U}}

\newcommand{\cX}{{\cal X}}
\newcommand{\cY}{{\cal Y}}
\newcommand{\cZ}{{\cal Z}}

\newcommand{\fD}{{\frak D}}

\newcommand{\fX}{{\frak X}}
\newcommand{\fY}{{\frak Y}}

\newcommand{\LNM}{{\mathrm{LNM}}}

\newcommand{\rk}{{\mathrm{rk}}\,}

\newcommand{\br}{\breve}
\newcommand{\Bl}{{\mathrm{Bl}}}

\renewcommand{\sp}{{\mathrm{sp}}}
\newcommand{\rig}{{\mathrm{rig}}}
\newcommand{\logrig}{{\mathrm{log}\text{-}\rig}}

\newcommand{\sk}{{\mathrm{sk}}}
\newcommand{\cosk}{{\mathrm{cosk}}}

\newcommand{\wh}{\widehat}

\renewcommand{\res}{{\mathrm{res}}}
\renewcommand{\d}{\dagger}

\newcommand{\dd}{\ddagger}
\newcommand{\lam}{\lambda}

\newcommand{\lr}{\langle \b \rangle}
\newcommand{\lmr}{\langle m \rangle}

\newcommand{\lur}{\langle u \rangle}

\newcommand{\dvpl}{\displaystyle\varprojlim}
\newcommand{\ve}{{\mathrm{v}}}
\newcommand{\ho}{{\mathrm{h}}}
\newcommand{\pa}{\partial}

\newcommand{\tiol}[1]{{\widetilde{\ol{#1}}}}
\newcommand{\brol}[1]{{\breve{\ol{#1}}}}

\renewcommand{\End}{{\mathrm{End}}}

\begin{document}

\title[Relative Log Convergent Cohomology III]
{Relative Log Convergent Cohomology and Relative Rigid Cohomology III}

\author{Atsushi Shiho}
\address{Graduate School of Mathematical Sciences, 
University of Tokyo, 3-8-1 Komaba, Meguro-ku, Tokyo 153-8914, JAPAN.}
                          

\begin{abstract}
In this paper, we prove the generic overconvergence of 
relative rigid cohomology with coefficient, by using 
the semistable reduction conjecture for overconvergent 
$F$-isocrystals (which is recently shown by Kedlaya). 
\end{abstract}

\maketitle

\tableofcontents


\section*{Introduction}

Let $k$ be a perfect field of characteristic $p>0$, 
let $V$ be a complete discrete valuation ring of mixed characteristic 
with residue field $k$ endowed with a lift of Frobenius $\sigma:V \lra V$ 
and let $X$ be a scheme separated of finite type
over $k$. Then, as we wrote in the introduction in the previous paper 
\cite{shiho4}, it is expected that the correct $p$-adic analogue on $X$ 
of the notion of 
the local systems (or smooth $\Q_l$-sheaves) should be the notion of 
overconvergent 
$F$-isocrystals. Based on this expectation, Berthelot conjectured 
in \cite{berthelot2} that, for a proper smooth morphism $f:X \lra Y$ 
of schemes of finite type over $k$ and an overconvergent $F$-isocrystal 
$\cE$ on $X$, the higher direct image of $\cE$ by $f$ ($=$ relative rigid 
cohomology) should have a canonical structure of an overconvergent 
$F$-isocrystal. A version of this conjecture was proved in 
\cite{shiho4}. \par 
The purpose of this paper is to consider the analogue of Berthelot's 
conjecture in the case where the given morphism $f:X \lra Y$ is no longer 
proper nor smooth. In this case, we cannnot expect that the relative 
rigid cohomology of $\cE$ has a structure of an overconvergent $F$-isocrystal 
on $Y$, because the higher direct image of an local system on $X$ 
should be only a constructible sheaf on $Y$ in this case. 
However, since a constructible sheaf is a local system on a 
dense open subset, it is natural to conjecture 
that the relative rigid cohomology of $\cE$ has a structure of 
an overconvergent $F$-isocrystal on a dense open subset of $Y$. The purpose 
of this paper is to prove that a version of this 
conjecture is true, if we admit the validity of the semi-stable reduction 
conjecture for overconvergent $F$-isocrystals, 
which was first conjectured in \cite[3.1.8]{shiho2} as a higher-dimensional 
analogue of quasi-unipotent conjecture ($=$ a $p$-adic local monodromy 
theorem of Andr\'{e}(\cite{andre}), Mebkhout(\cite{mebkhout}) and 
Kedlaya(\cite{kedlaya-p})). 
We would like to note that this result is already proved in 
the previous paper \cite{shiho4} in the case where $f$ is proper, 
without any conjectural hypothesis. So the main point in this paper is 
to treat the case where $f$ is not proper. \par 
Now we explain our conjecture and the main result in this paper in 
detail. Let us put $S:=\Spec k, \cS := \Spf V$ and assume that 
all the pairs (resp. all the triples) are separated of finite type over 
$(S,S)$ (resp. $(S,S,\cS)$.) (As for the terminology concerning 
pairs and triples, see \cite{chts}, \cite{shiho3}, \cite{shiho4}.) 
Then the precise form of our conjecture is as follows: 

\begin{conj}\label{b+conj}
Assume given a morphism of pairs $f:(X,\ol{X}) \lra (Y,\ol{Y})$ such that 
$\ol{X} \lra \ol{Y}$ is proper. Then there exist an open dense subset 
$U \subseteq Y$ and a subcategory $\cC$ of the category of 
$(U,\ol{Y})$-triples over $(S,S,\cS)$ such that, for any overconvergent 
$F$-isocrystal $\cE$ on $(X,\ol{X})/\cS_K$ and for any $q \geq 0$, 
there exists uniquely an overconvergent isocrystal $\cF$ satisfying the 
following condition$:$ For any $(Z,\ol{Z},\cZ) \in \cC$ such that 
$\cZ$ is formally smooth over $\cS$ on a neighborhood of $Z$, the restriction 
of $\cF$ to $I^{\d}((Z,\ol{Z})/\cS_K,\cZ)$ is given 
functorially by 
$(R^qf_{(X \times_Y Z,\ol{X} \times_{\ol{Y}} \ol{Z})/\cZ, \rig *}\cE, 
\epsilon)$, where $\epsilon$ is given by 
$$ 
p_2^*R^qf_{(X \times_Y Z,\ol{X} \times_{\ol{Y}} \ol{Z})/\cZ, \rig *}\cE 
\os{\simeq}{\rightarrow} 
R^qf_{(X \times_Y Z,\ol{X} \times_{\ol{Y}} \ol{Z})/\cZ \times_{\cS} \cZ, 
\rig *}\cE \os{\simeq}{\leftarrow} 
p_1^*R^qf_{(X \times_Y Z,\ol{X} \times_{\ol{Y}} \ol{Z})/\cZ, \rig *}\cE. 
$$ 
$($Here $I^{\d}((Z,\ol{Z})/\cS_K, \cZ)$ denotes the category of 
overconvergent isocrystals on $(Z,\ol{Z})\allowbreak / \allowbreak \cS_K$ 
over $\cZ$, 
$R^qf_{(X \times_Y Z,\ol{X} \times_{\ol{Y}} \ol{Z})/\cZ, \rig *}\cE, 
R^qf_{(X \times_Y Z,\ol{X} \times_{\ol{Y}} \ol{Z})/\cZ \times_{\cS} \cZ, 
\rig *}\cE$ are relative rigid cohomologies and 
$p_i$ is the morphism 
$]\ol{Z}[_{\cZ \times_{\cS} \cZ} \lra \,]\ol{Z}[_{\cZ}$ induced by the 
$i$-th projection.$)$ 
Also, $\cF$ admits a Frobenius structure which is induced by the 
Frobenius structure of $\cE$. 
\end{conj} 

\begin{rem} 
\begin{enumerate}
\item 
One can also consider a stronger version, which requires the category 
$\cC$ to be the category of all the $(Y,\ol{Y})$-triples over $(S,S,\cS)$. 
\item 
If $f$ is strict (that is, $f^{-1}(Y)=X$ holds), then Conjecture \ref{b+conj} 
is shown in \cite{shiho4}. 
\end{enumerate}
\end{rem} 

Then the main result in this paper is as follows (for more 
precise form, see Corollary \ref{maincor}(1)): 

\begin{thm}\label{mainthm0}
Conjecture \ref{b+conj} is true. 
\end{thm}

In the proof of Theorem \ref{mainthm0}, we use the semi-stable reduction 
conjecture of overconvergent $F$-isocrystals, which is recently proved 
by Kedlaya (\cite{kedlayaI}, \cite{kedlayaII}, \cite{kedlayaIII}, 
\cite{kedlayaIV}) for $\cE$. In fact, this is the only place where we use 
the Frobenius structure on $\cE$. So 
we can also form a variant of Theorem \ref{mainthm0} without 
using Frobenius structure, 
by introducing the notion of `potentially semi-stable overconvergent 
isocrystals'. (See Theorem \ref{mainthm}.) We remark also that 
we can state the theorem in such a way that claims 
`a kind of constructibility' of relative rigid cohomology. 
(See Corollary \ref{constr}, Corollary \ref{maincor}(2).) \par 

Now we give an outline of the proof of Theorem \ref{mainthm0} and 
explain the content of this paper. Roughly speaking, the 
semi-stable reduction conjecture for overconvergent $F$-isocrystals 
claims that any overconvergent $F$-isocrystal should be 
extendable to the boundary logarithmically after taking pull-back 
by a certain proper surjective generically etale morphism. 
So, as a first step, we consider the coherence and the overconvergence 
of relative rigid cohomology in 
the case where the morphism $f:\ol{X} \lra \ol{Y}$ admits a nice log 
structure, that is, the case where $f$ comes from a nice 
morphism of log schemes $(\ol{X},M_{\ol{X}}) \lra (\ol{Y},M_{\ol{Y}})$. 
In Section 1, we give a review of the residues of 
log-$\nabla$-modules and isocrystals which are developed in 
\cite{kedlayaI} and prove some preliminary facts which we need 
in Section 2, where we prove the coherence and the 
overconvergence of the relative 
rigid cohomology of $(X,\ol{X}) \lra (Y,\ol{Y})$ in the above-mentioned case. 
Let $D$ be the closure of $f^{-1}(Y)-X$ in $\ol{X}$. 
(In `nice case' we are treating, 
$\ol{X}$ is smooth over $k$ and $D$ 
is a normal crossing divisor in $\ol{X}$.) 
In this case, we prove the overconvergence of the relative 
rigid cohomology of $(X,\ol{X}) \lra (Y,\ol{Y})$ 
by relating it to the relative 
log analytic cohomology of $(\ol{X},M_{\ol{X}}) \lra (\ol{Y},M_{\ol{Y}})$ 
via an intermediate cohomology, the relative log rigid cohomology of 
$((\ol{X}-D,M_{\ol{X}}|_{\ol{X}-D}),(\ol{X},M_{\ol{X}})) \lra 
(\ol{Y},M_{\ol{Y}})$: We compare the relative 
rigid cohomology of $(X,\ol{X}) \lra (Y,\ol{Y})$ with 
the relative log rigid cohomology of 
$((\ol{X}-D,M_{\ol{X}}|_{\ol{X}-D}),(\ol{X},M_{\ol{X}})) \lra 
(\ol{Y},M_{\ol{Y}})$ by the method which is 
a generalization to the relative case of what we have developped in 
\cite[2.4]{shiho2}, and we relate the relative log rigid cohomology of 
$((\ol{X}-D,M_{\ol{X}}|_{\ol{X}-D}),(\ol{X},M_{\ol{X}})) \lra 
(\ol{Y},M_{\ol{Y}})$ to the relative 
rigid cohomology of $(X,\ol{X}) \lra (Y,\ol{Y})$ by the method 
we have developped in \cite[\S 5]{shiho3}. In Section 3, we prove the 
invariance of relative log analytic cohomology under log blow-ups. 
We need this property in the proof of the main theorem, 
mainly by technical reason. In Section 4, 
we prove a result of altering a given 
proper morphism $\ol{X} \lra \ol{Y}$ 
of schemes over $k$ endowed with an open subscheme $X \subseteq \ol{X}$ 
to a certain simplicial morphism 
of schemes whose components admit nice log structures, by 
using results of de Jong \cite{dejong1}, \cite{dejong2}. 
The results in this section is a slight generalization of what we have 
shown in \cite[\S 6]{shiho4}. This result allows us 
to reduce the proof of the main theorem to the case we 
considered in Section 2. 
In Section 5, we prove the main result by combining the results in 
the previous sections, by the method analogous to \cite[\S 7]{shiho4}. \par 
The author is partly supported by Grant-in-Aid for Young Scientists (B), 
the Ministry of Education, Culture, Sports, Science and Technology of 
Japan and JSPS Core-to-Core program 18005 whose representative is 
Makoto Matsumoto. \par 

\section*{Convention}
\hspace{-16pt}
(1) \,\, Throughout this paper, $k$ is a field of characteristic 
$p>0$, $V$ is a complete discrete valuation ring of mixed characteristic 
with residue field $k$, and $K$ is the fraction field of $V$. 
Note that we will assume that $k$ is perfect in some results. 
In particular, as we will write in the text, we will assume 
that $k$ is perfect in Sections 4, 5. 
We put $S := \Spec k$ and $\cS:=\Spf V$. Let us fix a power $p^a$ of $p$ 
and assume that we have an endomorphism $\sigma: V \lra V$ which is a 
lift of $p^a$-th power Frobenius endomorphism on $k$. 
(However, we need the existence of $\sigma$ 
only when we speak of Frobenius structures on isocrystals.) 
All the log formal schemes are fine log 
(not necessarily $p$-adic) formal schemes which are separated and 
topologically of finite type over $\cS$. We call such a log formal scheme 
a fine log formal $\cS$-scheme. If it is $p$-adic, we call it a 
$p$-adic fine log formal $\cS$-scheme. If it is a log scheme, we call it 
a fine log $\cS$-scheme, and if it is defined over $S$, we call it a fine 
log formal $S$-scheme. If the log structure is fs, we replace `fine' by 
`fs' and if the log structure is trivial, we omit the word `fine log'. \par 
(2) \,\, For a formal $\cS$-scheme $T$, we denote the rigid analytic space 
assocated to $T$ by $T_K$. \\
(3) \,\, In this paper, we freely use the 
terminologies concerning log structures 
defined in \cite{kkato1}, \cite{shiho1} and \cite{shiho2}. 
For a fine log (formal) scheme $(X,M_X)$, $(X,M_X)_{\triv}$ denotes the 
maximal open sub (formal) scheme of $X$ on which the log structure $M_X$ is 
trivial. A morphism $f:(X,M_X) \lra (Y,M_Y)$ is said to be strict if 
$f^*M_Y = M_X$ holds. \\ 
(4) \,\, In this paper, we freely use the notations and 
terminologies concerning (log) pairs, triples, 
isocrystals on relative log convergent site, overconvegent isocrystals, 
relative log analytic cohomologies, relative rigid cohomologies and 
relative log rigid cohomologies which are developped in the previous 
papers \cite{shiho3}, \cite{shiho4}. All the pairs in this paper are 
separated of finite type over $(S,S)$ and all the triples in this paper
are separated of finite type over $(S,S,\cS)$. Frobenius structures on 
isocrystals are always the ones with repsect to $\sigma$ defined in (1). \\
(5) \,\, Fiber products of log formal schemes are taken in the category 
of fine log formal schemes.


\section{Residues of log-$\nabla$-modules and isocrystals}

In this section, we give a review of definitions and results by 
Kedlaya \cite{kedlayaI} on log-$\nabla$-modules and isocrystals and 
prove some preliminary results on the residue of log-$\nabla$-modules and 
isocrystals. Finally, we recall the semi-stable reduction conjecture for 
overconvergent $F$-isocrystals. \par 
First we recall the definition of log-$\nabla$-modules and 
the residues of them (\cite[2.3]{kedlayaI}). 
Let $f:\fX \lra \fY$ be a morphism of rigid analytic spaces over $K$ and 
let $x_1, \cdots, x_r \in \Gamma(\fX,\cO_{\fX})$. Then we define 
the relative differential module $\Omega^1_{\fX/\fY}$ by 
$\Omega^1_{\fX/K}/f^*\Omega^1_{\fY/K}$ (where $\Omega^1_{\fX/K}, 
\Omega^1_{\fY/K}$ denotes the sheaf of continuous differential forms of 
$\fX, \fY$, respectively) and we define the relative log differential 
module $\omega^1_{\fX/\fY}$ (which is denoted by $\Omega^{1,\log}_{\fX/\fY}$ 
in \cite{kedlayaI}) by 
$$ 
\omega^1_{\fX/\fY} := (\Omega^1_{\fX/\fY} \oplus 
\bigoplus_{i=1}^r \cO_{\fX} \dlog x_i )/N, $$
where $N$ is the sub $\cO_{\fX}$-module locally generated by 
$(dx_i,0)-(0,x_i\dlog x_i)$. A log-$\nabla$-module $(E,\nabla)$ on 
$\fX$ with respect to $x_1,\cdots, x_r$ relative to $\fY$ is 
a locally free $\cO_{\fX}$-module $E$ endowed with an integrable log 
connection $\nabla: E \lra E \otimes_{\cO_{\fX}} \omega^1_{\fX/\fY}$. 
When $\fY=\cS_K$ holds, we omit the word `relative to $\fY$'. 
If we put $\fD_i := \{x_i=0\} \subseteq \fX$ and 
$M_i := {\mathrm{Im}}(\Omega^1_{\fX/\fY}\oplus 
\bigoplus_{j \not=i} \cO_{\fX}\dlog x_j \lra \omega^1_{\fX/\fY})$ \, 
$(1 \leq i \leq r)$, the composite 
$$ 
E \os{\nabla}{\lra} E \otimes_{\cO_{\fX}} \omega^1_{\fX/\fY} 
\lra E \otimes_{\cO_{\fX}} (\omega^1_{\fX/\fY}/M_i) 
\cong E \otimes_{\cO_{\fX}} \cO_{\fD_i}\dlog x_i = 
E \otimes_{\cO_{\fX}} \cO_{\fD_i} 
$$ 
naturally induces an element of 
$\End_{\cO_{\fD_i}}(E \otimes_{\cO_{\fX}} \cO_{\fD_i})$. 
We call this element the residue of $(E,\nabla)$ along $\fD_i$. 
$(E,\nabla)$ is said to have nilpotent residues if 
the residue of $(E,\nabla)$ along $\fD_i$ is nilpotent for all 
$1 \leq i \leq r$. We denote the category of log-$\nabla$-modules on 
$\fX$ with respect to $x_1,\cdots, x_r$ relative to $\fY$ having 
nilpotent residues by $\LNM_{\fX/\fY}$. When $\fY=\cS_K$ holds, 
we write $\LNM_{\fX}$ instead of $\LNM_{\fX/\fY}$. 

\begin{rem} 
The log differential module $\omega^1_{\fX/\fY}$ and the residues 
of log-$\nabla$-modules are unchanged if we replace $x_i$'s by 
$u_ix_i$'s for some $u_i \in \Gamma(\fX,\cO_{\fX}^{\times}) \, (1 \leq i 
\leq r)$. 
\end{rem} 

Following \cite[3.1]{kedlayaI}, we call a subinterval $I$ of $[0,\infty)$ 
aligned if any endpoint at which it is closed is in $p^{\Q} \cup \{0\}$ 
(it is equal to $\Gamma^*$ in \cite{kedlayaI} in our case) and for 
an aligned interval $I$, we define the polyannulus $A_K^n(I)$ by 
$$ A_K^n(I) := \{(t_1,\cdots,t_n) \,\vert\, t_i \in I \,(1 \leq i \leq n)\}.$$ 

Let $\fX$ be a smooth rigid space endowed with sections 
$x_1,\cdots, x_m \in \Gamma(\fX,\cO_{\fX})$ whose zero loci are smooth and 
meet transversally. For an aligned interval $I$ with $0 \notin I$ and 
$\eta \in [0,\infty)$, a log-$\nabla$-module $(E,\nabla)$ on 
$\fX \times A^n_K(I)$ with respect to $x_1,\cdots, x_m$ is called 
$\eta$-convergent with respect to $t_1,\cdots, t_n$ 
if, for any $e \in \Gamma(\fX \times A^n_K(I),E)$, the multisequence 
$$ \left\| \dfrac{1}{i_1! i_2! \cdots i_n!} \dfrac{\pa^{i_1}}{\pa t_1^{i_1}} 
\cdots 
\dfrac{\pa^{i_n}}{\pa t_n^{i_n}}(e) \right\| \eta^{i_1+\cdots +i_n} $$ 
converges to zero. (cf. ~\cite[2.4.2]{kedlayaI}. We slightly changed 
the terminology in order to make it shorter.) 
For $a \in p^{\Q_{<0}} \cup \{0\}$, a 
log-$\nabla$-module $(E,\nabla)$ on $\fX \times A^n_K[a,1)$ or 
$\fX \times A^n_K(a,1)$ with respect to $x_1,\cdots, x_m, t_1,\cdots, t_n$ 
is called convergent if, for any $\eta \in (0,1)$, there exists 
some $b \in (a,1) \cap p^{\Q}$ such that, for any $c \in [b,1) \cap p^{\Q}$, 
$(E,\nabla)$ is $\eta$-convergent with respect to $t_1,\cdots, t_n$ 
on $\fX \times A^n_K[b,c]$ (\cite[3.6.6]{kedlayaI}). \par 
Let $\fX,x_1,\cdots,x_m$ be as above and let $I$ be an aligned interval. 
Then $(E,\nabla) \in \LNM_{\fX \times A^n_K(I)}$ (here $\fX \times A^n_K(I)$ 
is endowed with sections $x_1,\cdots,x_m,t_1,\cdots,t_n$) is unipotent 
relative to $\fX$ if there exists a filtration of $(E,\nabla)$ by 
sub log-$\nabla$-modules such that each graded quotient is again 
a log-$\nabla$-module and has the form $\pi^*(F,\nabla_F)$ 
($\pi$ denotes the projection $\fX \times A^n_K(I) \lra \fX$) for 
some log-$\nabla$-module $(F,\nabla_F)$ on $\fX$ with respect to 
$x_1,\cdots,x_m$ (\cite[3.2.5]{kedlayaI}). \par 
Then, by \cite[3.6.2, 3.6.9]{kedlayaI}, we have the following: 

\begin{prop}\label{nilpunip}
Let $a \in p^{\Q_{<0}}$, 
let $\fX$ be a smooth rigid space endowed with sections 
$x_1,\cdots, x_m \in \Gamma(\fX,\cO_{\fX})$ whose zero loci are smooth and 
meet transversally and let $(E,\nabla)$ be an object in 
$\LNM_{\fX \times A^n_K(a,1)}$ which is convergent. 
Then $(E,\nabla)$ is unipotent relative to $\fX$ if and only if 
it extends to an object $(\ti{E},\ti{\nabla})$ 
in $\LNM_{\fX \times A^n_K[0,1)}$. Moreover, the extension is unique and 
$(\ti{E},\ti{\nabla})$ is also unipotent relative to $\fX$. 
\end{prop} 

Next, let $X$ be a smooth scheme over $k$, let 
$D := \bigcup_{i=1}^r D_i$ be a simple normal crossing divisor in $X$ 
(each $D_i$ is assumed to be smooth) and let $M_X$ be the log structure 
on $X$ associated to $D$. Let $\cE$ be a locally free isocrystal on 
log convergent site $(X/\cS)^{\log}_{\conv} = ((X,M_X)/\cS)_{\conv}$. 
Zariski locally on $X$, we can form the Cartesian diagram 
\begin{equation}\label{liftable}
\begin{CD}
D @>{\subset}>> \cD \\ 
@V{\cap}VV @V{\cap}VV \\ 
X @>{\subset}>> \cX, 
\end{CD}
\end{equation} 
where all arrows are closed immersions, $\cX$ is a $p$-adic formal $\cS$ 
scheme formally smooth over $\cS$ with $\cX \times_{\cS} S = X$ and 
$\cD = \bigcup_{i=1}^r \cD_i$ is a relative simple normal crossing 
divisor on $\cX$ (each $\cD_i$ is assumed to be formally smooth over $\cS$) 
such that $\cD_i \times_{\cS} S = D_i$ holds and that each 
$\cD_i$ is defined by a single element $t_i \in \Gamma(\cX,\cO_{\cX})$. 
Then, $\cE$ induces a log-$\nabla$-module on 
$\cX_K = \,]X[_{\cX}$ with respect to 
$t_1,\cdots, t_r$. $\cE$ is said to have nilpotent residues 
if there exists an open covering $X=\bigcup_j X_j$ such that 
each $X_j$ admits a diagram like \eqref{liftable} (with 
$X$ replaced by $X_j$) and that the induced log-$\nabla$-module 
on $\cX_K$ has nilpotent residues (\cite[6.4.4]{kedlayaI}). \par 
Keep the notation in the diagram \eqref{liftable} and let $\cE$ be 
a locally free 
isocrystal on $(X/\cS)^{\log}_{\conv}$ having nilpotent residues 
such that the induced log-$\nabla$-module $(E,\nabla)$ on 
$\cX_K=\,]X[_{\cX}$ 
with respect to $t_1,\cdots, t_r$ has nilpotent residues. 
For a subset $I$ of $\{1,\cdots,r\}$, we define 
$D_I := \bigcap_{i \in I}D_i, \cD_I := \bigcap_{i \in I}\cD_i$. 
Then, we can restrict $(E,\nabla)$ to 
$$ ]D_I[_{\cX} \,\cong\, ]D_I[_{\cD_I} \times A_K^{|I|}[0,1) (= 
\cD_{I,K} \times A_K^{|I|}[0,1)). $$
(The coordinate of $A_K^{|I|}[0,1)$ is given by $t_i \,(i \in I)$.) 
Then we have the following: 

\begin{prop}\label{isocunip}
$(E,\nabla)|_{]D_I[_{\cD_I} \times A_K^{|I|}[0,1)}$ is unipotent 
relative to $]D_I[_{\cD_I}$. 
\end{prop} 

\begin{pf} 
We may work Zariski locally on $\cX$. So we may assume that 
$\cX:=\Spf A$ is affine and $\cX$ admits a formally etale morphism 
$f: \cX \lra \Spf V \{t_1, \cdots, t_n\}$ for some $n \geq r$ 
such that the image of $t_i \in V \{t_1, \cdots, t_n\}$ in 
$\Gamma(\cX,\cO_{\cX})$ is the element $t_i$ defined above for $1 \leq i 
\leq r$. Let $M_{\cX}$ be the log structure on $\cX$ associated to $\cD$ and 
let $(\cX(1),M_{\cX(1)}) := (\cX,M_{\cX}) \times_{\cS} 
(\cX,M_{\cX})$. Then we have the morphism 
$$f(1): \cX(1) \lra \Spf V\{t^{(1)}_1, \cdots, t^{(1)}_n, 
t^{(2)}_1, \cdots, t^{(2)}_n\}$$ 
naturally induced by $f$. Let us define $\cX(1)'$ by 
$$ 
\cX(1)' := \cX(1) \times_{\Spf V\{t^{(1)}_1, \cdots, t^{(1)}_n, 
t^{(2)}_1, \cdots, t^{(2)}_n\}} 
\Spf V\{t^{(1)}_1, \cdots, t^{(1)}_n, 
u_1^{\pm 1}, \cdots, u_n^{\pm 1}\}, $$
where $u_i := t^{(2)}_i/t^{(1)}_i \,(1 \leq i \leq r)$, $u_i := 
t^{(2)}_i +1 \,(r+1 \leq i \leq n)$ and let $M_{\cX(1)'}$ be the log structure 
on $\cX(1)'$ associated to the pre-log structure 
$\N^r \lra \cO_{\cX(1)'}; \, (m_1,\cdots,m_r) \mapsto 
{t^{(1)}_1}^{m_1}\cdots {t^{(1)}_r}^{m_r}$. Then the closed immersion 
$(X,M_X) \hra (\cX(1),M_{\cX(1)})$ (induced by $X \hra \cX$ and the 
diagonal map $\cX \hra \cX(1)$) naturally factors as 
$$ (X,M_X) \hra (\cX(1)',M_{\cX(1)'}) \lra (\cX(1),M_{\cX(1)}), $$
where the first map is an exact closed immersion and the second map 
is formally log etale. So we have 
$]X[^{\log}_{\cX(1)}=]X[_{\cX(1)'}$. \par 
Let $p_i:]X[^{\log}_{\cX(1)}=]X[_{\cX(1)'} 
\lra ]X[_{\cX} \,(i=1,2)$ be the projection. 
Then $\cE$ induces the isomorphism 
$\epsilon: p_2^*E \os{\cong}{\lra} p_1^*E$. 
If we denote the completion of $\cX(1)'$ along $X$ by 
$\wh{\cX(1)}'$, the morphism $\wh{\cX(1)}' \lra \cX$ induced by the first 
projection $\cX(1) \lra \cX$ 
has the form $\Spf A[[u_1-1,\cdots,u_n-1]] \lra \Spf A$. So, the morphism 
$p_1$ can be identified with the projection 
$$ ]X[_{\cX} \times A^n_K[0,1) \lra ]X[_{\cX}, $$
where the coordinate of $A_K^n[0,1)$ is given by $u_i-1 \,(1 \leq i \leq n)$. 
By this identification, the isomorphism $\epsilon$ can be written as 
$$ 1 \otimes e \mapsto 
\sum_{i_1,\cdots,i_n=0}^{\infty} 
\left( \prod_{j=1}^r \prod_{l=0}^{i_j-1} \left( t_j \dfrac{\pa}{\pa t_j} - l 
\right) \prod_{j=r+1}^n \dfrac{\pa^{i_j}}{\pa t^{i_j}_j} 
(e)\right) 
\otimes 
\left( \prod_{j=1}^n \dfrac{(u_j-1)^{i_j}}{i_j!} \right) $$ 
by \cite[6.4.1]{kedlayaI}, \cite[6.7.1]{kkato1}. 
(There is a slight mistake in \cite[6.4.1]{kedlayaI}.) 
So the multisequence 
$$ \left\| \dfrac{1}{i_1! \cdots i_n!} 
\left( \prod_{j=1}^r \prod_{l=0}^{i_j-1} \left( t_j \dfrac{\pa}{\pa t_j} - l 
\right) \prod_{j=r+1}^n \dfrac{\pa^{i_j}}{\pa t^{i_j}_j} \right)
(e) \right\| \xi^{i_1+\cdots +i_n}
$$ 
converges to zero for any $e \in \Gamma(]X[_{\cX},E)$ and $\xi 
\in [0,1)$. \par 
Now we consider the restriction of $(E,\nabla)$ to 
$]D_I[_{\cX} \cong ]D_I[_{\cD_I} \times A_K^{|I|}[0,1)$. 
For any $\eta \in [0,1)$, take any $c \in (\sqrt{\eta},1) \cap p^{\Q}$. 
Then, on $]D_I[_{\cD_I} \times A_K^{|I|}[\sqrt{\eta},c]$, 
we have the following inequality for any $i_j \in \N \,(j \in I)$ and 
$e \in \Gamma(]X[_{\cX},E)$: 
\begin{align*}
\left\| \dfrac{1}{\prod_{j \in I}i_j!}
\prod_{j \in I} \prod_{l=0}^{i_j-1} \left( t_j \dfrac{\pa}{\pa t_j} - l 
\right) (e) \right\| \sqrt{\eta}^{\sum_{j \in I}i_j} 
& = 
\left\| \dfrac{1}{\prod_{j \in I}i_j!} \prod_{j \in I}
\dfrac{\pa^{i_j}}{\pa t_j^{i_j}}(e)\right\| \left\| \prod_{j \in I}t_j^{i_j} 
\right\| \sqrt{\eta}^{\sum_{j \in I}i_j} \\ 
& \geq 
\left\| \dfrac{1}{\prod_{j \in I}i_j!} \prod_{j \in I}
\dfrac{\pa^{i_j}}{\pa t_j^{i_j}}(e)\right\| \eta^{\sum_{j \in I}i_j}. 
\end{align*}
Since the multisequence (with respect to $i_j\,(j \in I)$) on the left hand 
side of the above inequality converges to zero, the multisequence 
on the right hand side of the above inequality also converges to zero. 
Since the trivial log-$\nabla$-module $(\cO,d)$ on 
$]D_I[_{\cD_I} \times A_K^{|I|}[\sqrt{\eta},c]$ is $\eta$-convergent 
with respect to $t_j\,(j \in I)$ by \cite[3.6.1]{kedlayaI}, 
we can show (by using Leibniz rule) that the multisequence 
on the right hand side of the above inequality converges to zero 
for any $e \in \Gamma(]D_I[_{\cD_I} \times A_K^{|I|}[\sqrt{\eta},c],E)$. 
So the log-$\nabla$-module 
$(E,\nabla)|_{]D_I[_{\cD_I} \times A_K^{|I|}[0,1)}$ 
is convergent. Since it has nilpotent residues, it is unipotent 
relative to $]D_I[_{\cD_I}$ by Proposition \ref{nilpunip}. 
\end{pf} 

We define the notion of `having nilpotent residues' for isocrystals 
in a slightly generalized situation. 

\begin{defn}\label{nilpres}
Let $X$ be a smooth scheme over $k$ and 
let $D$ be a normal crossing divisor on $X$. $(D$ is not necessarily a 
simple normal crossing divisor, that is, some irreducible component of 
$D$ may have self-intersection.$)$ A locally free 
isocrystal $\cE$ on $(X/\cS)^{\log}_{\conv}$ 
is said to have nilpotent residues 
if there exists an etale covering $\coprod_j X_j \lra X$ such that 
each $X_j$ admits a diagram like \eqref{liftable} $($with 
$X$ replaced by $X_j)$ and that the induced log-$\nabla$-module 
on $\cX_K$ has nilpotent residues. 
\end{defn} 

We prove 
the compatibility of the above definition with \cite[6.4.4]{kedlayaI} and 
the fact 
that the above definition is independent of the choice of 
an etale covering $\coprod X_j \lra X$ and diagrams like \eqref{liftable}. 
First we recall the notion of admissible closed immersion 
(\cite[2.1.7]{nakkshiho}). (We also introduce the notion of 
strongly admissible closed immersion.)

\begin{defn} 
Let $X$ be a smooth scheme over $k$, let $D = \bigcup_{i=1}^rD_i$ 
be a simple normal crossing divisor $($where $D_i$ are smooth 
divisors$)$
and denote the log structure on $X$ associated to $D$ by $M_X$. 
Then, an exact
closed immersion $(X,M_X) \hra (\cX,M_{\cX})$ into a $p$-adic fine log 
formal $\cS$-scheme $(\cX,M_{\cX})$ is called an admissible closed immersion 
if $\cX$ is formally log smooth over $\cS$, 
the log structure $M_{\cX}$ is induced by a relative 
simple normal crossing divisor 
$\cD = \bigcup_{i=1}^r\cD_i$ $($where $\cD_i$ are formally smooth over 
$\cS)$ such that $D_i = X \times_{\cX} \cD_i$ holds. Moreover, it is 
called a strongly admissible closed immersion if each $\cD_i$ is defined by 
a single element $t_i \in \Gamma(\cX,\cO_{\cX})$. 
\end{defn} 

Let the notation be as above and assume that 
$(X,M_X) \hra (\cX,M_{\cX})$ is a strongly admissible closed immersion. 
Then an isocrystal $\cE$ on $(X/\cS)^{\log}_{\conv}$ induces a 
log-$\nabla$-module on $]X[_{\cX}$ with respect to $t_1,\cdots,t_r$. 

\begin{lem}\label{lem1}
Let $X$ be a smooth scheme over $k$, let $D = \bigcup_{i=1}^rD_i$ 
be a simple normal crossing divisor $($where $\cD_i$ are smooth 
divisors$)$
and denote the log structure on $X$ associated to $D$ by $M_X$. 
Let us assume given the following diagram
\begin{equation}\label{lem1diag}
\begin{CD}
(X',M_{X'}) @>{i'}>> (\cX',M_{\cX'}) \\ 
@VfVV @VgVV \\
(X,M_X) @>i>> (\cX,M_{\cX}), 
\end{CD}
\end{equation}
where $f$ is a strict etale morphism, $i,i'$ are strictly admissible 
closed immersions and $M_{\cX}$ $($resp. $M_{\cX'})$ is induced by 
a relative simple normal crossing divisor $\cD = \bigcup_{i=1}^r \cD_i$ 
$($resp. $\cD' := \bigcup_{i=1}^r\cD'_i)$ with 
$g^*\cD_i=\cD'_i$ and let $t_i \in \Gamma(\cX,\cO_{\cX}) \,(1 \leq i \leq r)$ 
be an element defining $\cD_i$. Let $(E,\nabla)$ be a log-$\nabla$-module on 
$]X[_{\cX}$ with respect to $t_1,\cdots,t_r$ and denote the pull-back 
of $(E,\nabla)$ by the morphism $]X'[_{\cX'} \lra ]X[_{\cX}$ induced by $g$ 
by $(E',\nabla')$. Then we have the following$:$ 
\begin{enumerate} 
\item 
If $(E,\nabla)$ has nilpotent residues, $(E',\nabla')$ also has 
nilpotent residues. 
\item 
If the homomorphisms $\Gamma(]D_i[_{\cD_i},\cO_{]D_i[_{\cD_i}}) 
\lra \Gamma(]D'_i[_{\cD'_i},\cO_{]D'_i[_{\cD'_i}}) \, (1 \leq i \leq r)$ 
induced by $g$ are all injective, the converse is also true. 
\end{enumerate}
\end{lem} 

\begin{pf} 
The residue of $(E,\nabla)$ along $]D_i[_{\cD_i}$ ($=$ zero locus of $t_i$) 
is sent to the residue of $(E',\nabla')$ along 
$]D'_i[_{\cD'_i}$ by the homomorphism 
$$\End_{\cO_{]D_i[_{\cD_i}}}(E|_{]D_i[_{\cD_i}}) \lra 
\End_{\cO_{]D'_i[_{\cD'_i}}}(E|_{]D'_i[_{\cD'_i}})$$ induced by $g$, 
and this map is injective if the map 
$$\Gamma(]D_i[_{\cD_i},\cO_{]D_i[_{\cD_i}}) 
\lra \Gamma(]D'_i[_{\cD'_i},\cO_{]D'_i[_{\cD'_i}})$$ is injective. 
The assertion of the lemma is easily deduced by these facts. 
\end{pf} 

\begin{lem}\label{lem2}
The condition in Lemma \ref{lem1}$($2$)$ is satisfied in the following cases. 
\begin{enumerate} 
\item 
The case where $f$ is the identity map 
and $g$ is strict formally smooth. 
\item 
The case $g$ is formally etale, the diagram \eqref{lem1diag} is Cartesian 
and the images of 
the morphisms $D'_i:=X' \times_{X}D_i \lra D_i$ induced by $f$ are 
dense for all $1 \leq i \leq r$. 
\end{enumerate}
\end{lem} 

\begin{pf} 
First we prove (1). We may work Zariski locally. 
Let $\wh{\cD}_i, \wh{\cD}'_i$ be the completion of 
$\cD_i,\cD'_i$ along $D_i$ and assume that $\wh{\cD}_i=\Spf A$ is affine. 
Then, by \cite[2.31]{shiho3}, we have the 
isomorphism $\wh{\cD}'_i \isom \Spf A[[x_1,\cdots,x_n]]$ for some $n$ 
Zariski locally on $\wh{\cD}_i$. So the map 
\begin{align*}
\Gamma(]D_i[_{\cD_i},\cO_{]D_i[_{\cD_i}}) 
& = \Gamma(\wh{\cD}_{i,K}, \cO_{\wh{\cD}_{i,K}}) 
\lra 
\Gamma(\wh{\cD}_{i,K} \times A^n_K[0,1), 
\cO_{\wh{\cD}_{i,K} \times A^n_K[0,1)}) \\ & = 
\Gamma(\wh{\cD}'_{i,K}, \cO_{\wh{\cD}'_{i,K}}) = 
\Gamma(]D'_i[_{\cD'_i},\cO_{]D'_i[_{\cD'_i}})
\end{align*}
is injective. \par 
Next we prove (2). We may assume that $\cD_i$ is affine. 
Let $\wh{\cD}_i, \wh{\cD}'_i$ be the completion of 
$\cD_i,\cD'_i$ along $D_i,D'_i$ respectively. 
Zariski locally on $\cD_i$, we can take a Cartesian 
diagram 
\begin{equation*}
\begin{CD}
D_i @>{\subset}>> \cD_i \\ 
@VVV @VVV \\ 
\Spec k[x_1,\cdots,x_m] @>{\subset}>> \Spf V\{x_1,\cdots,x_n\} 
\end{CD}
\end{equation*} 
for some $m\leq n$, where 
the right vertical arrow is formally etale and 
the lower horizontal arrow is defined by 
$x_{m+1}=\cdots =x_n=0$ and the closed immersion $\Spec k \hra \Spf V$. 
Let us define $\cD_{i,0} \subset \cD_i$ as the zero locus 
$x_{m+1}=\cdots =x_n=0$. Then $\cD_{i,0}$ 
is formally smooth over $\cS$ and $D_i = \cD_{i,0} \times_{\cS} S$ 
holds. Moreover, we have the isomorphism 
$\wh{\cD}_i \cong \Spf A[[x_{m+1},\cdots,x_n]]$, where $\cD_{i,0}:=\Spf A$. 
So we have a retraction $\wh{\cD}_i \lra \cD_{i,0}$ defined by the natural 
inclusion $A \hra A[[x_{m+1},\cdots,x_n]]$. Let 
$\cD'_{i,0} \lra \cD_{i,0}$ be the unique formally etale morphism 
lifting $D'_i \lra D_i$. Then both $\wh{\cD'}_i \lra \wh{\cD}_i$ and 
$\cD'_{i,0} \times_{\cD_{i,0}} \wh{\cD}_i \lra \wh{\cD}_i$ are formally 
etale morphisms lifting $D'_i \lra D_i$. So they are isomorphic. So we have 
$\Gamma(\wh{\cD}'_i,\cO_{\wh{\cD}'_i}) = 
\Gamma(\cD'_{i,0},\cO_{\cD'_{i,0}}) \wh{\otimes}_A 
A[[x_{m+1},\cdots,x_n]]$. Hence, to prove that the map 
\begin{align*}
& \Gamma(]D_i[_{\cD_i},\cO_{]D_i[_{\cD_i}}) 
= \Gamma(\wh{\cD}_{i,K}, \cO_{\wh{\cD}_{i,K}}) 
= \Gamma(\cD_{i,0,K} \times A^{n-m}_K[0,1),
\cO_{\cD_{i,0,K} \times A^{n-m}_K[0,1)}) \\ 
\lra & 
\Gamma(\cD'_{i,0,K} \times A^{n-m}_K[0,1),
\cO_{\cD'_{i,0,K} \times A^{n-m}_K[0,1)}) = 
\Gamma(\wh{\cD}'_{i,K}, \cO_{\wh{\cD}'_{i,K}}) = 
\Gamma(]D'_i[_{\cD'_i},\cO_{]D'_i[_{\cD'_i}})
\end{align*}
is injective, it suffices to prove that the map 
$\Gamma(\cD_{i,0,K},
\cO_{\cD_{i,0,K}}) \lra 
\Gamma(\cD'_{i,0,K}, \cO_{\cD'_{i,0,K}})$ is injective. 
Using the formal etaleness of $\cD'_{i,0} \lra \cD_{i,0}$, 
the denseness of the image of $D'_i$ in $D_i$ and the $p$-torsion freeness 
of $\Gamma(\cD_{i,0},\cO_{\cD_{i,0}})$ and 
$\Gamma(\cD'_{i,0},\cO_{\cD'_{i,0}})$, we see easily the injectivity of 
the map $\Gamma(\cD_{i,0,K},
\cO_{\cD_{i,0,K}}) \lra 
\Gamma(\cD'_{i,0,K}, \cO_{\cD'_{i,0,K}})$. So we are done. 
\end{pf} 

\begin{lem}\label{lem3} 
Let us assume given a diagram \eqref{lem1diag} and assume that 
$g$ is formally smooth and that the images of 
the morphisms $D'_i:=X' \times_{X}D_i \lra D_i$ induced by $f$ are 
dense for all $1 \leq i \leq r$. Then $(E,\nabla)$ has nilpotent 
residues if and only if $(E',\nabla')$ has nilpotent residues. 
$($Here $(E,\nabla), (E',\nabla')$ are as in Lemma \ref{lem1}.$)$
\end{lem} 

\begin{pf} 
We may replace $(\cX',M_{\cX'})$ by a Zariski open covering of it. 
So we may assume by \cite[claim in p.81]{shiho2} that there exists a 
strict formally etale morphism $(\cX'',M_{\cX''}) \lra (\cX,M_{\cX})$ 
such that $(X',M_{X'}) = (X,M_X) \times_{(\cX,M_{\cX})} (\cX'',M_{\cX''})$ 
holds. 
Let us put $(\cX''',M_{\cX'''}) := 
(\cX',M_{\cX'}) \times_{(\cX,M_{\cX})} (\cX'',M_{\cX''})$, and 
let us denote the restriction of $(E,\nabla)$ to $]X'[_{\cX''}, 
]X'[_{\cX'''}$ by $(E'',\nabla''), \allowbreak 
(E''',\nabla''')$, respectively. Then 
$(E,\nabla)$ has nilpotent residues if and only if so is 
$(E'',\nabla'')$ by Lemma \ref{lem1} and Lemma \ref{lem2}(2), and 
$(E'',\nabla'')$ has nilpotent residues if and only if 
so is $(E''',\nabla''')$ by Lemma \ref{lem1} and Lemma \ref{lem2}(1), 
and $(E''',\nabla''')$ has nilpotent residues if and only if 
so is $(E',\nabla')$ again by Lemma \ref{lem1} and Lemma \ref{lem2}(1). 
\end{pf} 

\begin{prop}\label{prop1.9}
Let $X$ be a smooth scheme over $k$, 
let $D$ be a normal crossing divisor on $X$ and denote the log structure 
associated to $D$ by $M_X$. Let us assume given a diagram 
\begin{equation*}
\begin{CD}
(X',M_{X'}) @>i>> (\cX',M_{\cX'}) \\ 
@VfVV @. \\
(X,M_X), @. 
\end{CD}
\end{equation*}
where $f$ is a strict etale morphism and $i$ is a strongly admissible 
closed immersion. 
Let $\cD' = \bigcup \cD'_i$ be relative the simple normal crossing 
divisor on $\cX'$ associated to $M_{\cX'}$ such that each $\cD'_i$ is 
a formally smooth divisor defined by $t_i \in \Gamma(\cX',\cO_{\cX'}) 
\,(1 \leq i \leq r)$ and put $D'_i := X' \times_X \cD_i$. Let 
$\cE$ be a locally free isocrystal on $(X/\cS)^{\log}_{\conv}$ and denote the 
log-$\nabla$-module on $]X'[_{\cX'}$ with respect to $t_1,\cdots,t_r$ 
induced by $\cE$ by $(E,\nabla)$. Then, if $\cE$ has nilpotent residues in 
the sense of Definition \ref{nilpres}, 
$(E,\nabla)$ has nilpotent residues and it is unipotent on 
\begin{equation}\label{wf}
]D'_I[_{\cX'} \,\cong\, ]D'_I[_{\cD_I} \times A^{|I|}_K[0,1) 
\end{equation}
for any $I \subseteq \{1,\cdots,r\}$ $($where $D'_I:=\bigcap_{i\in I}D'_i$, 
$\cD'_I := \bigcap_{i \in I}\cD'_i)$. 
\end{prop} 

\begin{pf} 
Let us take a strict etale covering $\coprod_{j} (X_j,M_{X_j}) \lra (X,M_X)$ 
and a strongly admissible closed immersions $(X_j,M_{X_j}) \hra 
(\cX_j,M_{\cX_j})$ such that the log-$\nabla$-module $(E_j,\nabla_j)$ 
on $]X_j[_{\cX_j}$ induced by $\cE$ has nilpotent residues. 
(There exists such morphisms by the definition of `having nilpotent 
residues' for isocrystals.) Let us put 
$(X'_j,M_{X'_j}) := (X',M_{X'}) \times_{(X,M_X)} (X_j,M_{X_j})$. 
By replacing $(X'_j,M_{X'_j})$ by Zariski open covering of it, we may 
assume that there exist strict formally etale morphisms 
$(\cY',M_{\cY'}) \lra (\cX',M_{\cX'})$ and $(\cY_j,M_{\cY_j}) \lra 
(\cX_j,M_{\cX_j})$ such that $(X'_j,M_{X'_j}) = 
(X',M_{X'}) \times_{(\cX',M_{\cX'})} (\cY',M_{\cY'}) = 
(X_j,M_{X_j}) \times_{(\cX_j,M_{\cX_j})} (\cY_j,M_{\cY_j})$ holds. 
Moreover, we may assume that the closed immersions 
$(X'_j,M_{X'_j}) \hra (\cY',M_{\cY'}), 
(X'_j,M_{X'_j}) \hra (\cY_j,M_{\cY_j})$ are strongly admissible closed 
immersion and that the relative simple normal crossing divisor 
$\cC', \cC_j$ on $\cY',\cY_j$ corresponding to the log structure 
$M_{\cY'},M_{\cY_j}$ has the same number of irreducible smooth components. 
Let us put $\cC':=\bigcup_{i=1}^s\cC'_i, 
\cC_j:=\bigcup_{i=1}^s\cC_{j,i}$. Then let us define $(\cY'_j,M_{\cY'_j})$ 
to be 
`the logarithmic product' of $(\cY',M_{\cY'})$ and $(\cY_j,M_{\cY_j})$ over 
$\cS$: 
That is, we define $\cY'_j$ by 
\begin{align*}
\cY'_j & := \text{(blow-up of $\cY' \times_{\cS} \cY_j$ along 
$\bigcup_{i=1}^r(\pr_1^{-1}(\cC'_i) \cap \pr_2^{-1}(\cC_{j,i}))$)} \\ 
& - 
\bigcup_{i=1}^s\text{(strict transforms of $\pr_1^{-1}(\cC'_i)$ and 
$\pr_2^{-1}(\cC_{j,i})$)}
\end{align*}
and define the log structure $M_{\cY'_j}$ to be the log structure 
associated to the inverse image of $\cC'$ (which is equal to the inverse 
image of $\cC_j$). Then we have the commutative diagram 
\begin{equation*}
\begin{CD}
(X',M_{X'}) @<<< \coprod_j(X'_j,M_{X'_j}) @>>> \coprod_j (X_j,M_{X_j}) \\ 
@V{\cap}VV @V{\cap}VV @V{\cap}VV \\ 
(\cX',M_{\cX'}) @<<< (\cY'_j,M_{\cY'_j}) @>>> (\cX_j,M_{\cX_j}), 
\end{CD}
\end{equation*}
where the vertical arrows are (the disjoint union of) 
strongly admissible closed immersions, 
the upper horizontal arrows are strict etale surjective morphisms and 
the lower horizontal arrow are strict formally smooth morphisms. \par 
Let $(E'_j,\nabla'_j)$ be the restriction of $(E,\nabla)$ to 
$]X'_j[_{\cY'_j}$, which is the same as the restriction of $(E_j,\nabla_j)$ 
to $]X'_j[_{\cY'_j}$. Then, by Lemma \ref{lem3} and the nilpotence of the 
residues of $(E_j,\nabla_j)$, $(E'_j,\nabla'_j)$ has nilpotent residues, and 
by using Lemma \ref{lem3} again, we see that $(E,\nabla)$ has nilpotent 
residues. So we have proved the former assertion of the proposition. \par 
Next we prove the latter assertion. 
We may work Zariski locally on $\cX'$ by \cite[3.3.5]{kedlayaI}. 
So we may assume 
the existence of a Cartesian diagram 
\begin{equation*}
\begin{CD}
X' @>{\subset}>> \cX' \\ 
@VVV @VVV \\ 
\Spec k[t_1,\cdots,t_m] @>{\subset}>> \Spf V\{t_1,\cdots,t_n\} 
\end{CD}
\end{equation*} 
for some $r \leq m\leq n$, where 
the right vertical arrow is formally etale, 
the lower horizontal arrow is defined by 
$t_{m+1}=\cdots =t_n=0$ and the closed immersion $\Spec k \hra \Spf V$ and 
each $\cD'_i \, (1 \leq i \leq r)$ is the zero locus of $t_i$. 
Let $\cX'_0$ be the zero locus $t_{m+1}=\cdots = t_n=0$ in $\cX'$, 
denote the completion of $\cX', \cX'_0$ along $D'_I$ by $\wh{\cX}', 
\wh{\cX}'_0$, respectively and denote the pull-back of $\cD'_I$ to 
$\wh{\cX}', \wh{\cX}'_0$ by $\wh{\cD}'_I, \wh{\cD}'_{0,I}$, 
respectively. \par 
For disjoint subsets $L,L'$ of $\{1,\cdots, n\}$, let us define 
$\Af^L_k := \Spec k[t_i\,(i \in L)], 
\fAf^{L,L'}_V \allowbreak 
:= \Spf V\{t_i\,(i \in L)\}[[t_i\,(i \in L')]]$. 
Then, if we put $I^c := \{1,\cdots, m\}-I, J := \{m \allowbreak 
+1,\cdots, n\}$, 
we have the following Cartesian diagrams 
\begin{equation}\label{af1}
\begin{CD}
D'_I @>{\subset}>> \wh{\cX}'_0 @>{\subset}>> \wh{\cX}' \\
@VVV @VVV @VVV \\ 
\Af^{I^c}_k @>{\subset}>> \fAf^{I^c,I}_V @>{\subset}>> 
\fAf^{I^c,I\cup J}_V, 
\end{CD}
\,\,\,\,\,\,
\begin{CD}
D'_I @>{\subset}>> \wh{\cD}'_{0,I} @>{\subset}>> \wh{\cD}'_I \\ 
@VVV @VVV @VVV \\  
\Af^{I^c}_k @>{\subset}>> \fAf^{I^c,\emptyset}_V @>{\subset}>> 
\fAf^{I^c,J}_V, 
\end{CD}
\end{equation}
where the vertical arrows are formally etale and the horizontal arrows are 
natural closed immersions. Let us define the morphisms 
\begin{equation*}
\begin{CD}
\fAf^{I^c,I}_V @<{s'_J}<< \fAf^{I^c,I\cup J}_V \\
@V{s_I}VV @V{s'_I}VV \\ 
\fAf^{I^c,\emptyset}_V @<{s_J}<< \fAf^{I^c,J}_V 
\end{CD}
\end{equation*}
as the morphism induced by the natural inclusion of rings of the form 
$V\{t_i\,(i \in L)\}[[t_i\,(i \in L')]] \hra 
V\{t_i\,(i \in L)\}[[t_i\,(i \in L'')]] \,(L' \subseteq L'')$, and let 
us put $s:=s_I \circ s'_J = s_J \circ s'_I$. Then, both 
$\wh{\cD}'_{0,I} \times_{\fAf^{I^c,\emptyset}_V,s} 
\fAf^{I^c,I\cup J}_V \lra \fAf^{I^c,I\cup J}_V, \wh{\cX}' \lra 
\fAf^{I^c,I\cup J}_V$ are formally etale morphism lifting the 
morphism $D'_I \lra \Af^{I^c}_k$. So there exists an isomorphism 
$\wh{\cD}'_{0,I} \times_{\fAf^{I^c,\emptyset}_V,s} 
\fAf^{I^c,I\cup J}_V \cong \wh{\cX}'$ over $\fAf^{I^c,I\cup J}_V$, and 
it induces the isomorphisms 
$\wh{\cD}'_{0,I} \times_{\fAf^{I^c,\emptyset}_V,s_I} 
\fAf^{I^c,I}_V \cong \wh{\cX}'_0$, 
$\wh{\cD}'_{0,I} \times_{\fAf^{I^c,\emptyset}_V,s_J} 
\fAf^{I^c,J}_V \cong \wh{\cD}'_I$. Thus we have the isomorphisms 
$$ 
\varphi_0: \wh{\cX}'_0 \os{=}{\lra} 
\wh{\cD}'_{0,I} \times_{\fAf^{I^c,\emptyset}_V,s_I} \fAf^{I^c,I}_V 
\os{=}{\lra}
\wh{\cD}'_{0,I} \times_V \fAf^{\emptyset,I}_V, 
$$ 
{\small{ $$
\varphi: \wh{\cX}' \os{=}{\lra} 
\wh{\cD}'_{0,I} \times_{\fAf^{I^c,\emptyset}_V,s} \fAf^{I^c,I\cup J}_V 
\os{=}{\lra}
(\wh{\cD}'_{0,I} \times_{\fAf^{I^c,\emptyset}_V,s_J} \fAf^{I^c,J}_V) 
\times_V \fAf^{\emptyset,I}_V 
\os{=}{\lra} \wh{\cD}'_I \times_V \fAf^{\emptyset,I}_V.  
$$ }}
Let us define the morphisms $\pi, \pi'$ by 
$$ 
\pi: \wh{\cD}'_I \os{=}{\lra} 
\wh{\cD}'_{0,I} \times_{\fAf^{I^c,\emptyset}_V,s_J} \fAf^{I^c,J}_V 
\os{\id \times s_J}{\lra} 
\wh{\cD}'_{0,I}, 
$$ 
$$ 
\pi': \wh{\cX}' \os{=}{\lra} 
\wh{\cD}'_{0,I} \times_{\fAf^{I^c,\emptyset}_V,s} \fAf^{I^c,I\cup J}_V 
\os{\id \times s'_J}{\lra} 
\wh{\cD}'_{0,I} \times_{\fAf^{I^c,\emptyset}_V,s_I} 
\fAf^{I^c,I}_V \os{=}{\lra} \wh{\cX}'_0. 
$$
Then we have the commutative diagram 
\begin{equation*}
\begin{CD}
\wh{\cX}' @>{\varphi}>> \wh{\cD}'_I \times_V \fAf^{\emptyset,I} \\ 
@V{\pi'}VV @V{\pi \times \id}VV \\ 
\wh{\cX}'_0 @>{\varphi_0}>> 
\wh{\cD}'_{I,0} \times_V \fAf^{\emptyset,I} 
\end{CD}
\end{equation*}
and it induces the commutative diagram of rigid spaces 
\begin{equation*}
\begin{CD}
]D'_I[_{\cX'} @>{\cong}>> ]D'_I[_{\cD'_I} \times A_K^{|I|}[0,1) \\ 
@V{\pi'_K}VV @V{\pi_K \times \id}VV \\ 
]D'_I[_{\cX'_0} @>{\cong}>> ]D'_I[_{\cD'_{0,I}} \times A_K^{|I|}[0,1). 
\end{CD}
\end{equation*}
Let $(E_0,\nabla_0)$ be the log-$\nabla$-module on $]X'[_{\cX'_0}$ 
with respect to $t_1,\cdots,t_r$ induced by $\cE$. Then it is unipotent 
on $]D'_I[_{\cD'_{0,I}} \times A_K^{|I|}[0,1)$ relative to 
$]D'_I[_{\cD'_{0,I}}$ by Proposition \ref{isocunip}. By pulling-back by 
$\pi_K \times \id$, we see that $(E,\nabla)$ is unipotent on 
$]D'_I[_{\cD'_I} \times A_K^{|I|}[0,1)$ relative to 
$]D'_I[_{\cD'_I}$. So we have proved the latter assertion of the 
proposition. 
\end{pf}

\begin{cor} 
Let $X$ be a smooth scheme over $k$, 
let $D$ be a normal crossing divisor on $X$, let $M_X$ be the log 
structure on $X$ associated to $D$ and 
let $\cE$ be a locally free isocrystal on 
$(X/\cS)^{\log}_{\conv}=((X,M_X)/\cS)_{\conv}$. Then$:$ 
\begin{enumerate}
\item 
When $D$ is a simple normal crossing divisor, $\cE$ has nilpotent 
residues in the sense in Definition \ref{nilpres} if and only if 
it has nilpotent residues in the sense in \cite[6.4.4]{kedlayaI}. 
\item 
The definition of `having nilpotent residues' is independent of the 
choice of an etale covering $\coprod X_j \lra X$ and a diagram 
like \eqref{liftable} chosen in Definition \ref{nilpres}. 
\end{enumerate}
\end{cor} 

Next we prove the functoriality of the notion of `having nilpotent residues': 

\begin{prop}\label{functnilp}
Let $X,X'$ be smooth schemes over $k$, 
let $D$ $($resp. $D')$ 
be a normal crossing divisor on $X$ $($resp. $X')$, 
let $M_X$ $($resp. $M_{X'})$ be the log 
structure on $X$ $($resp. $X')$ associated to $D$ $($resp. $D')$ 
and let us assume given a morphism 
$f:(X',M_{X'}) \lra (X,M_X)$. Then, for a locally free isocrystal $\cE$ on 
$(X/\cS)^{\log}_{\conv}$ having nilpotent residues, the pull-back 
$f^*\cE$ of $\cE$ to $(X'/\cS)^{\log}_{\conv}$ also has nilpotent 
residues. 
\end{prop} 

\begin{pf} 
We may work etale locally on $X'$. So 
we may assume that $D,D'$ are simple normal crossing divisors and that 
there exist strongly admissible closed immersions 
$(X,M_X) \hra (\cX,M_{\cX}), (X',M_{X'}) \hra (\cX',M_{\cX'})$ 
satisfying $(X,M_X) = (\cX,M_{\cX}) \times_{\cS} S$, 
$(X',M_{X'}) = (\cX',M_{\cX'}) \times_{\cS} S$. Then, locally on $X'$, 
we have a morphism $g:(\cX',M_{\cX'}) \lra (\cX,M_{\cX})$ lifting $f$. 
Let $(E,\nabla)$ be the log-$\nabla$-module on $]X[_{\cX}$ induced by 
$\cE$ and let $(E',\nabla')$ be the log-$\nabla$-module on 
$]X'[_{\cX'}$ induced by $f^*\cE$. Then $(E',\nabla')$ is nothing but 
the pull-back of $(E,\nabla)$ by $g$. We may assume that $E$ is free 
on $]X[_{\cX}$. \par 
Let $\cD = \bigcup_{i=1}^r\cD_i$ (resp. $\cD':=\bigcup_{i=1}^{r'}\cD'_i$) 
be the relative simple normal crossing 
divisor on $\cX$ (resp. $\cX'$) corresponding to $M_{\cX}$ (resp. 
$M_{\cX'}$) and let $t_i \in \Gamma(\cX,\cO_{\cX})\,(1 \leq i \leq r)$ 
be an element defining $\cD_i$. 
For fixed $i \,(1 \leq i \leq r)$, we may replace $\cD',\cX$ by 
$\cD'_i - \bigcup_{j\not=i}\cD'_j, 
\cX - \bigcup_{j\not=i}\cD'_j$ to prove the nilpotence of the residue of 
$(E',\nabla')$ along $\cD'_i$. (See the proof of Lemma \ref{lem2}.) 
So we may assume that $\cD'$ is a smooth divisor. Let 
$s \in \Gamma(\cX',\cO_{\cX'})$ be an element defining $\cD'$. 
Since $g$ is a morphism of log schemes, $g^*t_i$ has the form 
$u_is^{n_i}$ for some $u_i \in \Gamma(\cX',\cO_{\cX'}^{\times})$ and 
$n_1 \in \N$. \par 
Let us denote the residue of $(E,\nabla)$ along $\cD_i$ by $\res_i$ and 
let us put $I:=\{i\,\vert\,n_i>0\}$. For $i \in I$, 
denote the map $\End(E|_{]D_i[_{\cD_i}}) \lra \End(E'|_{]\cD'[_{\cD'}})$
(induced by $f$ and $g$) by $f_i^*$. If we put 
$D_I:=\bigcap_{i\in I}D_i, \cD_I:=\bigcap_{i\in I}\cD_i$, 
the map $f_i^*$'s for $i \in I$ factors as 
$$ \End(E|_{]D_i[_{\cD_i}}) \lra \End(E|_{]\cD_I[_{\cD_I}}) \lra 
\End(E'|_{]\cD'[_{\cD'}}), $$ 
where the first map is the restriction and the second map is the one 
induced by $f$ and $g$ which is independent of $i \in I$. 
Since $(E,\nabla)$ is integrable, we have $\res_i \circ \res_j = 
\res_j \circ \res_i$ in $\End(E|_{]\cD_I[_{\cD_I}})$. So we have 
$f^*_i(\res_i) \circ f^*_j(\res_j) = f^*_j(\res_j) \circ f^*_i(\res_i)$. 
Simple computation implies that the residue $\res'$ of $(E',\nabla')$ along 
$]D'[_{\cD'}$ is given by 
$e \mapsto \sum_{i\in I}n_if_i^*(\res_i)(e)$. Then, if we take a positive 
integer $N$ satisfying 
$\res_i^N=0$ for $1 \leq i \leq r$, we have 
$$ (\res')^{N(|I|-1)+1} = 
\sum_{\scriptstyle \sum_{i\in I}m_i=N(|I|-1)+1 \atop 
\scriptstyle m_i \geq 0}
\prod_{i \in I}n_i^{m_i}f_i^*(\res_i^{m_i}) = 0. $$
So we are done. 
\end{pf} 

Next we prove a relation between Frobenius structure and the nilpotence 
of residues. (This is already remarked, for example, in 
\cite[7.1.4]{kedlayaI}.)

\begin{prop}\label{frobnilp}
Let $X$ be a smooth scheme over $k$, 
let $D$ be a normal crossing divisor on $X$, let $M_X$ be the log 
structure on $X$ associated to $D$ and 
let $(\cE,\alpha)$ be an $F$-isocrystal on 
$(X/\cS)^{\log}_{\conv}=((X,M_X)/\cS)_{\conv}$. 
Then $\cE$ is locally free and it has nilpotent residues. 
\end{prop} 

\begin{pf}
By \cite[2.4.3]{shiho2}, $\cE$ is locally free. To show that 
$\cE$ has nilpotent residues, we may assume that 
$D$ is a simple normal crossing divisor, and 
we may assume the existence of a strongly admissible closed immersion 
$(X,M_X) \hra (\cX,M_{\cX})$ with $(X,M_X) = (\cX,M_{\cX}) \times_{\cS} S$ 
($M_{\cX}$ is induced by a relative simple normal crossing divisor 
$\cD = \bigcup_{i=1}^r\cD_i$ and $\cD_i$ is defined by $t_i$) endowed with 
a lift $F:(\cX,M_{\cX}) \lra (\cX,M_{\cX})$ of Frobenius on $(X,M_X)$ 
compatible with $\sigma$ satisfying $F^*t_i=t_i^p \,(1 \leq i \leq r)$. 
We may assume that the log-$\nabla$-module $(E,\nabla)$ on 
$]X[_{\cX}$ induced by $\cE$ is free. Moreover, by the same argument as 
the proof of Proposition \ref{functnilp}, we may assume that $\cD$ is a 
smooth divisor (hence $r=1$). Let us denote the residue of 
$(E,\nabla)$ along $]D[_{\cD}$ by $\res$. Then, on 
$$ ]D[_{\cX} = ]D[_{\cD} \times A_K^1[0,1), $$ 
we have, by \cite[1.5.3]{bc}, a polynomial $P(x) \in K[x]$ such that 
$P(\res)=0$ holds. Take $P(x)$ to be minimal and monic. 
By considering the image of $\res$ in $\End(F^*E|_{]D[_{\cD}})$, 
we see that $P^{\sigma}(p^a \res')=0$ holds, where $\res'$ is the residue of 
$F^*(E,\nabla)$ along $]D[_{\cD}$. Since $F^*(E,\nabla)$ is isomorphic to 
$(E,\nabla)$ (via the Frobenius structure $\alpha$), we have
$P^{\sigma}(p^a \res)=0$. So $P^{\sigma}(p^ax)$ is equal to $P(x)$ up to a 
multiplication by a non-zero constant, 
and so we have $P(x)=x^m$ for some $m$. Hence $\res$ is nilpotent. 
\end{pf} 

Finally in this section, we recall the statement of the semi-stable 
reduction conjecture for overconvergent $F$-isocrystals, which is shown 
recently by Kedlaya (\cite{kedlayaI}, \cite{kedlayaII}, \cite{kedlayaIII}, 
\cite{kedlayaIV}). 

\begin{thm}[Semistable reduction conjecture, Kedlaya]\label{sconj}
Let $(X,\ol{X})$ be a pair over $k$ and let $\cE$ be an overconvergent 
$F$-isocrystal on $(X,\ol{X})/\cS_K$. Then there exists a strict morphism 
of pairs $\varphi:(X',\ol{X}') \lra (X,\ol{X})$ such that 
$\varphi: \ol{X}' \lra \ol{X}$ is proper surjective and generically etale, 
$D':=\ol{X}'-X'$ is a simple normal crossing divisor on $\ol{X}'$, and 
if we denote the log structure on $\ol{X}'$ associated to $D'$ by 
$M_{\ol{X}'}$, there exists an $F$-isocrystal $\cF$ on 
$(\ol{X}'/\cS)^{\log}_{\conv} = ((\ol{X}',M_{\ol{X}'})/\cS)_{\conv}$ 
satisfying $\varphi^*\cE = j^{\d}\cF$, where $j^{\d}$ is the functor of 
restriction to overconvergent isocrystals on $(X',\ol{X}')/\cS_K$. 
\end{thm} 

As pointed out in \cite[7.1.4]{kedlayaI}, the conjecture remains true if 
we require $\cF$ to be a locally free isocrystal (without Frobenius 
structure) having nilpotent residues: Indeed, by Proposition \ref{frobnilp}, 
any $F$-isocrystal on $(\ol{X}'/\cS)^{\log}_{\conv}$ is locally free 
and has nilpotent residues, 
and if we have a locally free isocrystal $\cF$ on 
$(\ol{X}'/\cS)^{\log}_{\conv}$ 
having nilpotent residues satisfying $\varphi^*\cE = j^{\d}\cF$, 
we can extend the Frobenius structure on $\varphi^*\cE$ to that on 
$\cF$ by using \cite[6.4.5]{kedlayaI}. \par 
To state the main result in this paper without using Frobenius structure, 
we introduce a terminology of `potential semi-stability' 
of overconvergent isocrystals: 

\begin{defn} 
Let $(X,\ol{X})$ be a pair over $k$ and let $\cE$ be an overconvergent 
isocrystal on $(X,\ol{X})/\cS_K$. $\cE$ is called potentially semi-stable 
if there exists a strict morphism 
of pairs $\varphi:(X',\ol{X}') \lra (X,\ol{X})$ such that 
$\varphi: \ol{X}' \lra \ol{X}$ is proper surjective and generically etale, 
$D':=\ol{X}'-X'$ is a simple normal crossing divisor on $\ol{X}'$, and 
if we denote the log structure on $\ol{X}'$ associated to $D'$ by 
$M_{\ol{X}'}$, there exists a locally free isocrystal $\cF$ on 
$(\ol{X}'/\cS)^{\log}_{\conv} = ((\ol{X},M_{\ol{X}})/\cS)_{\conv}$ 
having nilpotent residues satisfying $\varphi^*\cE = j^{\d}\cF$, 
where $j^{\d}$ is the functor of 
restriction to overconvergent isocrystals on $(X',\ol{X}')/\cS_K$.
\end{defn}


\section{Overconvergence in log smooth case} 

In this section, we prove the 
overconvergence of relative rigid 
cohomology for (not necessarily strict) 
morphisms of certain pairs which admit `nice' 
log structures. \par 
Let us assume given a diagram 
\begin{equation}\label{diag-S2}
(\ol{X},M_{\ol{X}}) \os{f}{\lra} (\ol{Y},M_{\ol{Y}}) \os{\iota}{\hra} 
(\cY,M_{\cY}), 
\end{equation}
where $f$ is proper log smooth, $\iota$ is an exact closed immersion, 
$\ol{X}, \ol{Y}$ are smooth over $S$, 
$M_{\ol{X}}$ (resp. $M_{\ol{Y}}$) is the log structure 
induced by some normal crossing divisor (resp. simple normal crossing 
divisor) and $(\cY,M_{\cY})$ is of Zariski type and 
formally log smooth over $\cS$. We denote the 
normal crossing divisor corresponding to $M_{\ol{X}}$ by $D$. 
We assume moreover 
that $D$ has the decomposition $D=D^{\ho} \cup D^{\ve}$ into sub normal 
crossing divisors $D^{\ho}, D^{\ve}$ satisfying the following condition 
$(*)$: \\
\quad \\
$(*)$: \,\,\, If we denote the log structure on $\ol{X}$ associated to 
$D^{\ve}$ by 
$M^{\ve}_{\ol{X}}$ and define 
$D_{[i]} \,(i \geq 0)$ by 
\begin{align*}
D_{[0]} & :=\ol{X}, \qquad D_{[1]} := \text{the normalization of $D^{\ho}$}, \\
D_{[i]} & := \text{$i$-fold fiber product of $D_{[1]}$ over $\ol{X}$}, 
\end{align*}
$(D_{[i]},M_{\ol{X}}^{\ve}|_{D_{[i]}})$ 
is log smooth over $(\ol{Y},M_{\ol{Y}})$ for any $i \geq 0$. \\
\quad \\
Then we have the following comparison 
theorem between the 
log analytic cohomology of $(\ol{X},M_{\ol{X}})/(\cY,M_{\cY})$ 
and the log rigid cohomology of $((\ol{X}-D^{\ho},
M_{\ol{X}}|_{\ol{X}-D^{\ho}}), (\ol{X}, \allowbreak 
M_{\ol{X}}))/\allowbreak (\cY,\allowbreak M_{\cY})$. 

\begin{thm}\label{thm2.1}
In the situation above, let $\cE$ be a locally free isocrystal on 
$(\ol{X}/\cS)^{\log}_{\conv}\allowbreak = 
\allowbreak ((\ol{X},M_{\ol{X}})/\cS)_{\conv}$ 
having nilpotent residues and denote the restriction of $\cE$ to 
$I_{\conv}\allowbreak ((\ol{X}/\cY)^{\log})$ by the same symbol. 
Assume moreover the following condition $(\star):$ \\
\quad \\
$(\star):$ \,\,\, $R^qf_{\ol{X}/\cY,\an *}\cE$ is a coherent 
$\cO_{]\ol{Y}[_{\cY}}$-module for any $q \in \N$. \\
\quad \\
Then, if we denote the restriction functor 
$$ I_{\conv}((\ol{X}/\cY)^{\log}) \lra I^{\d}(((\ol{X}-D^{\ho},
\ol{X})/\cY)^{\log}) $$ 
by $j^{\d}$, we have the isomorphism 
$$ 
R^qf_{\ol{X}/\cY,\an *}\cE \os{=}{\lra} 
R^qf_{(\ol{X}-D^{\ho},\ol{X})/\cY, \logrig *}j^{\d}\cE \,\,\,\, (q \geq 0). 
$$ 
\end{thm} 

\begin{pf} 
First we prove the existence of a `nice' chart of the morphism 
$f$. (We take a chart Zariski locally on $\ol{Y}$ and etale locally on 
$\ol{X}$.) 
Let us fix a chart $\ti{\beta}:Q \ra M_{\cY}$ of $M_{\cY}$ and 
restrict it to the chart $\beta:Q \ra M_{\ol{Y}}$ of $M_{\ol{Y}}$. 
(Note that we can take this chart Zariski locally on $\cY$.) 
Let us consider etale locally and assume that $D^{\ho}:=
\bigcup_{i=1}^rD_i$ is a simple normal crossing divisor and that $D_i$ 
is defined by $t_i \in \Gamma(\ol{X},\cO_{\ol{X}})$. 
Put $D_o := \bigcap_{i=1}^rD_i$, 
$M_{D_o} := M^{\ve}_{X}|_{D_o}$ and $\cI := \Ker(\cO_X \ra \cO_{D_o})$. 
Then, since $(D_o,M_{D_o})$ is strict etale over $(D_{[r]},M_{D_{[r]}})$, 
it is log smooth over $(\ol{Y},M_{\ol{Y}})$. So we have the following 
locally split exact sequence
$$ 0 \lra \cI/\cI^2 \lra 
\omega^1_{(\ol{X},M^{\ve}_{\ol{X}})/(\ol{Y},M_{\ol{Y}})}|_{D_o} 
\lra \omega^1_{D_o/\ol{Y}} \lra 0 $$
by \cite[2.1.3]{nakkshiho}. Since $\cI/\cI^2$ is freely 
generated by $t_1,\cdots, t_r$ locally, we see that there exists a basis of 
$\omega^1_{(\ol{X},M^{\ve}_{\ol{X}})/(\ol{Y},M_{\ol{Y}})}$ locally, 
consisting of a lift of a basis 
of $\omega^1_{D_o/\ol{Y}}$ of the form $\{\dlog m_i\}_{i=1}^{d}$ and 
$\{dt_i\}_{i=1}^r$. Using this and argueing as in \cite[3.13]{kkato1} 
(see also \cite[2.1.4, 2.1.6]{nakkshiho}), 
one can see that there exists a chart 
$(P \os{\alpha}{\ra} M_{\ol{X}}^{\ve}, Q \os{\beta}{\ra} M_{\ol{Y}}, 
Q \os{\gamma}{\ra} P)$ 
of the morphism $(\ol{X},M_{\ol{X}}^{\ve}) \lra (\ol{Y},M_{\ol{Y}})$ 
extending the given chart $Q \os{\beta}{\ra} M_{\ol{Y}}$ such that 
$\gamma$ is injective, $|\Coker(\gamma)_{\tor}|$ is prime to $p$ and that 
the morphism 
$(\ol{X},M_{\ol{X}}^{\ve}) \lra (\ol{Y},M_{\ol{Y}})$ factors as 
$$ (\ol{X},M_{\ol{X}}^{\ve}) \lra 
(\ol{Y},M_{\ol{Y}}) 
\times_{(\Spec k[Q], M_{Q})} 
(\Spec k[P \oplus \N^r], M_{P}|_{\Spec k[P \oplus \N^r]}) \lra 
(\ol{Y},M_{\ol{Y}}), $$
where the first map, which is induced by 
$\alpha$ and $\N \ni e_i \mapsto t_i$, is 
strict etale and the second map is the one induced by $\beta$ and the 
composite 
$Q \os{\gamma}{\ra} P \hra P \oplus \N^r$ (which we denote by $\ti{\gamma}$). 
(For a fine monoid $N$ and a ring $R$ or an adic topological ring $R$, 
$M_{N}$ denotes the log structure on $\Spec R[N]$ or $\Spf R\{N\}$ associated 
to the pre-log structure $N \lra R[N]$ or $N \lra R\{N\}$.) 
By adding the log structure $M^{\ho}_{\ol{X}}$ associated to $D^{\ho}$, 
the above factorization induces the following factorization of the 
morphism $f$ 
\begin{equation}\label{factor}
(\ol{X},M_{\ol{X}}) \lra 
(\ol{Y},M_{\ol{Y}}) 
\times_{(\Spec k[Q], M_{Q})} 
(\Spec k[P \oplus \N^r], M_{P \oplus \N^r}) \lra 
(\ol{Y},M_{\ol{Y}}),
\end{equation}
where the first map is strict etale and the second map is induced by 
$\beta$ and $\ti{\gamma}$. So we obtain a chart 
$(P \oplus \N^r \ra M_{\ol{X}}, Q \os{\beta}{\ra} M_{\ol{Y}}, 
Q \os{\ti{\gamma}}{\ra} P \oplus \N^r)$ of $f$. \par 
Using this chart, we can form the 
Cartesian diagram 
{\small{
\begin{equation}\label{lift1}
\begin{CD}
(\ol{X},M_{\ol{X}}) @>>> 
(\ol{Y},M_{\ol{Y}}) \times_{(\Spec k[Q],M_Q)} 
(\Spec k[P \oplus \N^r],M_{P \oplus \N^r}) @>>> 
(\ol{Y},M_{\ol{Y}}) \\ 
@VVV @VVV @V{\iota}VV \\ 
(\cX,M_{\cX}) @>>> 
(\cY,M_{\cY}) \times_{(\Spf V\{Q\}, M_Q)} 
(\Spf V\{P \oplus \N^r\},M_{P \oplus \N^r}) 
@>>> (\cY,M_{\cY}),
\end{CD}
\end{equation}}}
where the top horizontal line 
is as in \eqref{factor}, the vertical arrows are 
exact closed immersions, the lower right arrow 
is defined by $\ti{\beta}$ and $\ti{\gamma}$, and the lower left arrow is 
a strict formally etale morphism which makes the left square Cartesian. 
(The existence of such a morphism follows from \cite[claim in p.81]{shiho2}). 
Since $\gamma$ is injective and $|(\Coker \ti{\gamma})_{\tor}|$ is 
prime to $p$, 
the lower right arrow in \eqref{lift1} is formally log smooth. So 
the morphism $(\cX,M_{\cX}) \lra (\cY,M_{\cY})$ is a formally log smooth 
lift of $(\ol{X},M_{\ol{X}}) \lra (\ol{Y},M_{\ol{Y}})$. 
So we have shown the existence of a formally log smooth lift of 
$f$ etale locally on $\ol{X}$. \par 
So we can take a strict etale surjective map 
$(\ol{X}^{(0)},M_{\ol{X}^{(0)}}) \lra (\ol{X},M_{\ol{X}})$ and an exact 
closed immersion 
$(\ol{X}^{(0)},M_{\ol{X}^{(0)}}) \hra (\cX^{(0)},M_{\cX^{(0)}})$ over 
$(\cY,M_{\cY})$ such that $(\cX^{(0)},M_{\cX^{(0)}})$ is of Zariski type, 
formally log smooth 
over $(\cY,M_{\cY})$ and that $(\ol{X}^{(0)},\allowbreak M_{\ol{X}^{(0)}}) = 
(\ol{Y},M_{\ol{Y}}) \times_{(\cY,M_{\cY})} (\cX^{(0)},M_{\cX^{(0)}})$ holds. 
For $m \in \N$, let $(\ol{X}^{(m)},M_{\ol{X}^{(m)}})$ 
(resp. $(\cX^{(m)}, \allowbreak M_{\cX^{(m)}}$)) be the $(m+1)$-fold fiber product of 
$(\ol{X}^{(0)},M_{\ol{X}^{(0)}})$ (resp. $(\cX^{(0)},M_{\cX^{(0)}})$) over 
$(\ol{X},M_{\ol{X}})$ (resp. $(\cY,M_{\cY})$). Then we have the embedding 
system 
\begin{equation}\label{embsys1}
(\ol{X},M_{\ol{X}}) \lla (\ol{X}^{(\b)},M_{\ol{X}^{(\b)}}) \hra 
(\cX^{(\b)},M_{\cX^{(\b)}}). 
\end{equation} 

Let us take a diagram 
\begin{equation}\label{homeo1}
\begin{CD}
(\ol{Y}',M_{\ol{Y}'}) @>{\iota'}>> (\cY',M_{\cY'}) \\ 
@VVV @V{\varphi}VV \\ 
(\ol{Y},M_{\ol{Y}}) @>{\iota}>> (\cY,M_{\cY}) 
\end{CD}
\end{equation} 
such that $\iota'$ is a homeomorphic exact closed immersion and that 
the induced morphism of rigid analytic spaces 
$\cY'_K = 
]\ol{Y}'[_{\cY'} \lra ]\ol{Y}[_{\cY}$ is an admissible open immersion. 
Let us denote the pull-back of the diagram \eqref{embsys1} by 
the morphism $\varphi$ by 
\begin{equation}\label{embsys2}
(\ol{X}',M_{\ol{X}'}) \lla ({\ol{X}'}^{(\b)},M_{{\ol{X}'}^{(\b)}}) \hra 
({\cX'}^{(\b)},M_{{\cX'}^{(\b)}}) 
\end{equation} 
and put ${D'}^{\ho}:= D^{\ho} \times_{\ol{X}} \ol{X}'$. 
Let us denote the morphism 
$({\cX'}^{(\b)},M_{{\cX'}^{(\b)}}) \lra (\cY',M_{\cY'})$ by $g$, 
the morphism $]{\ol{X}'}^{(\b)}[^{\log}_{{\cX'}^{(\b)}} \lra {\ol{X}'}^{(\b)}$ 
by $\sp^{(\b)}$ and the morphism $]\ol{Y}'[_{\cY'} \lra \cY'$ by $\sp$. 
We prove the following claim: \\
\quad \\
{\bf claim.} The morphism 
\begin{align*}
R\sp_*Rf_{\ol{X}'/\cY',\an *}\cE & = 
Rg_*R\sp^{(\b)}_*\DR(]{\ol{X}'}^{(\b)}[^{\log}_{{\cX'}^{(\b)}}/\cY'_K,\cE) \\ 
& \lra 
Rg_*R\sp^{(\b)}_*\DR^{\d}(]{\ol{X}'}^{(\b)}[^{\log}_{{\cX'}^{(\b)}}/\cY'_K,
j^{\d}\cE) \\ & = 
R\sp_*Rf_{(\ol{X}'-{D'}^{\ho},\ol{X}')/\cY', \logrig *}j^{\d}\cE
\end{align*} 
is a quasi-isomorphism. \\
\quad \\
To prove the claim, it suffices to prove that the map 
$$R\sp^{(m)}_*\DR(]{\ol{X}'}^{(m)}[^{\log}_{{\cX'}^{(m)}}/\cY'_K,\cE)
\lra 
R\sp^{(m)}_*\DR^{\d}(]{\ol{X}'}^{(m)}[^{\log}_{{\cX'}^{(m)}}/\cY'_K, 
j^{\d}\cE) $$ 
is a quasi-isomorphism. 
In the following (until the end of the proof of the 
claim), we omit to write the superscript ${}^{(m)}$. 
(So we prove that 
\begin{equation}\label{an-logrig}
R\sp_*\DR(]{\ol{X}'}[^{\log}_{{\cX'}}/\cY'_K,\cE)
\lra 
R\sp_*\DR^{\d}(]{\ol{X}'}[^{\log}_{{\cX'}}/\cY'_K, 
j^{\d}\cE) 
\end{equation}
is a quasi-isomorphism.) \par 
To prove that \eqref{an-logrig} is a quasi-isomorphism, 
we may work Zariski locally on 
${\cX'}$, and by \cite[4.5]{shiho3} and \cite[claim in 4.7]{shiho4}, 
the both hand sides of \eqref{an-logrig} are unchanged if 
we replace the closed immersion $(\ol{X}', M_{\ol{X}'}) 
\hra (\cX',M_{\cX'})$ by another closed immersion 
into a $p$-adic fine log formal $\cS$-scheme which is formally log smooth
 over $(\cY',M_{\cY'})$. So we may assume the existence of the Cartesian 
 diagram 
\begin{equation*}
\begin{CD}
(\ol{X}',M_{\ol{X}'}) @>>> (\cX',M_{\cX'}) \\
@VVV @VVV \\ 
(\ol{X}^{(0)},M_{\ol{X}^{(0)}}) @>>> (\cX^{(0)},M_{\cX^{(0)}}) 
\end{CD}
\end{equation*}
such that the right vertical arrow is written as a composite morphism 
$$ (\cX',M_{\cX'}) \lra ({\cX'}^{(0)},M_{{\cX'}^{(0)}}) \lra (\cX^{(0)}, 
M_{\cX^{(0)}}), $$
where the first morphism is strict formally etale and the second morphism is 
the canonical morphism. Since  we may work Zariski locally on 
${\cX'}$, we may shrink $\cX^{(0)}$ in order that the exact closed immersion 
$(\ol{X}^{(0)},M_{\ol{X}^{(0)}}) \hra (\cX^{(0)},M_{\cX^{(0)}})$ is a 
strongly admissible closed immersion. Let 
$\cD^{(0)} = \bigcup_{i=1}^r\cD^{(0)}_i$ 
be the relative simple normal crossing divisor 
corresponding to $M_{\cX^{(0)}}$ and for $I \subseteq \{1,\cdots,r\}$, 
let us put $\cD^{(0)}_I := \bigcap_{i \in I}\cD^{(0)}_I$. Let 
us denote the pull-back of $\cD^{(0)}_I$ to $\ol{X}^{(0)}, \cX', \ol{X}'$ by 
$D^{(0)}_I, \cD'_I, D'_I$, respectively. \par 
Note that 
$R^q\sp_*\DR(]\ol{X}'[^{\log}_{{\cX'}}/\cY'_K,\cE)$ is the sheaf associated to 
the presheaf 
$$ \cX' \supseteq \cU \mapsto H^q(\cU_K, 
\DR(]\ol{X}'[^{\log}_{{\cX'}}/\cY'_K,\cE)) $$ 
and that 
$R^q\sp_*\DR^{\d}(]\ol{X}'[^{\log}_{{\cX'}}/\cY'_K, 
j^{\d}\cE)$ is the sheaf associated to the presheaf 
$$ \cX' \supseteq \cU \mapsto H^q(\cU_K, 
\DR^{\d}(]\ol{X}'[^{\log}_{{\cX'}}/\cY'_K,j^{\d}\cE)). 
$$ 
Since we may replace $\cX'$ by $\cU$, we see that it suffices to prove 
that the map 
$$ 
H^q(\cX'_K, 
\DR(]\ol{X}'[^{\log}_{{\cX'}}/\cY'_K,\cE)) \lra 
H^q(\cX'_K, 
\DR^{\d}(]\ol{X}'[^{\log}_{{\cX'}}/\cY'_K,j^{\d}\cE))
$$ 
is an isomorphism. If we put the admissible open immersion 
$\cX'_K-[{D'}^{\ho}]_{\cX',\lam} \hra \cX'_K$ by $j_{\lam}$, 
we have 
\begin{align*}
H^q(\cX'_K, 
\DR^{\d}(]\ol{X}'[^{\log}_{{\cX'}}/\cY'_K,j^{\d}\cE))
& = 
H^q(\cX'_K, 
\varinjlim_{\lam\to 1}
j_{\lam,*}j_{\lam}^{*}
\DR(]\ol{X}'[^{\log}_{{\cX'}}/\cY'_K,\cE)) \\ 
& = 
\varinjlim_{\lam\to 1}
H^q(\cX'_K, 
j_{\lam,*}j_{\lam}^{*}
\DR(]\ol{X}'[^{\log}_{{\cX'}}/\cY'_K,\cE)). 
\end{align*}
So it suffices to prove that the map 
$$ H^q(\cX'_K, 
\DR(]\ol{X}'[^{\log}_{{\cX'}}/\cY'_K,\cE)) \lra 
H^q(\cX'_K, 
j_{\lam,*}j_{\lam}^{*}
\DR(]\ol{X}'[^{\log}_{{\cX'}}/\cY'_K,\cE)) $$ 
is an isomorphism for $\lam \in (0,1)$. 
In the following, we denote $\DR(]\ol{X}'[^{\log}_{{\cX'}}/\cY'_K,\cE)$ 
simply by $\DR$. 
\par 
Let $(E,\nabla)$ be the log-$\nabla$-module on 
$]\ol{X}'[^{\log}_{\cX'} = ]\ol{X}'[_{\cX'} = \cX'_K$ 
induced by $\cE$, and let $(\ol{E},\ol{\nabla})$ be the  
log-$\nabla$-module on $\cX'_K$ relative to $\cY'_K$ 
induced by $(E,\nabla)$. Then $\DR$ is nothing but the de Rham complex 
associated to $(\ol{E},\ol{\nabla})$. On the other hand, let 
$(E^{(0)},\nabla^{(0)})$ be the log-$\nabla$-module on 
$]\ol{X}^{(0)}[_{\cX^{(0)}}$ induced by $\cE$. 
Then, by Proposition \ref{prop1.9}, 
$(E^{(0)},\nabla^{(0)})$ is unipotent on 
\begin{equation*}
]D^{(0)}_I[_{\cX^{(0)}} \,\cong\, 
]D^{(0)}_I[_{\cD_I^{(0)}} \times A^{|I|}_K[0,1). 
\end{equation*} 
Then $(E,\nabla)$ is unipotent on 
\begin{equation}\label{wf3}
]D'_I[_{\cX'} \,\cong\, 
]D'_I[_{\cD'_I} \times A^{|I|}_K[0,1). 
\end{equation} 

Next we consider a certain admissible covering of $\cX'_K$ 
introduced by Baldassarri and Chiarellotto (\cite[4.2]{bc2}). 
(See also \cite[2.4]{shiho2}.) Let us fix $\eta \in (\lam,1) \cap p^{\Q}$ and 
for $I \subseteq \{1,\cdots,r\}$, put 
$$ P_{I,\eta}:= \{x \in \cX'_K \,\vert\, 
\vert t_i(x) \vert < 1 \,(i \in I), 
\vert t_i(x) \vert \geq \eta \, (i \notin I) \}, $$
$$ V_{I,\eta}:= \{x \in \cX'_K \,\vert\, 
\vert t_i(x) \vert = 0 \,(i \in I), 
\vert t_i(x) \vert \geq \eta \,(i \notin I) \}. $$
Then we have the admissible open covering $\cX'_K = \bigcup_{I} P_{I,\eta}$. 
For $m \in \N$, let $\cI_m$ be the set 
$\{(I_0, \cdots, I_m) \,\vert\, I_j \subseteq \{1,\cdots,r\}\}.$ 
For $\I:=(I_0, \cdots, I_m) \in \cI_m$, let us put 
$\II := \bigcap_{j=0}^mI_j$ and $P_{\I,\eta}:=\bigcap_{j=0}^m P_{I_j,\eta}$. 
Then, since we have the spectral 
sequence 
$$ E^{s,t}_1 = \bigoplus_{\I \in \cI_s} H^t(P_{\I,\eta}, -) 
\,\Longrightarrow\, H^{s+t}(\cX'_K,-), $$
it suffices to prove the map 
\begin{equation}\label{qis1}
H^q(P_{\I,\eta}, 
\DR) \lra 
H^q(P_{\I,\eta}, 
j_{\lam,*}j_{\lam}^{*}
\DR)
\end{equation} 
is an isomorphism. Note that we have 
$$ P_{\I,\eta}:= \{x \in \cX'_K \,\vert\, 
\vert t_i(x) \vert < 1 \,(i \in \bigcup_{j=0}^m I_j), \,
\vert t_i(x) \vert \geq \eta \,(i \notin \II) \}. $$
Let us define $V_{\I,\eta}$ by 
\begin{align*}
V_{\I,\eta}:= \{x \in \cX'_K \,\vert\, 
\vert t_i(x) \vert = 0 \,(i \in \II),\, &
\vert t_i(x) \vert \geq \eta \, (i \notin \II), \\ & 
\vert t_i(x) \vert < 1 \,(i \in \bigcup_{j=0}^m I_j- \II) \}. 
\end{align*}
Then $P_{\I,\eta}$ is a quasi-Stein admissible open set of 
$]D'_{\II}[_{\cX'}$, 
$V_{\I,\eta}$ is a quasi-Stein admissible open set of 
$]D'_{\II}[_{\cD'_{\II}}$ and the isomorphism 
\eqref{wf3} ($I$ replaced by $\II$) induces the isomorphism 
$$ P_{\I,\eta} \cong V_{\I,\eta} \times A^{|\II|}_K[0,1). $$
Let us take an admissible open covering $V_{\I,\eta} = 
\bigcup_{j=1}^{\infty} V_j$ 
by increasing affinoid admissible open sets. Then we have 
\begin{align*}
H^q(P_{\I,\eta},\DR) &
= H^q(R\dvpl_j R\Gamma(V_j \times A^{|\II|}_K[0,1),\DR)), \\
H^q(P_{\I,\eta},j_{\lam,*}j^*_{\lam}\DR) & = 
H^q(R\dvpl_j R\Gamma(V_j \times A^{|\II|}_K[0,1),j_{\lam,*}j_{\lam}^*\DR)) 
\end{align*}
and we are reduced to showing that the map 
\begin{equation}\label{asia}
H^q(V_j \times A^{|\II|}_K[0,1),\DR) \lra 
H^q(V_j \times A^{|\II|}_K[0,1),j_{\lam,*}j_{\lam}^*\DR)
\end{equation} 
is an isomorphism for each $q \geq 0$. Then, since 
$V_j \times A^{|\II|}_K[0,1)$ and the map 
$j_{\lam}$ are quasi-Stein, the map \eqref{asia} is identical with the map 
\begin{equation}\label{qis2}
H^q(\Gamma(V_j \times A^{|\II|}_K[0,1),\DR)) \lra 
H^q(\Gamma(V_j \times A^{|\II|}_K(\lam,1),\DR)). 
\end{equation} 
So it suffices to prove that the map \eqref{qis2} 
is an isomorphism. Let $(\ol{\ol{E}},\ol{\ol{\nabla}})$ be the  
log-$\nabla$-module on $V_j \times A^{|\II|}_K[0,1)$ 
relative to $V_j$ induced by $(\ol{E},\ol{\nabla})$. Then, 
by using the spectral seuquence of Katz-Oda type 
for $V_j \times A^{|\II|}_K[0,1) \lra V_j \lra \cY'_K$, 
we may replace $\DR$ by the log de Rham complex associated to 
$(\ol{\ol{E}},\ol{\ol{\nabla}})$. Since the restriction of 
$(E,\nabla)$ to $V_j \times A^{|\II|}_K[0,1)$ 
is unipotent relative to $V_j$, $(\ol{\ol{E}},\ol{\ol{\nabla}})$ 
is a successive extension of trivial log-$\nabla$-module $(\cO,d)$ 
on $V_j \times A^{|\II|}_K[0,1)$ relative to $V_j$. So, to prove that 
the map \eqref{qis2} is isomorphic, we may replace $\DR$ by the 
log de Rham complex associated to $(\cO,d)$. In this case, the map 
is shown to be isomorphic in \cite[\S 6]{bc2}. So we have proved the 
claim.

Now we prove the theorem. We prove that the map 
\begin{equation}\label{qis4}
R^qf_{X/\cY,\an *}\cE \lra 
R^qf_{(\ol{X}-{D}^{\ho},\ol{X})/\cY, \logrig *}j^{\d}\cE
\end{equation} 
is an isomorphism, by induction on $q$, using the claim. 
Assume that \eqref{qis4} is an isomorphism up to $q-1$. Let us take any 
diagram as in \eqref{homeo1}. Then, by restricting the isomorphism 
to $]Y'[_{\cY'}=\cY'_K$, we have the isomorphism 
$$ R^tf_{X'/\cY',\an *}\cE \os{=}{\lra} 
R^tf_{(\ol{X}'-{D'}^{\ho},\ol{X}')/\cY', \logrig *}j^{\d}\cE $$ 
for $t<q$. In particular, both hand sides are coherent. So we have 
$$ R^s\sp_*R^tf_{X'/\cY',\an *}\cE = 
R^s\sp_*R^tf_{(\ol{X}'-{D'}^{\ho},\ol{X}')/\cY', \logrig *}j^{\d}\cE =0 $$ 
for $s>0,t<q$. By this and the claim, we can deduce that the map 
$$ 
\sp_*R^qf_{X'/\cY',\an *}\cE \lra 
\sp_*R^qf_{(\ol{X}'-{D'}^{\ho},\ol{X}')/\cY', \logrig *}j^{\d}\cE
$$ 
is an isomorphism. Since this isomorphism holds for any $\cY'$ as 
in the diagram \eqref{homeo1}, we can conclude that the map 
\eqref{qis4} is an isomorphism for $q$. So we are done. 
\end{pf} 

Next we compare the relative 
log analytic cohomology and the relative rigid cohomology, using 
Theorem \ref{thm2.1}. Assume that we are 
in the situation of Theorem \ref{thm2.1}, and assume given open 
immersions $X \hra \ol{X}, j_Y:Y \hra \ol{Y}$ such that 
$X \subseteq (\ol{X},M^{\ve}_{\ol{X}})_{\triv}$, 
$Y \subseteq (\ol{Y},M_{\ol{Y}})_{\triv}$ and $f^{-1}(Y)=X$ holds. 
Then we have $X-D^{\ho} \subseteq (\ol{X},M_{\ol{X}})_{\triv}$. So, 
by \cite[\S 5]{shiho3}, the open immersion $X-D^{\ho} \hra \ol{X}$ 
(which we denote by $j_X$) induces the functor 
$$ j_X^{\dd}: 
I_{\conv}((\ol{X}/\cY)^{\log}) \lra I^{\d}((X-D^{\ho},\ol{X})/\cY). $$
On the other hand, we have the exact functor 
$$ j_Y^{\d}: \Mod(\cO_{]\ol{Y}[_{\cY}}) \lra 
\Mod(j_Y^{\d}\cO_{]\ol{Y}[_{\cY}}) $$ 
which sends coherent modules to coherent modules. Then we have the 
following: 

\begin{thm}\label{thm2.2}
In the above situation, we have the isomorphism 
$$ 
j_Y^{\d}R^qf_{\ol{X}/\cY, \an *}\cE \os{=}{\lra} 
R^qf_{(X-D^{\ho},\ol{X})/\cY, \rig *}j_X^{\dd}\cE. 
$$ 
for any $q \geq 0$. 
In particular, 
$R^qf_{(X-D^{\ho},\ol{X})/\cY, \rig *}j_X^{\dd}\cE$ 
is a coherent $j_Y^{\d}\cO_{]\ol{Y}[_{\cY}}$-module. 
$($Note that we assume the condition $(\star)$ in Theorem \ref{thm2.1}.$)$ 
\end{thm} 

\begin{pf} 
The proof is similar to \cite[5.13]{shiho3}. 
Let us take the embedding system \eqref{embsys1} as in the proof of 
Theorem \ref{thm2.1}. Then we have the diagram of rigid analytic spaces 
\begin{equation*}
\begin{CD}
]\ol{X}^{(\b)}[^{\log}_{\cX^{(\b)}} @>{h^{\log}}>> 
]\ol{Y}[_{\cY} \\ 
@VVV @\vert \\ 
]\ol{X}^{(\b)}[_{\cX^{(\b)}} @>h>> ]\ol{Y}[_{\cY}. 
\end{CD}
\end{equation*} 
Also, us put $X^{(\b)}:=X \times_{\ol{X}} \ol{X}^{(\b)}, 
{D^{\ho}}^{(\b)}:= D^{\ho} \times_{\ol{X}} \ol{X}^{(\b)}$ and 
denote the open immersion $X^{(\b)} \hra \ol{X}^{(\b)}, 
\ol{X}^{(\b)} - {D^{\ho}}^{(\b)} \hra \ol{X}^{(\b)}, 
X^{(\b)} -{D^{\ho}}^{(\b)} \hra \ol{X}^{(\b)}$ by 
$j_{X,1}, j_{X,2}, j_X$, respectively. Then we have admissible open immersions 
\begin{align*}
& j_{X,1}^{\log}: \,\,]X^{(\b)}[^{\log}_{\cX^{(\b)}} \hra 
  ]\ol{X}^{(\b)}[^{\log}_{\cX^{(\b)}}, \\ 
& j_{X,2}^{\log}: \,\,
]\ol{X}^{(\b)} - {D^{\ho}}^{(\b)}[^{\log}_{\cX^{(\b)}} \hra 
  ]\ol{X}^{(\b)}[^{\log}_{\cX^{(\b)}}, \\
& j_X^{\log}: \,\,]X^{(\b)} -{D^{\ho}}^{(\b)}[_{\cX^{(\b)}} = 
  ]X^{(\b)} -{D^{\ho}}^{(\b)}[^{\log}_{\cX^{(\b)}} \hra 
  ]\ol{X}^{(\b)}[^{\log}_{\cX^{(\b)}}
\end{align*}
and it induces the exact functors $j_{X,1}^{\log,\d}, j_{X,2}^{\log,\d}, 
j_X^{\log,\d}$ in the standard way (see \cite[\S 5]{shiho3}) satisfying 
$j_X^{\log,\d} = j_{X,1}^{\log,\d} \circ j_{X,2}^{\log,\d}$. 
On the other hand, we have also an admissible open immersion 
$]X^{(\b)} -{D^{\ho}}^{(\b)}[_{\cX^{(\b)}} \hra 
]\ol{X}^{(\b)}[_{\cX^{(\b)}}$ induced by $j_X$ and it induces the functor 
$j_X^{\d}$ also in the standard way. \par 
Since $(\cX^{(0)},M_{\cX^{(0)}})$ is of Zariski type and the 
morphism $(\ol{X}^{(0)},M_{\ol{X}^{(0)}}) \hra 
(\cX^{(0)}, \allowbreak M_{\cX^{(0)}})$ is an 
exact closed immersion, $(\ol{X}^{(0)},M_{\ol{X}^{(0)}})$ is of Zariski type. 
Hence so is $(\ol{X}^{(m)}, \allowbreak M_{\ol{X}^{(m)}})$. Also, since 
$(\cX^{(0)},M_{\cX^{(0)}})$ and $(\cY,M_{\cY})$ are of Zariski type, so is 
$(\cX^{(m)},M_{\cX^{(m)}})$. Hence, Zariski locally on $\cX^{(m)}$, 
the closed immersion 
$(\ol{X}^{(m)},M_{\ol{X}^{(m)}}) \allowbreak 
\hra (\cX^{(m)},M_{\cX^{(m)}})$ admits a 
factorization 
$$ (\ol{X}^{(m)},M_{\ol{X}^{(m)}}) \hra ({\cX^{(m)}}',M_{{\cX^{(m)}}'}) \lra 
(\cX^{(m)},M_{\cX^{(m)}})$$ 
such that the first map is an exact closed immersion and the second map 
is a formally etale morphism. So, by the argument of \cite[5.9]{shiho3}, 
we see the equalities of functors $\varphi_{X,*}\circ j_X^{\log,\d} = 
j_X^{\d} \circ \varphi_{X,*}$ (which we denote by $j_X^{\dd}$), 
$Rj_X^{\dd} = j_X^{\dd}$ for coherent 
$\cO_{]\ol{X}^{(\b)}[^{\log}_{\cX^{(\b)}}}$-modules. 
Then we have the following diagram: 
{\allowdisplaybreaks{
\begin{align*}
j_Y^{\d}Rf_{\ol{X}/\cY,\an *}\cE & = 
j_Y^{\d}Rh^{\log}_*\DR(]\ol{X}^{(\b)}[^{\log}_{\cX^{(\b)}}/\cY_K,\cE) \\ 
& \os{=}{\lra} 
Rj_Y^{\d}Rh^{\log}_*\DR(]\ol{X}^{(\b)}[^{\log}_{\cX^{(\b)}}/\cY_K,\cE) \\
& \os{=}{\lra}
Rj_Y^{\d}Rh^{\log}_*j_{X,2}^{\log,\d}
\DR(]\ol{X}^{(\b)}[^{\log}_{\cX^{(\b)}}/\cY_K,\cE) 
\,\,\,\,\text{(Theorem \ref{thm2.1})} \\ 
& = 
Rj_Y^{\d}Rh^{\log}_*Rj_{X,2}^{\log,\d}
\DR(]\ol{X}^{(\b)}[^{\log}_{\cX^{(\b)}}/\cY_K,\cE) \\ 
& \os{(\spadesuit)}{\lra} 
Rh^{\log}_*Rj^{\log,\d}_{X,1}Rj_{X,2}^{\log,\d}
\DR(]\ol{X}^{(\b)}[^{\log}_{\cX^{(\b)}}/\cY_K,\cE) \\ 
& = 
Rh^{\log}_*Rj_X^{\log,\d}
\DR(]\ol{X}^{(\b)}[^{\log}_{\cX^{(\b)}}/\cY_K,\cE) \\ 
& = 
Rh_*Rj^{\dd}_X 
\DR(]\ol{X}^{(\b)}[^{\log}_{\cX^{(\b)}}/\cY_K,\cE) \\ 
& = 
Rh_*
j^{\dd}_X 
\DR(]\ol{X}^{(\b)}[^{\log}_{\cX^{(\b)}}/\cY_K,\cE) \\ 
& = 
Rh_*
\DR^{\d}(]\ol{X}^{(\b)}[_{\cX^{(\b)}}/\cY_K,j_X^{\dd}\cE) = 
Rf_{(X-D^{\ho},\ol{X})/\cY,\rig *}j_X^{\dd}\cE. 
\end{align*}}}
So it suffices to prove that $(\spadesuit)$ is a quasi-isomorphism. 
To prove this, we may replace $\b$ by $m \in \N$ and we may replace 
the closed immersion $(\ol{X}^{(m)},M_{\ol{X}^{(m)}}) \hra 
(\cX^{(m)},M_{\cX^{(m)}})$ so that we have 
$(\ol{X}^{(m)},M_{\ol{X}^{(m)}}) = 
(\ol{Y},M_{\ol{Y}}) \times_{(\cY,M_{\cY})} (\cX^{(m)},M_{\cX^{(m)}})$. 
In this case, the morphism $]\ol{X}^{(m)}[^{\log}_{\cX^{(m)}} = 
]\ol{X}^{(m)}[_{\cX^{(m)}} \lra ]\ol{Y}[_{\cY}$ induced by $h^{\log}=h$ 
is quasi-compact and quasi-separated, and so we have the equality of 
functors $j_Y^{\d} \circ h^{\log}_* = h^{\log}_* \circ j^{\log,\d}_{X,1}$. 
So the map $(\spadesuit)$ is a quasi-isomorphism and the theorem 
is proved. 
\end{pf} 

Using Theorem \ref{thm2.2}, we obtain the following theorem, 
which shows the overconvergence of relative rigid cohomology 
for (not necessarily strict) morphisms of pairs which admit `nice' 
log structures: 

\begin{thm}\label{thm2.3}
Let us assume given a proper log smooth morphism 
$f:(\ol{X},M_{\ol{X}}) \lra (\ol{Y},M_{\ol{Y}})$ 
such that 
$M_{\ol{X}}$ 
$($resp. $M_{\ol{Y}})$ is the log structure 
induced by a normal crossing divisor $D$ $($resp. simple normal crossing 
divisor $E)$. We assume 
that $D$ has the decomposition $D=D^{\ho} \cup D^{\ve}$ into sub normal 
crossing divisors $D^{\ho}, D^{\ve}$ satisfying the condition 
$(*)$ in the beginning of this section. Assume moreover that 
we are given open immersions $X \hra \ol{X}, Y \hra \ol{Y}$ 
satisfying $X \subseteq (\ol{X},M^{\ve}_{\ol{X}})_{\triv}, 
Y \subseteq (\ol{Y},M_{\ol{Y}})_{\triv}$ and $f^{-1}(Y)=X$. 
Denote the open immersion $X-D^{\ho} \hra \ol{X}$ by $j_X$. 
Then, 
for any $q \geq 0$ and 
for any locally free 
isocrystal $\cE$ on $(\ol{X}/\cS)^{\log}_{\conv}=
((\ol{X},M_{\ol{X}})/\cS)_{\conv}$ which has nilpotent residues and satisfies 
the condition $(\star)'$ below, 
there exists uniquely an overconvergent isocrystal $\cF$ on 
$(Y,\ol{Y})/\cS_K$ satisfying the following$:$ 
For any $(Y,\ol{Y})$-triple $(Z,\ol{Z},\cZ)$ over $(S,S,\cS)$ 
such that $\ol{Z}$ is smooth over $k$, $E \times_{\ol{Y}} \ol{Z}$ 
is a simple normal crossing divisor in $\ol{Z}$ 
and that 
$\cZ$ is formally smooth over $\cS$, the restriction of $\cF$ to 
$I^{\d}((Z,\ol{Z})/\cS_K,\cZ)$ is given 
functorially by 
$(R^qf_{((X-D^{\ho}) 
\times_Y Z,\ol{X} \times_{\ol{Y}} \ol{Z})/\cZ, \rig *}j_X^{\dd}\cE, 
\epsilon)$, where $\epsilon$ is given by 
\begin{align*}
p_2^*R^qf_{((X-D^{\ho})
\times_Y Z,\ol{X} \times_{\ol{Y}} \ol{Z})/\cZ, \rig *}j_X^{\dd}\cE 
& \os{\simeq}{\rightarrow} 
R^qf_{((X-D^{\ho})\times_Y Z,\ol{X} \times_{\ol{Y}} \ol{Z})/\cZ 
\times_{\cS} \cZ, \rig *}j_X^{\dd}\cE \\ & \os{\simeq}{\leftarrow} 
p_1^*R^qf_{((X-D^{\ho})\times_Y Z,\ol{X} \times_{\ol{Y}} \ol{Z})/\cZ, 
\rig *}j_X^{\dd}\cE. 
\end{align*}
$($Here $p_i$ is the morphism 
$]\ol{Z}[_{\cZ \times_{\cS} \cZ} \lra \,]\ol{Z}[_{\cZ}$ induced by the 
$i$-th projection.$)$ 
\end{thm} 

The condition $(\star)'$ for $\cE$ in the statement of the theorem is 
given as follows: \\
\quad \\
$(\star)'$: \,\,\, 
For any diagram 
\begin{equation*}
\begin{CD}
(\ol{X},M_{\ol{X}}) @<{\varphi_X}<< (\ol{X}_1,M_{\ol{X}_1}) \\ 
@VfVV @VVV \\ 
(\ol{Y},M_{\ol{Y}}) @<{\varphi_Y}<< (\ol{Y}_1,M_{\ol{Y}_1}) @>{i_1}>> 
(\cY_1,M_{\cY_1}) 
\end{CD}
\end{equation*}
such that $\varphi_Y$ is strict, the square is Cartesian, 
$i_1$ is an exact closed immersion and $(\cY_1,M_{\cY_1})$ is 
formally log smooth over $\cS$, the log analytic cohomology 
$R^qf_{\ol{X}_1/\cY_1, \an *}\varphi_X^*\cE$ is a coherent 
$\cO_{]\ol{Y}_1[_{\cY_1}}$-module. 

\begin{pf} 
Let us take a diagram 
$$ 
(\ol{Y},M_{\ol{Y}})  \os{g^{(0)}}{\lla} 
(\ol{Y}^{(0)},M_{\ol{Y}}) \os{i^{(0)}}{\hra} 
(\cY^{(0)},M_{\cY^{(0)}}), $$ 
where $g^{(0)}$ is a strict Zariski covering and $i^{(0)}$ is an exact 
closed immersion into a $p$-adic fine log formal $\cB$-scheme
$(\cY^{(0)},M_{\cY^{(0)}})$ such that $\cY^{(0)}$ is formally smooth over $\cS$  and $M_{\cY^{(0)}}$ is the log structure defined by a relative simple 
normal crossing divisor. For $n \in \N$, let 
$(\ol{Y}^{(n)},M_{\ol{Y}^{(n)}})$ be the $(n+1)$-fold fiber product of 
$(\ol{Y}^{(0)},M_{\ol{Y}^{(0)}})$ over $(\ol{Y},M_{\ol{Y}})$ 
and let $\cY^{(n)}$ be the $(n+1)$-fold fiber product of 
$\cY^{(0)}$ over $\cS$. Then, for each $n$, there exists a log structure 
$M_{\cY^{(n)}}$ (resp. $M_{\cY^{(n)}\times_{\cS}\cY^{(n)}}$) 
on $\cY^{(n)}$ (resp. $\cY^{(n)}\times_{\cS}\cY^{(n)}$) 
such that the canonical closed immersion 
$\ol{Y}^{(n)} \hra \cY^{(n)}$ (resp. $\ol{Y}^{(n)} \hra 
\cY^{(n)}\times_{\cS}\cY^{(n)}$) 
comes from an exact closed immersion 
$(\ol{Y}^{(n)},M_{\ol{Y}^{(n)}}) \hra (\cY^{(n)},M_{\cY^{(n)}})$ 
(resp. 
$(\ol{Y}^{(n)},M_{\ol{Y}^{(n)}}) \hra (\cY^{(n)}\times_{\cS}\cY^{(n)},
M_{\cY^{(n)}\times_{\cS}\cY^{(n)}})$) and $(\cY^{(n)},M_{\cY^{(n)}})$ 
(resp. $(\cY^{(n)}\times_{\cS}\cY^{(n)},
M_{\cY^{(n)}\times_{\cS}\cY^{(n)}})$) is formally log smooth over $\cS$. 
Let us denote the pull-back of $(\ol{X},M_{\ol{X}}), X, D^{\ho}$ by 
the morphism $(\ol{Y}^{(n)},M_{\ol{Y}^{(n)}}) \lra (\ol{Y},M_{\ol{Y}})$ by 
$(\ol{X}^{(n)},M_{\ol{X}^{(n)}}), X^{(n)}, {D^{\ho}}^{(n)}$, respectively. 
Then, by Theorem \ref{thm2.2} and the condition $(\star)'$, 
$R^qf_{(X^{(n)}-{D^{\ho}}^{(n)},\ol{X}^{(n)})/\cY^{(n)}, \rig *}j_X^{\dd}\cE$ 
(resp. 
$R^qf_{(X^{(n)}-{D^{\ho}}^{(n)},\ol{X}^{(n)})/
\cY^{(n)}\times_{\cS}\cY^{(n)}, \rig *}j_X^{\dd}\cE$) 
is a coherent $j_Y^{\d}\cO_{]\ol{Y}^{(n)}[_{\cY^{(n)}}}$-module 
(resp. a coherent $j_Y^{\d}\cO_{]\ol{Y}^{(n)}[_{\cY^{(n)} \times_{\cS} 
\cY^{(n)}}}$-module). 
Then, by the base change theorem of Tsuzuki (\cite[2.3.1]{tsuzuki3}), 
the arrows in the diagram 
\begin{align*} 
p_2^*R^qf_{(X^{(n)}-{D^{\ho}}^{(n)},\ol{X}^{(n)})/\cY^{(n)}, \rig *}
j_X^{\dd}\cE 
& \ra 
R^qf_{(X^{(n)}-{D^{\ho}}^{(n)},\ol{X}^{(n)})/\cY^{(n)}\times_{\cS}\cY^{(n)}, 
\rig *}j_X^{\dd}\cE \\ 
& \leftarrow 
p_1^*R^qf_{(X^{(n)}-{D^{\ho}}^{(n)},\ol{X}^{(n)})/\cY^{(n)}, \rig *}
j_X^{\dd}\cE 
\end{align*} 
(where $p_i$ denotes the $i$-th projection 
$]\ol{Y}^{(n)}[_{\cY^{(n)} \times_{\cS} \cY^{(n)}} \lra 
]\ol{Y}^{(n)}[_{\cY^{(n)}}$) are isomorphisms. 
If we denote the composite 
of the above arrows by $\epsilon^{(n)}$, one can see (by using 
\cite[2.3.1]{tsuzuki3} again) that $\cF^{(n)}:=
(R^qf_{(X^{(n)}-{D^{\ho}}^{(n)},
\ol{X}^{(n)})/\cY^{(n)}, \rig *}\cE,\epsilon^{(n)})$ 
defines an overconvergent isocrystal on $(Y^{(n)},\ol{Y}^{(n)})/\cS_K$ 
and that $\cF^{(n)}$ is compatible with respect to $n$. 
So $\{\cF^{(n)}\}_{n=0,1,2}$ descents to an overconvergent 
isocrystal $\cF$ on $(Y,\ol{Y})/\cS_K$. \par 
Since the triple $(Y^{(n)},\ol{Y}^{(n)},\cY^{(n)})$ satisfies the 
condition required for $(Z,\ol{Z},\cZ)$ in the theorem, we see that the 
image of $\cF$ in $I^{\d}((Y^{(n)},\ol{Y}^{(n)})/\cS_K)$ should be 
functorially isomorphic to $\cF^{(n)}$ in the previous paragraph. 
So we see that the condition in the statement of the theorem 
characterizes $\cF$ uniquely. \par 
Finally we prove that the overconvergent isocrystal $\cF$ satisfies the 
required condition. For a triple $(Z,\ol{Z},\cZ)$ as in the statement of 
the theorem, let us put $M_{\ol{Z}} := M_{\ol{Y}}|_{\ol{Z}}$ 
(this is nothing but the log structure associated to the 
simple normal crossing divisor $E \times_{\ol{Y}} \ol{Z}$). Then, 
Zariski locally on $\cZ$, there exists a fine log structure 
$M_{\cZ}$ on $\cZ$ such that the closed immersion $\ol{Z} \hra \cZ$ 
comes from an exact closed immersion $(\ol{Z},M_{\ol{Z}}) \hra 
(\cZ,M_{\cZ})$ and that $(\cZ,M_{\cZ})$ is formally log smooth over $\cS$. 
Then, by Theorem \ref{thm2.2}, 
$R^qf_{((X-D^{\ho}) 
\times_Y Z,\ol{X} \times_{\ol{Y}} \ol{Z})/\cZ, \rig *}j_X^{\dd}\cE$ is 
a coherent $j^{\d}\cO_{]\ol{Z}[_{\cZ}}$-module. 
Then, by 
\cite[2.3.1]{tsuzuki3} again, we see that the restriction of $\cF$ to 
$I^{\d}((Z,\ol{Z})/\cS_K,\cZ)$ is given by 
$(R^qf_{((X-D^{\ho}) \times_Y Z,
\ol{X} \times_{\ol{Y}} \ol{Z})/\cZ, \rig *}j_X^{\dd}\cE, 
\epsilon)$ as in the statement of the theorem. 
Finally, the functoriality of the expression above is proved as the 
proof of \cite[4.8]{shiho3}. So we are done. 
\end{pf} 

As for Frobenius structure, we have the following: 

\begin{thm}\label{thm2.4} 
With the situation in Theorem \ref{thm2.3}, assume moreover that 
$k$ is perfect and that 
$j_X^{\dd}\cE$ has a Frobenius structure. Then the overconvergent isocrystal 
$\cF$ on $(Y,\ol{Y})/\cS_K$ has a canonical Frobenius structure. 
\end{thm} 

\begin{pf} 
The proof is similar to \cite[5.16]{shiho3} (see also 
\cite[3.3.3, 4.1.4]{tsuzuki3}). First, by the argument developed in 
\cite[5.16]{shiho3}, we may reduce to the case that 
$Y=\ol{Y}$ is equal to $S$. (In the proof of \cite[5.16]{shiho3}, 
we once used the analytically flat base change theorem of log analytic 
cohomology. Here we use Tsuzuki's base change theorem 
(\cite[2.3.1]{tsuzuki3}) 
instead of it.) 
So, we are reduced to the claim that the endomorphism on 
$H^q_{\rig}(X-D^{\ho}/K, j_X^{\dd}\cE)$ induced by the Frobenius structure 
on $j_X^{\dd}\cE$ 
is an isomorphism. We can prove this by imitating the proof of 
\cite[9.1.2]{kedlaya2}: Let us consider the following 
two claims: \\
\quad \\
$(a)_n$ \,\,\, The Frobenius endomorphism on 
$H^q_{\rig}(X/K,\cG)$ is an isomorphism for any $X$ smooth over $k$ of 
dimension at most $n$ and any overconvergent $F$-isocrystal $\cG$ on $X$. \\
$(b)_n$ \,\,\, 
The Frobenius endomorphism on $H^q_{Z,\rig}(X/K,\cG)$ is an isomorphism 
for any closed immersion $Z \hra X$ of a geometrically reduced scheme $Z$ 
over $k$ of dimension at most $n$ into a scheme $X$ smooth over $k$ and 
any overconvergent $F$-isocrystal $\cG$ on $X$. \\
\quad \\
($(a)_0$ is obvious and $(b)_{-1}$ is vacuous.) 
Then, by a similar argument to the proof of 
\cite[9.1.2]{kedlaya2}, we can prove the implications 
$(b)_{n-1}+(a)_n \,\Longrightarrow (b)_n$ and $(a)_{n-1}+(b)_{n-1} 
\,\Longrightarrow\,(a)_n$: The main point is that all the arguments 
in the proof of \cite[9.1.2]{kedlaya2} are compatible with Frobenius. 
(The most important point is the compatibility of Gysin isomorphism 
with Frobenius, which is proved in \cite[4.1.1]{tsuzuki2}.) 
So the Frobenius endomorphism on the rigid cohomology of smooth variety 
with coefficient is always an isomorphism. So we are done. 
\end{pf} 

\begin{rem} 
The conditions $(\star)$ and $(\star)'$ are satisfied if we assume $f$ is 
integral, by \cite[4.7]{shiho3}. 
We will provide a slightly different situation where the 
conditions $(\star)$ and $(\star)'$ are satisfied, in the next section. 
\end{rem}


\section{Log blow-up invariance of relative log analytic cohomology} 

In this section, we prove an invariance property of 
relative log analytic cohomology under log blow-ups in certain case. \par 
First we recall the notion of log blow-up. (See \cite[6.1]{illusie}, 
\cite[2.1]{saito}, \cite[4]{niziol}.) Let $(X,M_X)$ be an fs log scheme. 
A sheaf of ideals $\cI \subseteq M_X$ is called a coherent ideal 
if, etale locally on $X$, there exists a chart $P \ra M_{X}$ of $M_X$ by a 
torsion-free fs monoid $P$ and an ideal $I$ of $P$ such that 
$IM_X = \cI$ holds. Then the log blow-up $(X_{\cI},M_{X_{\cI}}) \lra 
(X,M_X)$ of $(X,M_X)$ along $\cI$ is locally defined as follows: 
Let 
$I^{\sat} := \{a \in P^{\gp}\,\vert\, \exists n \geq 1, a^n \in I^n\}$ 
(the saturation of the ideal $I$). 
Then $X_{\cI}$ is defined to be $X \times_{\Spec \Z[P]} 
\Proj \Z[\oplus_{n=0}^{\infty}(I^{\sat})^n]$. It has an open covering 
by the schemes of the form $X \times_{\Spec \Z[P]} \Spec \Z[P_a^{\sat}] \, 
(a \in P)$, where $P_a := \bigcup a^{-n}I^n$ and 
$P_a^{\sat} := \{x \in P_a^{\gp}\,\vert\, \exists n \geq 1, x^n \in P_a\}$ 
(the saturation of $P_a$). We define the log structure $M_{X_{\cI}}$ 
as the fs log structure whose restriction to 
$X \times_{\Spec \Z[P]} \Spec \Z[P_a^{\sat}]$ is given as 
the pull-back of the canonical log structure $M_{P_a^{\sat}}$ on 
$\Spec \Z[P_a^{\sat}]$. Log blow-ups are log etale. 
We can define the notion of low blow-up also for formal fs log schemes. \par 
Let $R$ be a ring (resp. adic topolgical ring), let $P$ be a torsion-free 
fs monoid 
and let $I$ be an ideal of $P$. Let us denote the canonical log structure 
on $\Spec R[P]$ (resp. $\Spf R\{P\}$) by $M_P$. Then we denote the 
log blow-up of $(\Spec R[P],M_P)$ (resp. $(\Spf R\{P\},M_P)$) along 
the coherent ideal $IM_P$ by $(\Bl_I(P)_R, M_{P,I})$. \par 
Next let us recall the notion of log regular log scheme (\cite[2.2]{niziol}, 
\cite[2]{kkato2}) and several properties of it. 
For a log scheme $(X,M_X)$ and $x \in X$, let us denote the 
ideal of $\cO_{X,\ol{x}}$ generated by the image of $M_{X,\ol{x}} - 
\cO^{\times}_{X,\ol{x}}$ by $I(M_{X},\ol{x})$. An fs log scheme 
$(X,M_X)$ is called log regular if, for any $x \in X$, 
$\cO_{X,\ol{x}}/I(M_X,\ol{x})$ is regular and the equality 
$\dim (\cO_{X,\ol{x}}) = \dim (\cO_{X,\ol{x}}/I(M_X,\ol{x})) + 
\rk (M_{X,\ol{x}}^{\gp}/\cO^{\times}_{X,\ol{x}})$ holds. 
If $(X,M_X)$ is log regular, $X$ is normal and 
$M_X = j_*\cO^{\times}_{(X,M_X)_{\triv}} \cap \cO_{X}$ holds, where 
$j$ denotes the open immersion $(X,M_X)_{\triv} \hra X$.
An fs log scheme which is log smooth over a log regular scheme is 
again log regular. In particular, any fs log scheme which is log smooth 
over $k$ is log regular. In this case, the converse is also true if $k$ is 
perfect. 
For a log regular log scheme $(X,M_X)$, $X$ is regular if and only if, 
for any $x \in X$, $M_{X,\ol{x}}/\cO^{\times}_{X,\ol{x}}$ is isomorphic to 
$\N^{r(x)}$ for some $r(x)$ (depending on $x$). In this case, 
the log structure $M_X$ is the one 
associated to some normal crossing divisor $D$ on $X$. \par 
It is known that, for a log regular scheme $(X,M_X)$, there exists 
a log blow-up $(\ti{X},M_{\ti{X}}) \lra (X,M_X)$ such that 
$\ti{X}$ is regular (\cite[5.7,5.10]{niziol}). \par 
Now we prove a log blow-up invariance property for log analytic 
cohomology: 

\begin{thm}\label{lbui}
Let us assume given the diagram 
\begin{equation*}
(\ol{X}',M_{\ol{X}'}) \os{\varphi}{\lra} 
(\ol{X},M_{\ol{X}}) \os{f}{\lra} 
(\ol{Y},M_{\ol{Y}}) \os{\iota}{\hra} 
(\cY,M_{\cY}), 
\end{equation*}
where $(\ol{X},M_{\ol{X}}), (\ol{Y},M_{\ol{Y}})$ 
are fs log schemes, $(\cY,M_{\cY})$ is an fs log formal $\cS$-scheme 
of Zariski type formally log smooth over $\cS$, 
$f$ is log smooth, $\iota$ is an exact closed immersion and 
$\varphi$ is a log blow-up of $(\ol{X},M_{\ol{X}})$ by some coherent ideal 
of $M_X$. Then, for a locally free isocrystal $\cE$ on 
$(\ol{X}/\cY)^{\log}_{\conv}$, we have the quasi-isomorphism 
$$ Rf_{\ol{X}/\cY, \an *}\cE \os{=}{\lra} 
Rf_{\ol{X}'/\cY, \an *}\varphi^*\cE. $$ 
\end{thm} 

\begin{pf} 
We may work etale locally on $\ol{X}$ and Zariski locally on $\cY$. 
So we may assume that there exists a chart 
$(Q \ra M_{\cY}, R \ra M_X, Q \os{\alpha}{\ra} R)$ of $\iota \circ f$ 
such that $\alpha$ is injective, 
$|\Coker(\alpha^{\gp})_{\tor}|$ is prime to $p$ and that $f$ factors as 
$$ 
(\ol{X},M_{\ol{X}}) \lra (\ol{Y},M_{\ol{Y}}) \times_{(\Spec k[Q],M_Q)} 
(\Spec k[R],M_R) \lra (\ol{Y},M_{\ol{Y}})
$$ 
with the first map strict etale. Then we can form the Cartesian diagram 
\begin{equation*}
\begin{CD}
(\ol{X},M_{\ol{X}}) @>>> (\ol{Y},M_{\ol{Y}}) \times_{(\Spec k[Q],M_Q)} 
(\Spec k[R],M_R) @>>> (\ol{Y},M_{\ol{Y}}) \\ 
@VVV @VVV @VVV \\ 
(\cX,M_{\cX}) @>>> (\cY,M_{\ol{Y}}) \times_{(\Spf V\{Q\}, M_Q)}
(\Spf V\{R\},M_R) @>>> (\cY,M_{\cY}) 
\end{CD}
\end{equation*}
with the lower left horizontal arrow strict formally etale, 
by shrinking $\ol{X}$. Then $(\cX,M_{\cX}) \lra (\cY,M_{\cY})$ is 
a log smooth lift of the morphism $f$. 
We may also assume that the log structure of $\cX$ is fs. \par 
Next note that, etale locally on $\ol{X}$, there exists a 
chart $\psi: P \ra M_{\ol{X}}$ of $M_{\ol{X}}$ by a 
torsion-free fs monoid $P$, 
an ideal $I \subseteq P$ and a Cartesian diagram 
\begin{equation}\label{lbu}
\begin{CD}
(\Bl_I(P)_k,M_{P,I}) @<<< (\ol{X}',M_{\ol{X}'}) \\ 
@VVV @V{\varphi}VV \\ 
(\Spec k[P],M_P) @<<< (\ol{X},M_{\ol{X}}), 
\end{CD}
\end{equation}
where the lower horizontal arrow is the strict 
morphism induced by the chart. Note that, for $x \in \ol{X}$, 
we can lift the homomorphism 
$\psi^{\gp}:P^{\gp} \ra M_{\ol{X},\ol{x}}^{\gp}$ to a homomorphism 
$\psi':P^{\gp} \ra M^{\gp}_{\cX,\ol{x}}$ such that the composite 
$P^{\gp} \os{\psi'}{\ra} 
M^{\gp}_{\cX,\ol{x}} \ra M^{\gp}_{\cX,\ol{x}}/\cO_{\cX,\ol{x}}^{\times} 
(=M^{\gp}_{\ol{X},\ol{x}}/\cO_{\ol{X},\ol{x}}^{\times})$ is surjective. 
So, if we put $P' := {\psi'}^{-1}(M_{\cX,\ol{x}})$, 
$\psi':P' \ra M_{\cX,\ol{x}}$ extends to a chart of $M_{\cX}$ on an 
etale neighborhood of $x$ (\cite[2.10]{kkato1}). Moreover, by 
the exactness of $(X,M_X) \hra (\cX,M_{\cX})$, we have $P \subseteq P'$ and 
$P' \ra M_{\cX} \ra M_{\ol{X}}$ is again a chart of $M_{\ol{X}}$. Also, 
the Cartesian diagram \eqref{lbu} induces the similar Cartesian diagram 
with $P,I$ replaced by $P',IP'$. As a conclusion, we may assume that 
(by putting $P := P', I := IP'$) the chart $\psi:P \ra M_{\ol{X}}$ 
of $M_{\ol{X}}$ extends 
to a chart $\psi':P \ra M_{\cX}$ of $M_{\cX}$. Now we put 
$$(\cX',M_{\cX'}) := (\cX,M_{\cX}) \times_{(\Spf V\{P\},M_P)} 
(\Bl_I(P)_V,M_{P,I}) \os{\varphi'}{\lra} (\cX,M_{\cX}). $$ 
Then $\varphi'$ is a formally log etale lift of $\varphi$. 
Let us denote the map of rigid analytic spaces 
$]\ol{X}'[_{\cX'} \lra ]\ol{X}[_{\cX}, ]\ol{X}[_{\cX} \lra ]\ol{Y}[_{\cY}$ 
by $\varphi'_K, h$, respectively. Then we have the map 
{\small{ $$ 
Rf_{\ol{X}/\cY,\an *\cE} = 
Rh_*\DR(]\ol{X}[_{\cX}/\cY_K,\cE) \lra 
Rh_*R\varphi'_{K,*}{\varphi'}_K^*\DR(]\ol{X}[_{\cX}/\cY_K,\cE) = 
Rf_{\ol{X}'/\cY, \an *}\cE.$$}}
We prove that this is a quasi-isomorphism. 
To prove this, 
it suffices to prove the equality 
$R\varphi'_{K,*}{\varphi'_K}^*\cO_{]\ol{X}[_{\cX}} = \cO_{]\ol{X}[_{\cX}}$, 
and it is reduced to the equality 
$R\varphi'_{*}{\varphi'}^*\cO_{\cX} = \cO_{\cX}$. 
Let us denote the uniformizer of $V$ by $\pi$ and 
let $\varphi'_n: \cX'_n \lra \cX_n$ be $\varphi' \otimes_V V/\pi^nV$. 
Then it suffices to prove the equality 
\begin{equation}\label{vanish}
R\varphi'_{n,*}{\varphi'}_n^*\cO_{\cX_n} = \cO_{\cX_n}. 
\end{equation}
Note that we have $\cX'_n = \cX_n \times_{\Spec \Z[P]} \Bl_I(P)_{\Z}$, 
and note also that we have the vanishing 
$\Tor_i^{\Z[P]}(\Z[P'],\cO_{\cX_n}) =0$ for 
any $i \geq 1, n \in \N$ and injective homomorphism $P \hra P'$ of integral 
monoids: Indeed, if $n=1$, this follows from \cite[6.1(ii)]{kkato2}, 
since $(\cX_1,M_{\cX}|_{\cX_1})$, being log smooth over $k$, is 
log regular with a chart $P \ra M_{\cX} \ra M_{\cX_1}$. 
For general $n$, we can see by induction, using the exact sequence 
$$ 0 \lra \cO_{\cX_{n-1}} \os{\pi}{\lra} \cO_{\cX_n} \lra \cO_{\cX_1} 
\lra 0. $$
(Note that, since $(\cX_n,M_{\cX_n})$ is log smooth over $\Spf V/\pi^nV$, 
it is flat over $V/\pi^nV$.) By this vanishing of $\Tor$, \eqref{vanish} is 
reduced to the equality $R\ul{\varphi}_* \ul{\varphi}^* 
\cO_{\Spec \Z[P]} = \allowbreak 
\cO \allowbreak {}_{\Spec \Z[P]}$, where $\ul{\varphi}$ is the map 
$\Bl_I(P)_{\Z} \lra \Spec \Z[P]$. This is true by 
\cite[Ch I, \S 3, Cor. 1]{kkms}. (There they treat the case of $\Spec k[P]$, 
but the arguments work for $\Spec \Z[P]$. See also \cite[11.3]{kkato2}.) 
So the proof is finished. 
\end{pf} 

Next we deduce a consequence of Theorem \ref{lbui}, which is useful 
in later sections. To describe this, first we 
define the notion of log normal crossing divisor as follows: 

\begin{defn}\label{lncd}
Let $(X,M_X^{\ve})$ be a log regular log scheme. Then a Cartier divisor 
$D \subseteq X$ is called a log normal crossing divisor if, for any 
$x \in D$, $D \times_X \Spec \cO_{X,\ol{x}} \subseteq \Spec \cO_{X,\ol{x}}$ 
is defined by the equation of the form $t_1t_2\cdots t_r=0$ 
$(t_i \in \cO_{X,\ol{x}}-\cO_{X,\ol{x}}^{\times})$ such that 
the loci $\{t_i=0\} \,(1 \leq i \leq r)$ is a regular divisor 
in $\Spec \cO_{X,\ol{x}}/I(M^{\ve}_X,\ol{x})$ and that the locus 
$\{t_1\cdots t_r=0\}$ is a simple normal crossing divisor 
in $\Spec \cO_{X,\ol{x}}/I(M^{\ve}_X,\ol{x})$. 
\end{defn} 

\begin{rem}\label{lncdrem}
The conditions required for $t_i$'s in the above definition is equivalent 
to the following condition: For any subset $I$ of $\{1,\cdots, r\}$, 
the locus $\{t_i=0 \,(i \in I)\}$ is a regular closed subscheme of 
$\Spec \cO_{X,\ol{x}}/I(M^{\ve}_X,\ol{x})$ of codimension $|I|$. 
\end{rem} 

\begin{prop}\label{lncdprop}
In the situation in Definition \ref{lncd}, let us define new log structures 
$M_X^{\ho},M_X$ on $X$ by 
$M_X^{\ho}:=j_*\cO^{\times}_{X-D} \cap \cO_{X}$ and 
$M_X := M_X^{\ho} \oplus_{\cO_{X}^{\times}} M_X^{\ve}$. Then$:$ 
\begin{enumerate} 
\item 
Etale locally, We have $M_X^{\ho} = \cO_X^{\times}t_1^{\N}\cdots t_r^{\N}$ and 
it is associated to the monoid homomorphism $\varphi: \N^r \lra \cO_X; \, 
e_i \mapsto t_i$. $($Here $e_i\,(1 \leq i \leq r)$ is a canonical basis 
of $\N^r.)$ 
\item 
The log scheme $(X,M_X)$ is log regular. Also, if we define 
$D_{[i]} \,(i \geq 0)$ by 
$$ D_{[0]}:=X, \qquad D_{[1]} := \text{the normalization of $D$}, $$ 
$$ D_{[i]}:= \text{$i$-fold fiber product of $D_{[1]}$ over $X$}, $$ 
the log scheme $(D_{[i]},M^{\ve}_{D_{[i]}}):=
(D_{[i]},M^{\ve}_X |_{D_{[i]}})$ is also log regular. 
\end{enumerate}
\end{prop} 

\begin{pf} 
We may work etale locally. 
For $I \subseteq \{1,\cdots,r\}$, let us denote the closed subscheme 
of $X$ defined by the equation $t_i=0\,(i \in I)$ by $D_I$. 
First we prove that $(D_I,M_X^{\ve}|_{D_I})$ is log regular by 
induction on $|I|$. Write $I = \{i\} \coprod I'$ with $i \notin I'$. 
Then we have 
$$\dim (\cO_{D_{I'},\ol{x}}) = 
\dim (\cO_{D_{I'},\ol{x}}/I(M_X^{\ve}|_{D_{I'}},\ol{x})) 
+ \rk (M_X^{\ve}|_{{D_{I'}},\ol{x}}^{\gp}/\cO^{\times}_{D_{I'},\ol{x}})$$ 
by induction hypothesis. Then, since $t_i$ is non-zero and non-invertible in 
$\cO_{D_{I'},\ol{x}}/I(M_X^{\ve}\allowbreak |_{D_{I'}},\ol{x})$, 
so is $t_i$ is in $\cO_{D_{I'},\ol{x}}$. Since 
$\cO_{D_{I'},\ol{x}}$ is an integral domain (by induction hypothesis), 
we have 
\begin{align*}
\dim (\cO_{D_I,\ol{x}}) & = \dim (\cO_{D_{I'},\ol{x}}) -1, \\
\dim (\cO_{D_I,\ol{x}}/I(M_X^{\ve}|_{D_I},\ol{x})) & = 
\dim (\cO_{D_{I'},\ol{x}}/I(M_X^{\ve}|_{D_{I'}},\ol{x})) -1.
\end{align*} 
So we see the equality 
$$ 
\dim (\cO_{D_I,\ol{x}}) = \dim (\cO_{D_I,\ol{x}}/I(M_X^{\ve}|_{D_I},\ol{x})) 
+ \rk (M_X^{\ve}|_{{D_I},\ol{x}}^{\gp}/\cO^{\times}_{D_I,\ol{x}}),
$$
and $\cO_{D_I,\ol{x}}/I(M_X^{\ve}|_{D_I},\ol{x})$ is regular by assumption. 
So $(D_I,M_X|_{D_I})$ is log regular. 
In particular, $D_{\{i\}}$ is normal. Then we see that any local 
section $f$ in $M_X^{\ho}$ can be written uniquely as 
$f=ut_1^{n_1}\cdots t_r^{n_r}$ for some $u \in \cO_X^{\times}, n_i \in \N$. 
So we have the assertion (1). Also, we see that $D_{[i]}$ is locally 
written as the disjoint union of $D_I$'s with $|I|=i+1$. So 
$(D_{[i]},M^{\ve}_{D_{[i]}})$ is log regular. \par 
Finally we prove that $(X,M_X)$ is log regular. If we take a chart 
$\psi:P \ra \cO_{X}$ of $M_X^{\ve}$ with $P$ an fs monoid, we see by (1) that 
the homomorphism $P \oplus \N^r \lra \cO_X; \, (p,n) \mapsto 
\psi(p)\varphi(n)$ is a chart of $M_X$. So $M_X$ is an fs log structure. 
Then, since we have 
$$ \dim (\cO_{X,\ol{x}}/I(M_X,\ol{x})) = 
\dim (\cO_{X,\ol{x}}/I(M_X^{\ve},\ol{x})+(t_1,\cdots,t_r)) = 
\dim (\cO_{X,\ol{x}}/I(M_X^{\ve},\ol{x}))-r, $$
$$ 
\rk (M_{X,\ol{x}}^{\gp}/\cO^{\times}_{X,\ol{x}}) = 
\rk (M_{X,\ol{x}}^{\ve,\gp}/\cO^{\times}_{X,\ol{x}})+r, $$
we have the equality 
$$\dim (\cO_{X,\ol{x}}) = \dim (\cO_{X,\ol{x}}/I(M_X,\ol{x})) + 
\rk (M_{X,\ol{x}}^{\gp}/\cO^{\times}_{X,\ol{x}}).$$ 
Since $\cO_{X,\ol{x}}/I(M_X,\ol{x}) = 
\cO_{X,\ol{x}}/(I(M_X^{\ve},\ol{x})+(t_1,\cdots,t_r))$ is regular by 
assumption, we see that $(X,M_X)$ is log regular. 
\end{pf} 

\begin{exam} 
Let $X$ be a smooth scheme over $k$ and let $D^{\ho}, D^{\ve}$ are normal 
crossing divisors on $X$ such that $D:=D^{\ho} \cup D^{\ve}$ is again 
a normal crossing divisor. Then, if we denote the log structure on $X$ 
associated to $D^{\ve}$ by $M^{\ve}_X$, $D^{\ho}$ is a log normal 
crossing divisor on $(X,M_X^{\ve})$. Moreover, $M_X^{\ho}$ 
(resp. $M_X$) 
above is nothing but the log structure associated to $D^{\ho}$ (resp. $D$). 
\end{exam} 

Now let us assume given a proper log smooth integral 
morphism of fs log schemes 
$f':(\ol{X}',M^{\ve}_{\ol{X}'}) \lra (\ol{Y},M_{\ol{Y}})$ over $S$ 
such that 
$\ol{Y}$ is smooth over $S$ and that $M_{\ol{Y}}$ is the log structure 
associated to a simple normal crossing divisor on $\ol{Y}$. 
(Then $(\ol{X}',M^{\ve}_{\ol{X}'})$, being log smooth over $S$, is 
log regular.)
Let $D' \subseteq \ol{X}'$ be a log normal crossing divisor, and 
let us define $M_{\ol{X}'}^{\ho}, M_{\ol{X}'}, D'_{[i]} \,(i \geq 0)$ 
as above. Assume moreover the following conditions: 
\begin{enumerate} 
\item 
$(\ol{X}',M_{\ol{X}'})$ is log smooth over $(\ol{Y},M_{\ol{Y}})$. 
\item 
For any $i\geq 0$, $(D'_{[i]},M_{\ol{X}'}^{\ve}|_{D'_{[i]}})$ is 
log smooth over $(\ol{Y},M_{\ol{Y}})$. 
\end{enumerate} 
Take a log blow-up $\varphi^{\ve}:(\ol{X},M_{\ol{X}}^{\ve}) \lra 
(\ol{X}',M^{\ve}_{\ol{X}'})$ of $(\ol{X}',M^{\ve}_{\ol{X}'})$ such that 
$\ol{X}$ is regular, and put 
$$ (\ol{X},M_{\ol{X}}) := (\ol{X},M_{\ol{X}}^{\ve}) 
\times_{(\ol{X}',M_{\ol{X}'}^{\ve})} (\ol{X}',M_{\ol{X}'}) 
\os{\varphi}{\lra} (\ol{X}',M_{\ol{X}'}),$$
where $\varphi=\varphi^{\ve} \times \id$. (It is also a log blow-up.) 
Let us put $D^{\ho}:=\varphi^*D'$ and let $D^{\ve}$ be a normal crossing 
divisor corresponding to $M_{\ol{X}}^{\ve}$. Then 
$D := D^{\ho} \cup D^{\ve}$ is a normal crossing divisor corresponding to 
$M_{\ol{X}}$ and they satisfies the condition $(*)$ in the beginning of 
Section 2. (So the morphism 
$f := \varphi \circ f':(\ol{X},M_{\ol{X}}) \lra (\ol{Y},M_{\ol{Y}})$ is in the 
situation we treated in Section 2.) With this notation, we have the 
following: 

\begin{prop}\label{starok}
With the above notation, let $\cE_0$ be a locally free isocrystal 
on $(\ol{X}'/\cS)^{\log}_{\conv} = 
((\ol{X}',M_{\ol{X}'})/\cS)_{\conv}$ and put $\cE := \varphi^*\cE_0$. 
Then$:$ 
\begin{enumerate}
\item 
For any exact closed immersion $(\ol{Y},M_{\ol{Y}}) \hra (\cY,M_{\cY})$ 
such that $(\cY,M_{\cY})$ is formally log smooth over $\cS$, 
$\cE$ satisfies the condition $(\star)$ in Theorem \ref{thm2.1}. 
\item 
$\cE$ satisfies the condition $(\star)'$ in Theorem \ref{thm2.3}. 
\end{enumerate} 
\end{prop} 

\begin{pf} 
We only prove (1). (The proof of (2) is similar.) 
Since $(\ol{X}',M_{\ol{X}'}) \lra (\ol{Y},M_{\ol{Y}})$ is proper 
log smooth integral and $(\ol{Y},M_{\ol{Y}})$ is log smooth over $k$, 
the relative log analytic cohomology $R^qf_{\ol{X}'/\cY, \an *}\cE_0$ is 
coherent by \cite[4.7]{shiho3}. Moreover, by Theorem \ref{lbui}, 
we have the isomorphism 
$R^qf_{\ol{X}'/\cY, \an *}\cE_0 = R^qf_{\ol{X}/\cY, \an *}\cE$. 
So $R^qf_{\ol{X}/\cY, \an *}\cE$ is coherent. 
\end{pf} 

\begin{rem}\label{starokrem}
As a consquence of Proposition \ref{starok}, we have the following: 
In the situation of Proposition \ref{starok}, assume moreover that 
$\cE$ has nilpotent residues. Then we have Theorems \ref{thm2.1}, 
\ref{thm2.2} and \ref{thm2.3} for $\cE$. 
\end{rem} 

\begin{rem} 
The reason we need Proposition \ref{starok} is that it is not 
always true that the morphism $f$ is integral. 
\end{rem}

Finally, we prove a proposition concerning 
log normal crossing divisors which we use in later sections: 

\begin{prop}\label{lncdpullback}
Let $f: (X,M_X^{\ve}) \lra (Y,M_Y)$ be a log smooth integral universally 
saturated morphism of fs log schemes such that $(Y,M_Y)$ is log smooth over 
$S$ and let $D \subseteq X$ be a log 
normal crossing divisor. Define $M_X^{\ho},M_X,D_{[i]}$ as in 
Proposition \ref{lncdprop} and assume that $(X,M_X), 
(D_{[i]},M_X^{\ve}|_{D_{[i]}})$ are log smooth over $(Y,M_Y)$. Let us assume 
given the Cartesian diagram 
\begin{equation*}
\begin{CD}
(X',M_{X'}^{\ve}) @>>> (Y',M_{Y'}) \\ 
@VVV @VVV \\ 
(X,M_X^{\ve}) @>f>> (Y,M_Y)
\end{CD}
\end{equation*}
such that $(Y',M_{Y'})$ is again an fs log scheme which is log smooth 
over $S$. Put $D':=D \times_X X'$. Then $D'$ is a log normal 
crossing divisor in $X'$ and if we define $M_{X'}^{\ho},M_{X'},D'_{[i]}$ 
using $D'$, $(X',M_{X'})$ and $(D'_{[i]},M_{X'}^{\ve}|_{D'_{[i]}})$ are
log smooth over $(Y',M_{Y'})$. 
\end{prop} 

\begin{pf} 
Let $x' \in X'$, let $x$ be its image in
$X$ and let us take $t_1,\cdots,t_r \in \cO_{X,\ol{x}}-
\cO_{X,\ol{x}}^{\times}$ such that $D \times_X \Spec \cO_{X,\ol{x}}$ is 
defined by $t_1\cdots t_r=0$ and that they satisfy the condition required 
in Remark \ref{lncdrem}. Then $D' \times_{X'} \Spec \cO_{X',\ol{x}'}$ 
is also defined by 
$t_1\cdots t_r=0$ in $\Spec \cO_{X',\ol{x}'}$. To prove the proposition, 
it suffices to prove that the elements $t_i\,(1 \leq i \leq r)$ 
satisfy the condition required 
in Remark \ref{lncdrem} (with $D,X,x$ replaced by $D',X',x'$): 
Indeed, this assertion implies that $D'$ is a log normal crossing divisor 
in $X'$ and that we have the isomorphisms 
$(X',M_{X'}) = (X,M_X) \times_{(Y,M_Y)} (Y',M_{Y'}), 
(D'_{[i]},M_{X'}^{\ve}|_{D'_{[i]}}) = 
(D_{[i]},M_{X}^{\ve}|_{D_{[i]}}) \times_{(Y,M_Y)} (Y',M_{Y'})$, and these 
immediately implies the assertions of the proposition. So we prove 
the above claim. \par 
For $I \subseteq \{1,\cdots,r\}$, 
let us put $D_I := \{t_i=0\,(i \in I)\}$ and put 
$D'_I := D_I \times_{Y} Y'$. 
Then, since $(D'_I,M_{X'}^{\ve}|_{D'_I})$ 
is log smooth over $(Y',M_{Y'})$, it is log regular. 
In particular, 
$\cO_{D'_I,\ol{x}'}/I(M_{X'}^{\ve}|_{D'_I},\ol{x}')$ is regular. 
So it suffices to prove the equality 
\begin{equation}\label{dim1}
\dim \cO_{D'_J,\ol{x}'}/I(M_{X'}^{\ve}|_{D'_J},\ol{x}') = 
\dim \cO_{D'_I,\ol{x}'}/I(M_{X'}^{\ve}|_{D'_I},\ol{x}') -1
\end{equation}
for any $J=\{j\}\cup I$ with $j \notin I$. 
Since $(D'_I,M_{X'}^{\ve}|_{D'_I}), (D'_J,M_{X'}^{\ve}|_{D'_J})$ are 
log regular, we have the equalities 
\begin{equation}\label{dim2}
\dim \cO_{D'_I,\ol{x}'} = 
\dim \cO_{D'_I,\ol{x}'}/I(M_{X'}^{\ve}|_{D'_I},\ol{x}') + 
\rk M_{X',\ol{x}'}^{\gp}/\cO_{X',\ol{x}'}^{\times}, 
\end{equation}
\begin{equation}\label{dim3}
\dim \cO_{D'_J,\ol{x}'} = 
\dim \cO_{D'_J,\ol{x}'}/I(M_{X'}^{\ve}|_{D'_J},\ol{x}') + 
\rk M_{X',\ol{x}'}^{\gp}/\cO_{X',\ol{x}'}^{\times}. 
\end{equation}
Now let us note the equality 
\begin{equation}\label{dim4}
\dim \cO_{D'_I,\ol{x}'} = \dim \cO_{D'_J,\ol{x}'} +1: 
\end{equation}
This follows from the equalities 
$\reldim (D'_I/Y')=\reldim (D_I/Y), \reldim (D'_J/\allowbreak Y')=
\reldim (D_J/Y)$, which are true because $D_I \lra Y, D_J \lra Y$ are 
smooth (in classical sense) on dense open subset. Combining 
\eqref{dim2}, \eqref{dim3} and \eqref{dim4}, we obtain \eqref{dim1}. 
So we are done. 
\end{pf}


\section{Alteration and hypercovering}

In this section, we prove the existence of certain diagrams 
involving hypercovering, which is a slight generalization of that 
treated in \cite[\S 6]{shiho4}. In this section, all the schemes are 
assumed to be defined over $S$ and from now on, we always assume that 
the field $k$ is perfect. \par 
First let us prove the following 
proposition, which is a slight generalization of \cite[6.4]{shiho3} 
and a consequence of the papers \cite{dejong1}, \cite{dejong2}: 

\begin{prop}\label{dejongup}
Let $f:X \lra Y$ be a proper morphism of integral schemes whose generic 
fiber is geometrically irreducible and let $D \subseteq X$ be a Cartier 
divisor. Then there exists a diagram 
\begin{equation}\label{alt1}
\begin{CD}
X @<{\psi}<< X' \\ 
@VfVV @VgVV \\ 
Y @<<< Y'
\end{CD}
\end{equation}
such that horizontal arrows are alterations with $X',Y'$ regular 
and a normal crossing divisor $D'$ $($resp. $E')$ on $X'$ $($resp. $Y')$ with 
a decomposition $D'={D'}^{\ho} \cup {D'}^{\ve}$ into sub normal 
crossing divisors ${D'}^{\ho}, {D'}^{\ve}$, satisfying the following 
conditions$:$ 
\begin{enumerate}
\item ${D'}^{\ve} = g^{-1}(E)_{\red}$. 
\item $\psi^{-1}(X-D)\cap (X'-{D'}^{\ve}) = (X'-{D'}^{\ho}) \cap 
(X'-{D'}^{\ve})$. 
\item 
If we denote the log structure on $X'$ $($resp. $Y')$ 
associated to $D'$ (resp. $E')$ by $M_{X'}$ $($resp. $M_{Y'})$, 
the morphism $(X',M_{X'}) \lra (Y',M_{Y'})$ induced by $g$ is 
proper log smooth integral and universally saturated. 
\item 
If we denote the log structure on $X'$ associated to ${D'}^{\ve}$ by 
$M^{\ve}_{X'}$ and define 
$D'_{[i]} \,(i \geq 0)$ by 
\begin{align*}
D'_{[0]} & :=X', \qquad D'_{[1]} := \text{the normalization of 
${D'}^{\ho}$}, \\
D'_{[i]} & := \text{$i$-fold fiber product of $D'_{[1]}$ over $X'$}, 
\end{align*}
the morphism $(D'_{[i]},M_{X'}^{\ve}|_{D'_{[i]}}) \lra (Y',M_{Y'})$ 
induced by $g$ is proper log smooth integral and universally saturated 
for any $i \in \N$. 
\end{enumerate}
\end{prop} 

\begin{pf} 
By using \cite[Thm 5.9]{dejong2} and \cite[Prop 5.11]{dejong2} and 
argueing as \cite[6.4]{shiho4}, we see that there exist 
a diagram \eqref{alt1} and a normal crossing divisor 
$D''$ $($resp. $E')$ on $X'$ $($resp. $Y')$ satisfying the following 
conditions: 
\begin{enumerate}
\item 
$\psi^{-1}(D)_{\red} \subseteq D''$. 
\item 
${D'}^{\ve} := g^{-1}(E)_{\red}$ is a sub normal crossing divisor of $D''$. 
\item 
For any sub normal crossing divisor $D''' \subseteq D''$ containing 
${D'}^{\ve}$, 
the condition (3) in the statement of the proposition is true 
if we replace $D'$ by $D'''$. 
\item 
For any sub normal crossing divisor $D''' \subseteq D''$ containing 
${D'}^{\ve}$, 
the condition (4) in the statement of the proposition is true 
if we replace ${D'}^{\ho}$ by 
${D'''}^{\ho} = \text{closure of ${D'''}-{D'}^{\ve}$}$. 
\end{enumerate}
Now let us define $D'$ as $(\psi^{-1}(D))_{\red} \cup {D'}^{\ve}$ and 
let ${D'}^{\ho}$ be the closure of $D'-{D'}^{\ve}$. 
Then we see that all the required conditions are satisfied. 
\end{pf} 

Next we prove the following lemma, which is a slight generalization of 
\cite[6.9]{shiho4}: 

\begin{lem}\label{hc1}
Let $f:\ol{X} \lra \ol{Y}$ be a proper morphism of schemes, let 
$U \subseteq Y$ be a dense open subscheme and 
let $X \subseteq \ol{X}$ be an open subscheme. Then we have the diagram 
consisting of strict morphism of pairs 
\begin{equation}\label{*}
\begin{CD}
(U_{\ol{X}},\ol{X}) @<{g_X}<< (U_{\ol{X}'},\ol{X}') \\ 
@VfVV @V{f'}VV \\ 
(U_{\ol{Y}},\ol{Y}) @<{g_Y}<< (U_{\ol{Y}'},\ol{Y}')
\end{CD}
\end{equation}
with $(U_{\ol{Y}'},\ol{Y}') = \coprod_{j=1}^b (U_{\ol{Y}'_j}, \ol{Y}'_j)$, 
$(U_{\ol{X}'},\ol{X}') = 
\coprod_{j=1}^b(\coprod_{i=1}^{r_j} (U_{\ol{X}'_{ji}}, \ol{X}'_{ji}))$ 
$($decomposition into connected components with 
$f'(\ol{X}'_{ji}) \subseteq \ol{Y}'_j)$ 
satisfying the following conditions$:$ 
\begin{enumerate} 
\item $U_{\ol{Y}}$ is contained in $U$ and dense in $\ol{Y}$. 
\item $\ol{X}'_{ji}$'s, $\ol{Y}'_j$'s are integral, and 
if we put $X' := X \times_{\ol{X}} \ol{X}'$, 
$\ol{X}'_{ji} - X' \subseteq \ol{X}'_{ji}$ is a Cartier divisor $($possibly 
empty$)$. 
\item 
The map $\ol{g}_X: (X' \cap U_{\ol{X}'},\ol{X}') \lra 
((X \cap U_{\ol{X}}) \times_{U_{\ol{Y}}} U_{\ol{Y}'}, \ol{X} \times_{\ol{Y}} 
\ol{Y}')$ 
induced by $(g_X,f')$ 
is a proper covering and the map $g_Y$ is a good 
strongly proper covering. 
\item 
For each $i,j$, the generic fiber of 
$f'|_{\ol{X}'_{ji}}: \ol{X}'_{ji} \lra \ol{Y}'_j$ is 
non-empty and geometrically 
irreducible. $($Attention$:$ It is possible that $r_j=0$ holds for some $j.)$
\end{enumerate} 
\end{lem} 

\begin{pf} 
Let $\ol{X}_{\red} = \bigcup_l \ol{X}_l$ be the decomposition of 
$\ol{X}_{\red}$ into irreducible components and put $X_l := X \cap \ol{X}_l$. 
For $l$ with $X_l \not= \emptyset$, let us take an alteration 
$\ol{X}_{l} \lla \ol{X}''_l$ so that $\ol{X}''_l$ is smooth over $S$ 
and that the complement of $X''_l := X_l \times_{\ol{X}_l} \ol{X}''_l$ 
is a simple normal crossing divisor in $\ol{X}''_l$. Let us put 
$\ol{X}'' := \coprod_{l,X_l \not= \emptyset} \ol{X}''_l, 
X'' := X \times_{\ol{X}} \ol{X}''$. Then the natural morphism of pairs 
$(X'',\ol{X}'') \lra (X,\ol{X})$ is a proper covering. Then, 
by using \cite[6.9]{shiho4} for the morphism $\ol{X}'' \lra \ol{Y}$, 
we can take a diagram consisting of strict morphisms of pairs 
\begin{equation}\label{hyper1}
\begin{CD}
(U_{\ol{X}},\ol{X}) @<<< (U_{\ol{X}''},\ol{X}'') @<<< (U_{\ol{X}'},\ol{X}') \\ 
@VVV @VVV @VVV \\ 
(U_{\ol{Y}},\ol{Y}) @= (U_{\ol{Y}},\ol{Y}) @<<< (U_{\ol{Y}'},\ol{Y}') 
\end{CD}
\end{equation} 
such that the right square satisfies the conclusion of \cite[6.9]{shiho4}. 
Then one can see that the diagram \eqref{*} induced by \eqref{hyper1} 
satisfies the required conditions. 
\end{pf} 

Next we prove the following proposition, which is a generalization of 
\cite[6.10]{shiho4}: 

\begin{prop}\label{hc2}
Let $f:\ol{X} \lra \ol{Y}$ be a proper morphism of schemes, let 
$U \subseteq \ol{Y}$ be a dense open subscheme and let $X \subseteq \ol{X}$ 
be an open subscheme. Then we have the diagram 
consisting of strict morphism of pairs 
\begin{equation}\label{**}
\begin{CD}
(U_{\ol{X}},\ol{X}) @<{g_X}<< (U_{\tiol{X}},\tiol{X}) \\ 
@VfVV @V{\ti{f}}VV \\ 
(U_{\ol{Y}},\ol{Y}) @<{g_Y}<< (U_{\tiol{Y}},\tiol{Y})
\end{CD}
\end{equation}
satisfying the following conditions$:$ 
\begin{enumerate} 
\item $U_{\ol{Y}}$ is contained in $U$ and dense in $\ol{Y}$. 
\item $\tiol{Y}$ is regular. 
\item 
The map $\ol{g}_X: (\ti{X} \cap U_{\tiol{X}},\tiol{X}) 
\lra ((X \cap U_{\ol{X}}) \times_{U_{\ol{Y}}} U_{\tiol{Y}}, 
\ol{X} \times_{\ol{Y}} \tiol{Y})$ 
induced by $(g_X,\ti{f})$ 
is a proper covering $($where we put 
$\ti{X} := X \times_{\ol{X}} \tiol{X})$ 
and the map $g_Y$ is 
a good strongly proper covering. 
\item 
There exist fs log structures $M^{\ve}_{\tiol{X}}, M_{\tiol{Y}}$ on 
$\tiol{X}, \tiol{Y}$ 
respectively such that $M_{\tiol{Y}}$ is associated to a simple normal 
crossing divisor on $\tiol{Y}$, 
$U_{\tiol{X}} \subseteq (\tiol{X},M^{\ve}_{\tiol{X}})_{\triv}, 
U_{\tiol{Y}} \subseteq (\tiol{Y},M_{\tiol{Y}})_{\triv}$ holds and that 
there exists a proper log smooth integral universally saturated morphism 
$(\tiol{X},M^{\ve}_{\tiol{X}}) \lra (\tiol{Y},M_{\tiol{Y}})$ 
whose underlying morphism of schemes is the same as $\ti{f}$. 
\item 
There exists a log normal crossing divisor $D \subseteq \tiol{X}$ such that, 
if we define $M_{\tiol{X}}, D_{[i]} \,(i \geq 0)$ as in Proposition 
\ref{lncdprop} using $D$, the induced morphisms 
$$ (\tiol{X},M_{\tiol{X}}) \lra (\tiol{Y},M_{\tiol{Y}}), 
\,\,\,\, 
(D_{[i]},M^{\ve}_{\tiol{X}}|_{D_{[i]}}) \lra 
(\tiol{Y},M_{\tiol{Y}}) \, (i \geq 0) $$
are proper log smooth integral universally saturated, and we have 
$\ti{X} \cap (\tiol{X}, \allowbreak M^{\ve}_{\tiol{X}})_{\triv} \allowbreak = 
\allowbreak (\tiol{X}-D) \cap (\tiol{X},M^{\ve}_{\tiol{X}})_{\triv}$. 
$($In particular, we have 
$\ti{X} \cap U_{\tiol{X}} = (\tiol{X}-D) \cap U_{\tiol{X}}.)$
\end{enumerate}
\end{prop} 

\begin{pf} 
The proof is similar to \cite[6.10]{shiho4}, although the proof here 
is slightly more complicated. 
First let us take the diagram 
consisting of strict morphism of pairs 
\begin{equation*}
\begin{CD}
(V_{\ol{X}},\ol{X}) @<<< (U_{\ol{X}'},\ol{X}') \\ 
@VfVV @V{f'}VV \\ 
(V_{\ol{Y}},\ol{Y}) @<<< (U_{\ol{Y}'},\ol{Y}')
\end{CD}
\end{equation*}
satisfying the conclusion of Lemma \ref{hc1}, and let 
$\ol{Y}' = \coprod_{j=1}^b \ol{Y}'_j$, $\ol{X}' := 
\coprod_{j=1}^b(\coprod_{i=1}^{r_j} \ol{X}'_{ji})$ 
be the decomposition of $\ol{Y}',\ol{X}'$ 
into connected components with $f'(\ol{X}'_{ji}) \subseteq \ol{Y}'_j$. 
Also let us put $X' := X \times_{\ol{X}} \ol{X}'$. \par 
Let us fix $j$. We prove the following claim: \\
\quad \\
{\bf claim}\,\,\, 
For any $1 \leq s \leq r_j$, there exists a diagram 
consisting of strict morphism of pairs 
\begin{equation}\label{altdiag}
\begin{CD}
\coprod_{i=1}^{r_j}(\ol{X}'_{ji},O'_{\ol{X},ji}) 
@<<< \coprod_{i=1}^{r_j} (\ol{X}''_{ji},O''_{\ol{X},ji}) \\ 
@V{f'}VV @V{f''}VV \\
(\ol{Y}'_j,O'_{\ol{Y},j}) @<<< (\ol{Y}''_j,O''_{\ol{Y},j}) 
\end{CD}
\end{equation} 
(we put $X''_{ji} := X' \times_{\ol{X}'} \ol{X}''_{ji}$)
satisfying the following conditions: 
\begin{enumerate}
\item 
$O'_{\ol{Y},j}$ is dense in $\ol{Y}'_j$. 
\item 
$\ol{Y}''_j$ is regular and connected (so it is integral). 
\item 
The upper horizontal arrow is a proper covering 
and the lower horizontal arrow is a good strongly proper covering. 
\item 
There exist 
fs log structures $M^{\ve}_{\coprod_{i=1}^s\ol{X}''_{ji}}, 
M_{\ol{Y}''_j}$ on $\coprod_{i=1}^s\ol{X}''_{ji}, \ol{Y}''_j$ respectively and 
a proper log smooth integral universally saturated morphism 
$$ g:(\coprod_{i=1}^s\ol{X}''_{ji}, 
M^{\ve}_{\coprod_{i=1}^s\ol{X}''_{ji}}) \lra (\ol{Y}''_j, M_{\ol{Y}''_j})$$ 
whose underlying morphism of schemes is the same as 
$f''|_{\coprod_{i=1}^s \ol{X}''_{ji}}$ 
such that $M_{\ol{Y}''_j}$ is associated to a simple normal crossing 
divisor on $\ol{Y}''_j$ and that 
$g^{-1}((\ol{Y}''_j, \allowbreak 
M_{\ol{Y}''_j})_{\triv})$ is contained in 
$(\coprod_{i=1}^s\ol{X}''_{ji}, 
M^{\ve}_{\coprod_{i=1}^s\ol{X}''_{ji}})_{\triv}$. 
\item 
There exists a log normal crossing divisor $D \subseteq \coprod_{i=1}^s
\ol{X}''_{ji}$ such that, 
if we define $M_{\coprod_{i=1}^s\ol{X}''_{ji}}, 
D_{[i]} \,(i \geq 0)$ as in Proposition 
\ref{lncdprop} using $D$, the induced morphisms 
$$ (\coprod_{i=1}^s\ol{X}''_{ji}, 
M_{\coprod_{i=1}^s\ol{X}''_{ji}}) \lra (\ol{Y}''_j,M_{\ol{Y}''_j}), 
\,\,\,\, 
(D_{[i]},M^{\ve}_{\coprod_{i=1}^s\ol{X}''_{ji}}|_{D_{[i]}}) 
\lra (\ol{Y}''_j,M_{\ol{Y}''_j}) \, (i \geq 0) $$
are proper log smooth integral universally saturated, and we have 
$(\coprod_{i=1}^sX''_{ji}) \cap (\coprod_{i=1}^s\ol{X}''_{ji},
M^{\ve}_{\coprod_{i=1}^s\ti{X}''_{ji}})_{\triv} = 
(\coprod_{i=1}^s\ol{X}''_{ji}-D) \cap 
(\coprod_{i=1}^s\ol{X}''_{ji},
M^{\ve}_{\coprod_{i=1}^s\ti{X}''_{ji}})_{\triv}$. 
\item 
For each $i \geq s+1$, $\ol{X}''_{ji}$ is integral, $\ol{X}''_{ji}-X''_{ji}$ 
is a Cartier divisor and 
the generic fiber of 
$f''|_{\ol{X}''_{ji}}: \ol{X}''_{ji} \lra \ol{Y}''_j$ is 
non-empty, geometrically irreducible. 
\end{enumerate}

We prove the claim by induction on $s$: In the case $s=1$, 
we apply Proposition \ref{dejongup} to the morphism 
$\ol{X}'_{j1} \lra \ol{Y}'_j$ to obtain the diagram 
\begin{equation*}
\begin{CD}
\ol{X}'_{j1} @<<< \ol{X}''_{j1} \\ 
@VVV @VVV \\ 
\ol{Y}'_j @<<< \ol{Y}''_j 
\end{CD}
\end{equation*}
satisfying the conclusion of Proposition \ref{dejongup}. 
Then there exist log structures $M^{\ve}_{\ol{X}''_{j1}}, \allowbreak 
M_{\ol{Y}''_j}$ 
on $\ol{X}''_{j1}, \ol{Y}''_j$ respectively and a log normal crossing 
divisor $D$ in $\ol{X}''_{ji}$ satisfying (2), (4) and (5) (for $s=1$). 
Then, by the argument of the proof of 
the claim in \cite[6.10]{shiho4}, we can take an open subscheme 
$O'_{Y,j} \subseteq \ol{Y}'_j$ and the diagram \eqref{altdiag} 
satisfying the other assertions. So the claim is true for $s=1$. \par 
If the claim is true for $s-1$, we have the diagram 
consisting of strict morphisms of pairs 
\begin{equation}
\begin{CD}
\coprod_{i=1}^{r_j} (\ol{X}'_{ji},O'_{\ol{X},ji}) @<<< 
\coprod_{i=1}^{r_j} (\ol{X}''_{ji},O''_{\ol{X},ji}) \\ 
@V{f'}VV @V{f''}VV \\
(\ol{Y}'_j,O'_{\ol{Y},j}) @<<< (\ol{Y}''_j,O''_{\ol{Y},j}), 
\end{CD}
\end{equation} 
fs log structures $M^{\ve}_{\coprod_{i=1}^{s-1}\ol{X}''_{ji}}, 
M_{\ol{Y}''_j}$ on $\coprod_{i=1}^{s-1}\ol{X}''_{ji}, \ol{Y}''_j$ 
respectively,  
a proper log smooth integral universally saturated morphism 
$g:(\coprod_{i=1}^{s-1}\ol{X}''_{ji}, 
M^{\ve}_{\coprod_{i=1}^{s-1}\ol{X}''_{ji}}) \lra (\ol{Y}''_j, M_{\ol{Y}''_j})$ 
 and a log normal crossing divisor 
$D_1 \subseteq \coprod_{i=1}^{s-1}\ol{X}''_{ji}$ satisfying the properties 
stated in the claim (for $s-1$). 
Let us apply Proposition \ref{dejongup} to the morphism 
$\ol{X}''_{js} \lra \ol{Y}''_j$ to obtain the diagram 
\begin{equation*}
\begin{CD}
\ol{X}''_{js} @<<< \ol{X}'''_{js} \\ 
@VVV @VVV \\ 
\ol{Y}''_j @<{\varphi}<< \ol{Y}'''_j, 
\end{CD}
\end{equation*}
log structures $M^{\ve}_{\ol{X}'''_{js}}, M_{\ol{Y}'''_j}$ on 
$\ol{X}'''_{js}, \ol{Y}'''_j$ respectively and a log normal 
crossing divisor $D_2 \subseteq \ol{X}'''_{js}$ which  
satisfy the analogue of the conditions (4), (5) for `$s$-component'. 
Let us define the reduced closed subscheme $Z \subseteq \ol{Y}'''_j$ 
as the scheme whose underlying space is 
the union of the simple normal crossing divisor corresponding to 
$M_{\ol{Y}'''_j}$ and the inverse image of the simple normal crossing 
divisor corresponding to $M_{\ol{Y}''_j}$, and let us take 
an alteration $\psi: \ol{Y}''''_j \lra \ol{Y}'''_j$ such that 
$\ol{Y}''''_j$ is smooth over $S$ and that $\psi^{-1}(Z)_{\red}$ is a 
simple normal crossing divisor on $\ol{Y}''''_j$. Then, if we denote 
the log structure associated to $\psi^{-1}(Z)_{\red}$ by $M_{\ol{Y}''''_j}$, 
the morphisms $\varphi \circ \psi, \psi$ induce the morphisms of log schemes 
$\varphi \circ \psi: (\ol{Y}''''_j,M_{\ol{Y}''''_j}) \lra 
(\ol{Y}'',M_{\ol{Y}''_j}), 
\psi: (\ol{Y}''''_j,M_{\ol{Y}''''_j}) \lra (\ol{Y}'''_j,M_{\ol{Y}'''_j})$. 
Now let us put 
$$ 
(\coprod_{j=1}^{s-1}\ol{X}''''_{ji},
M^{\ve}_{\coprod_{j=1}^{s-1}\ol{X}''''_{ji}}) 
:= 
(\coprod_{j=1}^{s-1}\ol{X}''_{ji},M^{\ve}_{\coprod_{j=1}^{s-1}\ol{X}''_{ji}}) 
\times_{(\ol{Y}''_j,M_{\ol{Y}''_j})} (\ol{Y}''''_j,M_{\ol{Y}''''_j}), $$ 
$$ 
(\ol{X}''''_{js},
M^{\ve}_{\ol{X}''''_{js}}) 
:= 
(\ol{X}'''_{js},M^{\ve}_{\ol{X}'''_{js}}) 
\times_{(\ol{Y}'''_j,M_{\ol{Y}'''_j})} (\ol{Y}''''_j,M_{\ol{Y}''''_j}), $$ 
$$ 
(\coprod_{j=1}^{s}\ol{X}''''_{ji},
M^{\ve}_{\coprod_{j=1}^{s}\ol{X}''''_{ji}}) := 
(\coprod_{j=1}^{s-1}\ol{X}''''_{ji},
M^{\ve}_{\coprod_{j=1}^{s-1}\ol{X}''''_{ji}}) \coprod 
(\ol{X}''''_{js},
M^{\ve}_{\ol{X}''''_{js}}) $$ 
and let $D \subseteq \coprod_{j=1}^s\ol{X}''''_{ji}$ be the disjoint 
union of the inverse images of $D_1, D_2$ to 
$\coprod_{j=1}^s\ol{X}''''_{ji}$. 
Then these data satisfy the analogue of (2), (4) and (5). 
(To see that $D$ is a log normal crossing divisor, 
we use Proposition \ref{lncdpullback}.) 
Then, by the argument in the proof of the claim in \cite[6.10]{shiho4}, 
we see that we can form the diagram like \eqref{altdiag} for $s$. So 
the proof of the claim is finished. \par 
By the claim for $s=r_j$ (for each $j$), we see that there exists a diagram 
consisting of strict morphism of pairs 
\begin{equation*}
\begin{CD}
(\ol{X}',O'_{\ol{X}}) @<<< (\tiol{X},O_{\tiol{X}}) \\ 
@V{f'}VV @V{\ti{f}}VV \\ 
(\ol{Y}',O'_{\ol{Y}}) @<<< (\tiol{Y},O_{\tiol{Y}}) 
\end{CD}
\end{equation*}
(we put $\ti{X} := X' \times_{\ol{X}'} \tiol{X}$)
satisfying the following conditions: 
\begin{enumerate}
\item 
$\ol{O}'_Y \subseteq \ol{Y}'$ is dense open. 
\item 
$\tiol{Y}$ is regular. 
\item 
The upper horizontal arrow is a proper covering 
and the lower horizontal arrow is a good strongly proper covering. 
\item  
There exist fs log structures $M^{\ve}_{\tiol{X}}, M_{\tiol{Y}}$ on 
$\tiol{X}, \tiol{Y}$ respectively and 
a proper log smooth integral universally saturated morphism 
$(\tiol{X},M^{\ve}_{\tiol{X}}) \lra (\tiol{Y},M_{\tiol{Y}})$ 
whose underlying morphism of schemes is $\ti{f}$ 
such that $M_{\tiol{Y}}$ is associated to a simple normal crossing 
divisor on $\tiol{Y}$ and that 
$\ti{f}^{-1}((\tiol{Y}, M_{\tiol{Y}})_{\triv})$ is contained in 
$(\tiol{X}, M_{\tiol{X}})_{\triv}$. 
\item 
There exists a log normal crossing divisor $D \subseteq \tiol{X}$ such that, 
if we define $M_{\tiol{X}}, 
D_{[i]} \,(i \geq 0)$ as in Proposition 
\ref{lncdprop} using $D$, the induced morphisms 
$$ (\tiol{X}, M_{\tiol{X}}) \lra (\tiol{Y},M_{\tiol{Y}}), 
\,\,\,\, 
(D_{[i]},M^{\ve}_{\tiol{X}}|_{D_{[i]}}) 
\lra (\tiol{Y},M_{\tiol{Y}}) \,\, (i \geq 0) $$
are proper log smooth integral universally saturated, and we have 
the equality 
$\ti{X} \cap (\tiol{X},M^{\ve}_{\ti{X}})_{\triv} = 
(\tiol{X}-D) \cap (\tiol{X}, M^{\ve}_{\tiol{X}})_{\triv}$. 
\end{enumerate}
Then, by the argument in the proof of \cite[6.10]{shiho4}, we can 
take an open subset $U_{\ol{Y}} \subseteq U$ and the diagram 
\eqref{**} in order that all the required conditions are satisfied. 
So we are done. 
\end{pf} 

Using Proposition \ref{hc2}, we can prove the following theorem: 

\begin{thm}\label{hc3}
Let $f:\ol{X} \lra \ol{Y}$ be a proper morphism of schemes, 
let $X \subseteq \ol{X}$ be an open subscheme and let $q \in \N$. 
Then we have the diagram 
consisting of strict morphism of pairs 
\begin{equation}\label{***}
\begin{CD}
(U_{\ol{X}},\ol{X}) @<<< (U_{\tiol{X}},\tiol{X}) @<h<< 
(U_{\tiol{X}^{(\b)}}, \tiol{X}^{(\b)})\\ 
@VfVV @V{\ti{f}}VV @V{\ti{f}^{(\b)}}VV \\ 
(U_{\ol{Y}},\ol{Y}) @<{g_Y}<< (U_{\tiol{Y}},\tiol{Y}) @= (U_{\tiol{Y}},
\tiol{Y})
\end{CD}
\end{equation}
$($we put $\ti{X} := X \times_{\ol{X}} \tiol{X}, 
\ti{X}^{(\b)}:=X \times_{\ol{X}} \tiol{X}^{(\b)})$ 
satisfying the following conditions$:$ 
\begin{enumerate} 
\item $U_{\ol{Y}}$ is dense in $\tiol{Y}$. 
\item $\tiol{Y}$ is regular. 
\item 
The left square is Cartesian, $g_Y$ is a good strongly proper covering 
and the morphism 
$$ (\ti{X} \cap U_{\tiol{X}}, \tiol{X}) \lla 
(\ti{X}^{(\b)} \cap U_{\tiol{X}^{(\b)}}, \tiol{X}^{(\b)}) $$ 
induced by $h$ is a $q$-truncated proper hypercovering by a 
$q$-truncated split simplicial pair. 
\item 
There exist fs log structures 
$M^{\ve}_{\tiol{X}^{(n)}}, M_{\tiol{Y}}$ on $\tiol{X}^{(n)}, \tiol{Y}$ 
respectively $(n \leq q)$ 
such that $M_{\tiol{Y}}$ is associated to a simple normal 
crossing divisor on $\tiol{Y}$, 
$U_{\tiol{X}^{(n)}} \subseteq 
(\tiol{X}^{(n)},M^{\ve}_{\tiol{X}^{(n)}})_{\triv}, 
U_{\tiol{Y}} \subseteq (\tiol{Y},M_{\tiol{Y}})_{\triv}$ holds and that 
for each $n \leq q$, there exists a proper log smooth integral 
universally saturated morphism 
$(\tiol{X}^{(n)},M^{\ve}_{\tiol{X}^{(n)}}) \lra (\tiol{Y},M_{\tiol{Y}})$ 
whose underlying morphism of schemes is the same as $\ti{f}^{(n)}$. 
\item 
For each $n \leq q$, 
there exists a log normal crossing divisor $D^{(n)} 
\subseteq \tiol{X}^{(n)}$ such that, 
if we define $M_{\tiol{X}^{(n)}}, D^{(n)}_{[i]} 
\,(i \geq 0)$ as in Proposition 
\ref{lncdprop} using $D^{(n)}$, the induced morphisms 
$$ (\tiol{X}^{(n)},M_{\tiol{X}^{(n)}}) \lra (\tiol{Y},M_{\tiol{Y}}), 
\,\,\,\, 
(D^{(n)}_{[i]},M^{\ve}_{\tiol{X}^{(n)}}|_{D^{(n)}_{[i]}}) 
\lra (\tiol{Y},M_{\tiol{Y}}) \, (i \geq 0) $$
are proper log smooth integral universally saturated, and we have 
$\ti{X}^{(n)} \cap (\tiol{X}^{(n)},M^{\ve}_{\tiol{X}^{(n)}})_{\triv} = 
(\tiol{X}^{(n)}-D^{(n)}) \cap (\tiol{X}^{(n)},
M^{\ve}_{\tiol{X}^{(n)}})_{\triv}$. 
$($In particular, we have 
$\ti{X}^{(n)} \cap U_{\tiol{X}^{(n)}} = (\tiol{X}^{(n)}-D^{(n)}) 
\cap U_{\tiol{X}^{(n)}}.)$
\end{enumerate}
\end{thm} 

\begin{pf} 
We can prove the theorem in the same way as \cite[6.11]{shiho4}, 
by using Proposition \ref{hc2} instead of \cite[6.10]{shiho4}. 
The detail is left to the reader. 
\end{pf}


\section{Generic overconvergence} 

In this section, we prove the main theorem of this paper, 
that is, the generic overconvergence of relative rigid 
cohomology for any morphism when the coefficient is a potentially semistable 
overconvergent isocrystal. In particular, we obtain the 
generic overconvergence of relative rigid cohomology when the 
coefficient is an overconvergent $F$-isocrystal, by using 
Theorem \ref{sconj}. \par 
Note that we keep the assumption that $k$ is perfect also in this section. 
The precise statement of the main theorem is as follows: 

\begin{thm}\label{mainthm}
Let us assume given a morphism of pairs 
$f: (X,\ol{X}) \lra (Y,\ol{Y})$ such that $f:\ol{X} \lra \ol{Y}$ 
is proper. 
Then, there exists a dense open set 
$U_{\ol{Y}}$ of $Y$, a proper surjective map $\tiol{Y} \lra \ol{Y}$ 
and a simple normal crossing 
divisor $E$ on $\tiol{Y}$ with $E \cap U_{\tiol{Y}} = \emptyset$ 
$($where we put $U_{\tiol{Y}} := U_{\ol{Y}} \times_{\ol{Y}} 
\tiol{Y})$ such that, 
for any potentially semi-stable overconvergent isocrystal $\cE$ on 
$(X,\ol{X})/\cS_K$ and $q \in \N$, there exists 
the unique overconvergent isocrystal $\cF$ on $(U_Y,\ol{Y})/\cS_K$ 
satisfying the following condition$:$ 
For any $(U_{\tiol{Y}}, \tiol{Y})$-triple $(Z,\ol{Z},\cZ)$ over 
$(S,S,\cS)$ such that $\ol{Z}$ is smooth over $S$, $(E \times_{\tiol{Y}} 
\ol{Z})_{\red}$ is contained in a simple normal crossing divisor on 
$\ol{Z}$ and that $\cZ$ is formally smooth over $\cS$, 
the restriction of $\cF$ to $I^{\d}((Z,\ol{Z})/\cS_K, \cZ)$ 
is given functorially by 
$(R^qf_{(X \times_{Y} Z, \ol{X} \times_{\ol{Y}} \ol{Z})/\cZ, \rig *}\cE, 
\epsilon)$, where $\epsilon$ is an isomorphism 
{\small{ $$ 
p_2^*R^qf_{(X \times_{Y} Z, 
\ol{X} \times_{\ol{Y}} \ol{Z})/\cZ, \rig *}\cE 
\os{\sim}{\ra} 
R^qf_{(X \times_{Y} Z, \ol{X} \times_{\ol{Y}} \ol{Z})/\cZ \times_{\cS} \cZ, 
\rig *}\cE \os{\sim}{\leftarrow}
p_1^*R^qf_{(X 
\times_{Y} Z, \ol{X} \times_{\ol{Y}} \ol{Z})/\cZ, \rig *}\cE
$$ }}
$(p_i$ denotes the $i$-th projection 
$]\ol{Z}[_{\cZ \times_{\cS} \cZ} \lra ]\ol{Z}[_{\cZ})$. 
\end{thm} 

\begin{pf} 
First, by the argument in the beginning of the proof of \cite[7.4]{shiho4}, 
there exists a positive integer $q_0$ such that, for any overconvergent 
isocrystal $\cE$ on $(X,\ol{X})/\cS_K$ and any $(Y,\ol{Y})$-triple 
$(Z,\ol{Z},\cZ)$ over $(S,S,\cS)$, we have 
\begin{equation*}\label{bound-rig}
R^qf_{(X \times_Y Z,\ol{X} \times_{\ol{Y}} \ol{Z})/\cZ, \rig *}\cE = 0
\end{equation*} 
for any $q > q_0$. \par 
Let us take a proper hypercovering $(X^{(\b)}, \ol{X}^{(\b)}) 
\lra (X,\ol{X})$ such that each $\ol{X}^{(n)}$ is smooth over $k$ 
and, each $\ol{X}^{(n)}-X^{(n)}$ is a simple normal crossing divisor on 
$\ol{X}^{(n)}$ and that the restriction of $\cE$ to 
$I^{\d}((X^{(n)},\ol{X}^{(n)})/\cS_K)$ extends to a locally free 
isocrystal $\ol{\cE}$ on 
$(\ol{X}^{(n)}/\cS)^{\log}_{\conv}=((\ol{X}^{(n)},
M_{\ol{X}^{(n)}})/\cS)_{\conv}$ having nilpotent residues, 
where $M_{\ol{X}^{(n)}}$ denotes the 
log structure on $\ol{X}^{(n)}$ associated to $\ol{X}^{(n)}-X^{(n)}$. 
(This is possible by \cite{dejong1} and the potential semi-stability of 
$\cE$). \par 
Next let us put $q_1:=q_0(q_0+1)/2, q_2:=q_1(q_1+1)/2$. 
By applying Theorem \ref{hc3} to the proper map 
$\coprod_{n=0}^{q_1}\ol{X}^{(n)} \lra \ol{Y}$ and the open subscheme 
$\coprod_{n=0}^{q_1}X^{(n)}$ of $\coprod_{n=0}^{q_1}\ol{X}^{(n)}$ 
for $q=q_2$, we have the following: For each $0 \leq n \leq q_1$, 
we have the following diagram consisting of strict morphisms of pairs 
\begin{equation}\label{hc3diag}
\begin{CD}
(U_{\ol{X}^{(n)}}, \ol{X}^{(n)}) @<<< 
(U_{\tiol{X}^{(n)}}, \tiol{X}^{(n)}) @<h<< 
(U_{\tiol{X}^{(n)\lr}}, \tiol{X}^{(n)\lr}) \\ 
@VVV @VVV @V{\ti{f}^{(n)\lr}}VV \\ 
(U_{\ol{Y}},\ol{Y}) @<<< (U_{\tiol{Y}},\tiol{Y}) @= 
(U_{\tiol{Y}},\tiol{Y}), \\ 
\end{CD}
\end{equation}
which satisfies the conclusion of Theorem \ref{hc3} for $q=q_2$ and that 
the lower horizontal line in \eqref{hc3diag} is independent of 
$n$. In particular, there exist fs log structures 
$M^{\ve}_{\tiol{X}^{(n)\lmr}}, M_{\tiol{Y}}$ on 
$\tiol{X}^{(n)\lmr}, \tiol{Y}$ respectively such that 
$M_{\tiol{Y}}$ is assciated to a simple normal crossing divisor, 
$U_{\tiol{X}^{(n)\lmr}} \subseteq (\tiol{X}^{(n)\lmr},
M^{\ve}_{\tiol{X}^{(n)\lmr}})_{\triv}$ and 
$U_{\tiol{Y}} \subseteq (\tiol{Y},M_{\tiol{Y}})_{\triv}$ hold, 
and that there exists a proper log smooth integral universally 
saturated morphism 
$(\tiol{X}^{(n)\lmr},M^{\ve}_{\tiol{X}^{(n)\lmr}}) \lra 
(\tiol{Y},M_{\tiol{Y}})$ 
whose underlying morphism of schemes is the same as $\ti{f}^{(n)\lmr}$ 
(which we also denote by $\ti{f}^{(n)\lmr}$). 
Moreover, there exists a log normal crossing divisor 
$D^{(n)\lmr}$ in $\tiol{X}^{(n)\lmr}$ such that, if we define 
$M_{\tiol{X}^{(n)\lmr}}, D^{(n)\lmr}_{[i]}$ as in Proposition \ref{lncdprop}, 
the morphisms 
$$ 
(\tiol{X}^{(n)\lmr},M_{\tiol{X}^{(n)\lmr}}) \lra 
(\tiol{Y},M_{\tiol{Y}}), $$ 
$$ 
(D^{(n)\lmr}_{[i]}, M^{\ve}_{\tiol{X}^{(n)\lmr}}|_{D^{(n)\lmr}_{[i]}}) \lra 
(\tiol{Y},M_{\tiol{Y}}) $$ 
induced by $\ti{f}^{(n)\lmr}$ are proper log smooth integral universally 
saturated, and that we have the equality 
$\ti{X}^{(n)\lmr} \cap (\tiol{X},M^{\ve}_{\tiol{X}^{(n)\lmr}})_{\triv} 
= (\tiol{X}^{(n)\lmr} - D^{(n)\lmr}) \cap 
(\tiol{X},M^{\ve}_{\tiol{X}^{(n)\lmr}})_{\triv}$, 
where we denoted the inverse image of $X^{(n)}$ in 
$\tiol{X}^{(n)\lmr}$ by $\ti{X}^{(n)\lmr}$. \par 
Note that the properties in Theorem \ref{hc3} for the diagram 
\eqref{hc3diag} remains true if we take a morphism of log schemes 
$\varphi: (\tiol{Y}',M_{\tiol{Y}'}) \lra (\tiol{Y},M_{\tiol{Y}})$ such that 
$\tiol{Y}'$ is smooth over $k$ and that $M_{\tiol{Y}'}$ is a simple normal 
crossing divisor in $\tiol{Y}'$ which does not meet 
$U_{\tiol{Y}'} := \varphi^{-1}(U_{\tiol{Y}})$, and replace the right 
square of the diagram \eqref{hc3diag} by the pull-back of it by 
$\varphi:(U_{\tiol{Y}'}, \tiol{Y}') \lra (U_{\tiol{Y}},\tiol{Y})$. 
(We use Proposition \ref{lncdpullback}.) We can take 
$\varphi$ to be the alteration in order that $\tiol{Y}'$ is 
quasi-projective smooth over $S$, 
$\tiol{Y}'-U_{\tiol{Y}'}$ is a simple normal crossing divisor and 
$M_{\tiol{Y}'}$ is the log structure associated to 
$\tiol{Y}'-U_{\tiol{Y}'}$. So, in the diagram \eqref{hc3diag}, 
we may assume that $\tiol{Y}$ is quasi-projective smooth over $S$, 
$E := \tiol{Y}-U_{\tiol{Y}}$ is a simple normal crossing 
divisor and the log structure $M_{\tiol{Y}}$ is the one associated to 
$\tiol{Y}-U_{\tiol{Y}}$. \par 
Next we replace $U_{\ol{Y}}$ by $U_{\ol{Y}} \cap Y$: 
Then all the conditions in Theorem \ref{hc3} except the denseness of 
$U_{\ol{Y}}$ in $\ol{Y}$ remains true. (Note that $U_{\ol{Y}}$ is still 
dense in $Y$. So, if $Y$ is dense 
in $\ol{Y}$, $U_{\ol{Y}}$ remains to be dense in $\ol{Y}$.) \par
Let us take an $(U_{\tiol{Y}},\tiol{Y})$-triple 
$(Z,\ol{Z},\cZ)$ over $(S,S,\cS)$ as in 
the statement of the theorem. Let us denote the open immersion 
$Z \hra \ol{Z}$ by $j_Z$. We prove the following claim: \\
\quad \\
{\bf claim.} \,\,\, 
$R^qf_{(X \times_Y Z, \ol{X} \times_{\ol{Y}} \ol{Z})/\cZ, \rig *}\cE$ 
is a coherent $j_Z^{\d}\cO_{]\ol{Z}[_{\cZ}}$-module. \\
\quad \\
Since we have 
$R^qf_{(X \times_Y Z, \ol{X} \times_{\ol{Y}} \ol{Z})/\cZ, \rig *}\cE=0$ for 
$q>q_0$, we may assume $0 \leq q \leq q_0$. Let us put 
$E_Z := \text{(inverse image of $E$)}_{\red} \subseteq \ol{Z}$. 
By assumption on the triple $(Z,\ol{Z},\cZ)$, it is a simple normal 
crossing divisor in $\ol{Z}$. Denote the log structure associated to 
$E_Z$ by $M_{\ol{Z}}$. Since the claim is Zariski local on $\cZ$, 
we may assume that there exists a fine log structure $M_{\cZ}$ on 
$\cZ$ such that the closed immersion $\ol{Z} \hra \cZ$ comes from an 
exact closed immersion $(\ol{Z},M_{\ol{Z}}) \hra (\cZ,M_{\cZ})$. \par 
Let us put $\tiol{X} := \ol{X} \times_{\ol{Y}} \tiol{Y}, 
U_{\tiol{X}} := \ol{X} \times_{\ol{Y}} U_{\tiol{Y}}$. Also, 
we put $\brol{X}^{(n)\lr} := \cosk_{q_2}^{\tiol{X}^{(n)}} \tiol{X}^{(n)\lr}$
and denote the inverse image of $U_{\tiol{X}^{(n)}}$ in 
$\brol{X}^{(n)\lr}$ by $U_{\brol{X}^{(n)\lr}}$. We denote the inverse image 
of $X$ in $\tiol{X}, \tiol{X}^{(n)}, \tiol{X}^{(n)\lr}, \brol{X}^{(n)\lr}$ 
by $\ti{X}, \ti{X}^{(n)},\ti{X}^{(n)\lr}, \br{X}^{(n)\lr}$, 
respectively. With this notation, we have 
$(X \times_Y Z, \ol{X} \times_{\ol{Y}} \ol{Z}) = 
((\ti{X} \cap U_{\tiol{X}}) \times_{\tiol{Y}} Z, 
\tiol{X} \times_{\tiol{Y}} \ol{Z})$. \par 
Since $(X^{(\b)},\ol{X}^{(\b)}) \lra (X,\ol{X})$ is a proper hypercovering, 
so is 
$(\ti{X}^{(\b)} \cap U_{\tiol{X}^{(\b)}}, \tiol{X}^{(\b)}) \allowbreak 
\lra 
(\ti{X} \cap U_{\tiol{X}},\tiol{X})$. Hence 
$$((\ti{X}^{(\b)} \cap U_{\tiol{X}^{(\b)}}) \times_{\tiol{Y}} Z, 
\tiol{X}^{(\b)} \times_{\tiol{Y}} \ol{Z}) \lra 
((\ti{X} \cap U_{\tiol{X}}) \times_{\tiol{Y}} Z, \tiol{X} 
\times_{\tiol{Y}} \tiol{Z}) = (X \times_Y Z, \ol{X} \times_{\ol{Y}} \ol{Z})$$ 
is a proper hypercovering. So, by \cite{tsuzuki4}, 
we have the spectral sequence 
$$ 
E_1^{s,t} := 
R^tf_{((\ti{X}^{(s)} \cap U_{\tiol{X}^{(s)}}) \times_{\tiol{Y}} Z, 
\tiol{X}^{(s)} \times_{\tiol{Y}} \ol{Z})/\cZ, \rig*}\cE 
\,\Longrightarrow\, 
R^{s+t}f_{(X \times_Y Z, \ol{X} \times_{\ol{Y}} \ol{Z})/\cZ, \rig *}\cE. $$ 
So, to prove the claim, it suffices to prove that 
$R^tf_{((\ti{X}^{(s)} \cap U_{\tiol{X}^{(s)}}) \times_{\tiol{Y}} Z, 
\tiol{X}^{(s)} \times_{\tiol{Y}} \ol{Z})/\cZ, \rig*}\cE$ is a 
coherent $j_Z^{\d}\cO_{]\ol{Z}[_{\cZ}}$-module for 
$0 \leq s,t \leq q_1$. Next, since  
$(U_{\brol{X}^{(s)\lr}}, \brol{X}^{(s)\lr}) \lra 
(U_{\tiol{X}^{(s)}},\tiol{X}^{(s)})$ is a proper hypercovering for 
for $0 \leq s \leq q_1$, so is 
$(\br{X}^{(s)\lr} \cap U_{\brol{X}^{(s)\lr}}, \allowbreak \brol{X}^{(s)\lr}) 
\allowbreak \lra 
(\ti{X}^{(s)}\cap U_{\tiol{X}^{(s)}}, \tiol{X}^{(s)}).$ Hence 
$$
((\br{X}^{(s)\lr} \cap U_{\brol{X}^{(s)\lr}}) \times_{\tiol{Y}} Z, 
\brol{X}^{(s)\lr} \times_{\tiol{Y}} \ol{Z}) \lra 
(\ti{X}^{(s)}\cap U_{\tiol{X}^{(s)}} \times_{\tiol{Y}} Z, 
\tiol{X}^{(s)} \times_{\tiol{Y}} \ol{Z})$$ is a proper hypercovering. 
So we have the spectral sequence 
\begin{align*} 
E_1^{u,v} & := 
R^vf_{((\br{X}^{(s)\lur} \cap U_{\brol{X}^{(s)\lur}}) \times_{\tiol{Y}} Z, 
\brol{X}^{(s)\lur} \times_{\tiol{Y}} \ol{Z})/\cZ, \rig*}\cE \\ 
& \,\Longrightarrow\, 
R^{u+v}f_{((\ti{X}^{(s)}\cap U_{\tiol{X}^{(s)}}) \times_{\tiol{Y}} Z, 
\tiol{X}^{(s)} \times_{\tiol{Y}} \ol{Z})/\cZ, \rig *}\cE.
\end{align*}
So, to prove the claim, it suffices to prove that
$
R^vf_{((\br{X}^{(s)\lur} \cap U_{\brol{X}^{(s)\lur}}) \times_{\tiol{Y}} Z, 
\brol{X}^{(s)\lur} \times_{\tiol{Y}} \ol{Z})/\cZ, \rig*}\cE 
$
is a 
coherent $j_Z^{\d}\cO_{]\ol{Z}[_{\cZ}}$-module for 
$0 \leq s \leq q_1, 0 \leq u,v \leq q_2$. For such $s,u$, we have 
$\br{X}^{(s)\lur} \cap U_{\brol{X}^{(s)\lur}} = 
\ti{X}^{(s)\lur} \cap U_{\tiol{X}^{(s)\lur}}, 
\brol{X}^{(s)\lur}= \tiol{X}^{(s)\lur}$. So 
it suffices to prove that
$
R^vf_{((\ti{X}^{(s)\lur} \cap U_{\tiol{X}^{(s)\lur}}) \times_{\tiol{Y}} Z, 
\tiol{X}^{(s)\lur} \times_{\tiol{Y}} \ol{Z})/\cZ, \rig*}\cE 
$
is a 
coherent $j_Z^{\d}\cO_{]\ol{Z}[_{\cZ}}$-module for 
$0 \leq s \leq q_1, 0 \leq u,v \leq q_2$. \par 
Now let us take a log blow-up 
$\psi': (\ol{\X},M^{\ve}_{\ol{\X}}) \lra 
(\tiol{X}^{(s)\lur},M^{\ve}_{\tiol{X}^{(s)\lur}}) 
\times_{(\tiol{Y},M_{\tiol{Y}})} (\ol{Z},M_{\ol{Z}})$ 
such that $\ol{\X}$ is regular. Then $M^{\ve}_{\ol{\X}}$ corresponds to 
a normal crossing divisor on $\X$, which we denote by $\D^{\ve}$. 
Let us define 
$(\ol{\X},M_{\ol{\X}})$ in order that the following diagram is Cartesian: 
\begin{equation*}
\begin{CD}
(\ol{\X},M_{\ol{\X}}) @>{\psi}>> 
(\tiol{X}^{(s)\lur},M_{\tiol{X}^{(s)\lur}}) 
\times_{(\tiol{Y},M_{\tiol{Y}})} (\ol{Z},M_{\ol{Z}}) \\ 
@VVV @VVV \\ 
(\ol{\X},M^{\ve}_{\ol{\X}}) @>{\psi'}>> 
(\tiol{X}^{(s)\lur},M^{\ve}_{\tiol{X}^{(s)\lur}}) 
\times_{(\tiol{Y},M_{\tiol{Y}})} (\ol{Z},M_{\ol{Z}}). 
\end{CD}
\end{equation*}
Let us denote the pull-back of $D^{(s)\lur} \times_{\tiol{Y}} \ol{Z}$ in 
$\ol{\X}$ by $\D^{\ho}$. Since $D^{(s)\lur}$ is a log normal crossing 
divisor in $\tiol{X}^{(s)\lur}$, 
$D^{(s)\lur} \times_{\tiol{Y}} \ol{Z}$ is a log normal crossing divisor 
in $\tiol{X}^{(s)\lur} \times_{\tiol{Y}} \ol{Z}$ by Proposition 
\ref{lncdpullback}. Then, 
by definition of $\psi$ above, we see that 
$\D^{\ho} \cup \D^{\ve}$ is the normal crossing divisor corresponding to 
$M_{\ol{\X}}$, $\D^{\ho}, \D^{\ve}$ are sub normal crossing divisors of it. 
Moreover, $(\ol{\X},M_{\ol{\X}})$ is log smooth over $(\ol{Z},M_{\ol{Z}})$, 
and if we define $\D_{[i]}$ by 
\begin{align*}
\D_{[0]} & :=\ol{\X}, \qquad 
\D_{[1]} := \text{the normalization of $\D^{\ho}$}, \\
\D_{[i]} & := \text{$i$-fold fiber product of $\D_{[1]}$ over $\ol{\X}$}, 
\end{align*}
$(\D_{[i]},M_{\ol{\X}}^{\ve}|_{\D_{[i]}})$ is log smooth over 
$(\ol{Z},M_{\ol{Z}})$ for any 
$i \in \N$. That is, the morphism $(\ol{\X},M_{\ol{\X}}) \lra 
(\ol{Z},M_{\ol{Z}})$ and $\D^{\ho},\D^{\ve}$ are as in the situation in 
the beginning of Section 2. \par 
Next, note that we have 
$\ti{X}^{(s)\lur} \cap 
(\tiol{X}^{(s)\lur},M^{\ve}_{\tiol{X}^{(s)\lur}})_{\triv} = 
(\tiol{X}^{(s)\lur}-D^{(s)\lur}) \cap 
(\tiol{X}^{(s)\lur},M^{\ve}_{\tiol{X}^{(s)\lur}})_{\triv} = 
(\tiol{X}^{(s)\lur},M_{\tiol{X}^{(s)\lur}})_{\triv}.$ 
This implies that the inverse image of $X^{(s)}$ by the morphism 
$\tiol{X}^{(s)\lur} \times_{\tiol{Y}} \ol{Z} \lra \ol{X}^{(s)}$ 
contains 
$
((\tiol{X}^{(s)\lur},M_{\tiol{X}^{(s)\lur}}) 
\times_{(\tiol{Y},M_{\tiol{Y}})} (\ol{Z},M_{\ol{Z}}))_{\triv}
$. So the diagram 
$$ \ol{\X} \os{\psi}{\lra} 
\tiol{X}^{(s)\lur} \times_{\tiol{Y}} \ol{Z} \lra \ol{X}^{(s)}$$
comes from the diagram of log schemes 
$$ 
(\ol{\X},M_{\ol{\X}}) \os{\psi}{\lra}  
(\tiol{X}^{(s)\lur},M_{\tiol{X}^{(s)\lur}}) 
\times_{(\tiol{Y},M_{\tiol{Y}})} (\ol{Z},M_{\ol{Z}}) \lra 
(\ol{X}^{(s)},M_{\ol{X}^{(s)}}). 
$$ 
Then the locally free isocrystal $\ol{\cE}$ on 
$(\ol{X}^{(s)}/\cS)^{\log}_{\conv}$ having nilpotent residues induces, by 
pull-back, a locally free isocrystal on 
$((\tiol{X}^{(s)\lur} \times_{\tiol{Y}} \ol{Z})/\cS)^{\log}_{\conv}$, 
and the induced locally free isocrystal 
on $(\ol{\X}/\cS)^{\log}_{\conv}$ (which 
we denote also by $\ol{\cE}$, by abuse of notation) 
has nilpotent residues by Proposition \ref{functnilp}. 
In particular, the diagram 
$$ 
(\ol{\X},M_{\ol{\X}}) \os{\psi}{\lra}  
(\tiol{X}^{(s)\lur},M_{\tiol{X}^{(s)\lur}}) 
\times_{(\tiol{Y},M_{\tiol{Y}})} (\ol{Z},M_{\ol{Z}}) \lra 
(\ol{Z},M_{\ol{Z}}) \hra (\cZ,M_{\cZ})
$$ 
and $\ol{\cE}$ are in the situation of Proposition \ref{starok} and 
Remark \ref{starokrem}. \par 
Finally, note that we have 
$$ 
\ti{X}^{(s)\lur} \cap U_{\tiol{X}^{(s)\lur}} = 
(\tiol{X}^{(s)\lur}-D^{(s)\lur}) \cap U_{\tiol{X}^{(s)\lur}} 
\subseteq (\tiol{X}^{(s)\lur},M_{\tiol{X}^{(s)\lur}})_{\triv} 
$$ 
and $Z \subseteq (\ol{Z},M_{\ol{Z}})_{\triv}$ (which follows from 
the equality $U_{\tiol{Y}} \cap E = \emptyset$.) So, on 
$
\ti{X}^{(s)\lur} \cap U_{\tiol{X}^{(s)\lur}} \times_{\tiol{Y}} Z
$, 
the log structure of 
$(\tiol{X}^{(s)\lur},M_{\tiol{X}^{(s)\lur}}) 
\times_{(\tiol{Y},M_{\tiol{Y}})} (\ol{Z},M_{\ol{Z}})$ is trivial. Hence 
the morphism $\psi$ is identity on it. That is, we have the isomorphism 
$$ 
\ol{\X} \times_{\ol{Z}} Z - \D^{\ho} = 
(\ol{\X}-\D^{\ho}) \times_{\tiol{Y}} Z 
\os{=}{\lra} 
(\tiol{X}^{(s)\lur}-D^{(s)\lur}) \times_{\tiol{Y}} Z 
= 
(\ti{X}^{(s)\lur} \cap U_{\tiol{X}}) \times_{\tiol{Y}} Z. 
$$ 

Let us denote the open immersion 
$\ol{\X} \times_{\ol{Z}} Z - \D^{\ho} \hra \ol{\X}$ by $j_{\X}$. 
Then, by combining the above results, we obtain the equalities 
{\allowdisplaybreaks{
\begin{align*}
& \phantom{=} 
R^vf_{((\ti{X}^{(s)\lur} \cap U_{\tiol{X}^{(s)\lur}}) \times_{\tiol{Y}} 
\ol{Z}, \tiol{X}^{(s)\lur} \times_{\tiol{Y}} \ol{Z})/\cZ, \rig *}\cE \\
& = 
R^vf_{(\ol{X}\times_{\ol{Z}}Z - \D^{\ho}, \ol{\X})/\cZ, \rig *}\cE \\ 
& = 
R^vf_{(\ol{X}\times_{\ol{Z}}Z - \D^{\ho}, \ol{\X})/\cZ, \rig *}
j_{\X}^{\dd}\ol{\cE} \\ 
& = 
j_Z^{\d}R^vf_{\ol{X}/\cZ, \an *}\ol{\cE} \,\,\,\, 
\text{(Theorem \ref{thm2.2}, 
Proposition \ref{starok}, Remark \ref{starokrem})}
\end{align*}}}
and it is a coherent $j^{\d}_Z\cO_{]\ol{Z}[_{\cZ}}$-module. 
So we have proved the claim. \par 
By using the claim, we see also that 
$R^q
f_{(X \times_Y Z, \ol{X} \times_{\ol{Y}} \ol{Z})/\cZ \times_{\cS} \cZ, \rig *}
\cE$ is a coherent $j_Z\cO_{]\ol{Z}[_{\cZ \times_{\cS}\cZ}}$-module. 
So, by \cite[2.3.1]{tsuzuki3}, we see that the morphisms 
$$ p_i^*R^qf_{(X \times_{Y} Z, 
\ol{X} \times_{\ol{Y}} \ol{Z})/\cZ, \rig *}\cE 
\lra 
R^qf_{(X \times_{Y} Z, \ol{X} \times_{\ol{Y}} \ol{Z})/\cZ \times_{\cS} \cZ, 
\rig *}\cE \,\, (i=1,2)$$ 
in the definition of $\epsilon$ in the statement of the 
theorem are isomorphisms, and the pair 
$(R^qf_{(X \times_{Y} Z, 
\ol{X} \times_{\ol{Y}} \ol{Z})/\cZ, \rig *}\cE, \epsilon)$ defines 
an object in $I^{\d}((Z,\ol{Z})/\cS_K, \cZ)$. We denote it by 
$\cF_{\cZ}$. \par 
Now let us take a hypercovering $(U_{\ol{Y}},\ol{Y}) \lla 
(U_{\tiol{Y}^{(\b)}},\tiol{Y}^{(\b)})$ with 
$(U_{\tiol{Y}^{(0)}},\tiol{Y}^{(0)}) = (U_{\tiol{Y}},\tiol{Y})$ such that 
each $\tiol{Y}^{(n)}$ is quasi-projective smooth over $S$ and that 
$\tiol{Y}^{(n)}-U_{\tiol{Y}^{(n)}}$ is a simple normal crossing divisor 
in $\tiol{Y}^{(n)}$. For each $n \leq 2$, let us take a closed immersion 
$\tiol{Y}^{(n)} \hra {\cY'}^{(n)}$ such that ${\cY'}^{(n)}$ is formally 
smooth over $\cS$. (It is possible because $\tiol{Y}^{(n)}$ is 
quasi-projective over $S$.) Then, following \cite[7.3.1]{tsuzuki5}, 
we define the simplicial formal scheme 
$\Gamma^{\cS}_m({\cY'}^{(m)})^{\leq m}$ by 
$$ [n] \mapsto \prod_{\gamma:[m] \ra [n]} {\cY'}^{(m)}, $$ 
where $[n]$ denotes the set $\{0,1,\cdots,n\}$, $\gamma$ runs through 
the non-decreasing maps $[m] \ra [n]$ and the transition maps are defined 
in natural way. Using this, we define the $2$-truncated simplicial formal 
scheme $\cY^{(\b)}$ by 
$$ \cY^{(\b)} := 
\sk_2(\prod_{0 \leq m \leq 2} \Gamma^{\cS}_m({\cY'}^{(m)})^{\leq m}) 
$$ 
with the transition maps induced by those of 
$\Gamma^{\cS}_m({\cY'}^{(m)})^{\leq m}$'s. Then, for $0 \leq n \leq 2$, 
the triple $\cY^{(n)}:=(U_{\tiol{Y}^{(n)}}, \tiol{Y}^{(n)}, \cY^{(n)})$ 
is a $(U_{\tiol{Y}},\tiol{Y})$-triple satisfying the condition in the 
statement of the theorem required for $(Z,\ol{Z},\cZ)$. So we have 
overconvergent isocrystals $\cF_{\cY^{(n)}}$ on 
$(U_{\tiol{Y}^{(n)}}, \tiol{Y}^{(n)})/\cS_K$ for $n=0,1,2$. Moreover, by 
\cite[2.3.1]{tsuzuki3}, it is compatible with respect to $n$. So, by 
proper descent for overconvergent isocrystals (\cite[7.3]{shiho4} and 
\cite[3.3.4.2]{saintdonat}), the comptible family 
$\{\cF_{\cY^{(n)}}\}_{n=0,1,2}$ descents to an overconvergent isocrystal 
on $(U_{\ol{Y}},\ol{Y})/\cS_K$, which we denote by $\cF$. \par 
The uniqueness of the overconvergent isocrystal $\cF$ follows from 
the uniqueness of the comptible family 
$\{\cF_{\cY^{(n)}}\}_{n=0,1,2}$ and proper descent for overconvergent 
isocrystals. We check that, for any $(Z,\ol{Z},\cZ)$ as in the statement of 
the theorem, the restrction of $\cF$ to $I^{\d}((Z,\ol{Z})/\cS_K,\cZ)$ 
is given by 
$\cF_{\cZ}:=
(R^qf_{(X \times_Y Z,\ol{X} \times_{\ol{Y}}\ol{Z})/\cZ,\rig *}\cE, \epsilon)$ 
as in the statement of the theorem. 
(The proof is similar to \cite[4.8]{shiho3}.) 
Let us take $(Z,\ol{Z},\cZ)$ as above and let us consider the following 
diagram 
\begin{equation}\label{express}
\begin{CD}
I^{\d}((U_{\ol{Y}},\ol{Y})/\cS_K) @>>> 
I^{\d}((U_{\tiol{Y}},\tiol{Y})/\cS_K) @>=>>
I^{\d}((U_{\tiol{Y}},\tiol{Y})/\cS_K, \cY^{(0)}) \\ 
@. @VVV @VVV \\ 
@. I^{\d}((Z,\ol{Z})/\cS_K,\cZ) @>=>> 
I^{\d}((Z,\ol{Z})/\cS_K, \cZ \times_{\cS} \cY^{(0)}), 
\end{CD}
\end{equation}
where all the functors are restrictions. Then, by definition, $\cF$ is 
sent to $\cF_{\cY^{(0)}}$ by the composite of horizontal arrows and 
it is sent to $\cF_{\cZ \times_{\cS} \cY^{(0)}}$ by the right vertical 
arrow (by functoriality of relative rigid cohomology). On the other hand, 
the overconvergent isocrystal $\cF_{\cZ}$ is sent to 
$\cF_{\cZ \times_{\cS} \cY^{(0)}}$ by the lower horizontal arrow 
(again by functoriality of relative rigid cohomology). So, by 
the commutativity of the diagram \eqref{express}, we see that the restriction 
of $\cF$ to $I^{\d}((Z,\ol{Z})/\cS_K,\cZ)$ is given by $\cF_{\cZ}$. \par 
Finally we prove the functoriality of the above expression. 
Let $\varphi:(Z',\ol{Z}',\cZ') \lra (Z,\ol{Z},\cZ')$ be a morphism 
of $(U_{\tiol{Y}},\tiol{Y})$-triples over $(S,S,\cS)$ satisfying the 
condition in the statement of the theorem required for $(Z,\ol{Z},\cZ)$. 
Then we have two morphisms of the form 
\begin{equation}\label{functex}
\varphi^*\cF_{\cZ} \lra \cF_{\cZ'}:
\end{equation}
One is the morphism induced by the functoriality of the relative rigid 
cohomology and the other is the morphism induced by the isocrystal 
structure of $\cF$. We should prove that they are equal. 
Let us consider the following commutative diagram 
\begin{equation*}
\begin{CD}
I^{\d}((Z,\ol{Z})/\cS_K, \cZ) @>{\pi^*}>> 
I^{\d}((Z,\ol{Z})/\cS_K,\cZ \times_{\cS} \cY^{(0)}) \\ 
@V{\varphi^*}VV @V{{\varphi'}^*}VV \\ 
I^{\d}((Z',\ol{Z}')/\cS_K, \cZ') @>{{\pi'}^*}>> 
I^{\d}((Z',\ol{Z}')/\cS_K,\cZ' \times_{\cS} \cY^{(0)}), 
\end{CD}
\end{equation*}
where $\varphi^*, {\varphi'}^*$ are the restriction functors induced by 
$\varphi$ and $\pi,\pi'$ are the restriction functors induced by 
$\cZ \times_{\cS} \cY^{(0)} \lra \cZ, 
\cZ' \times_{\cS} \cY^{(0)} \lra \cZ'$, respectively. Then, by looking at 
the argument in the previous paragraph, we see that either of the morphisms 
\eqref{functex} fits into the following diagram 
\begin{equation*}
\begin{CD}
{\pi'}^*\varphi^*\cF_{\cZ} @>>> {\varphi'}^*\cF_{\cZ \times_{\cS}\cY^{(0)}} \\ 
@VVV @VVV \\
{\pi'}^*\cF_{\cZ'} @>>> \cF_{\cZ' \times_{\cS}\cY^{(0)}}, 
\end{CD}
\end{equation*}
where the left vertical arrow is the pull-back of the morphism 
\eqref{functex} by ${\pi'}^*$ and the other arrows are induced by the 
functoriality of relative rigid cohomology. Since ${\pi'}^*$ is 
an equivalence of categories, we see that the two morphisms \eqref{functex} 
are equal. So we have proved the functoriality and the proof of the 
theorem is now finished. 
\end{pf} 

By using Theorem \ref{mainthm} repeatedly, we obtain the 
following corollary, which says a kind of constructibility 
of relative rigid cohomology: 

\begin{cor}\label{constr}
Let us assume given a morphism of pairs 
$f: (X,\ol{X}) \lra (Y,\ol{Y})$ such that $f:\ol{X} \lra \ol{Y}$ 
is proper. 
Then, there exists a stratification $Y = \coprod_{i=0}^d Y_i$ of $Y$ 
by locally closed subschemes $Y_i$ with $Y'_{i} = \coprod_{j\geq i}Y_{j}$ 
$($where $Y'_{i}$ denotes the closure of $Y_i$ in $Y)$
and proper surjective maps $\tiol{Y}_i \lra \ol{Y}_i \,(0 \leq i \leq d)$ 
$($where $\ol{Y}_i$ denotes the closure of $Y_i$ in $\ol{Y})$ 
and a simple normal crossing 
divisor $E_i$ on $\tiol{Y}_i$ with $E_i \cap \ti{Y}_i = \emptyset$ 
$($where we put $\ti{Y}_i := Y_i \times_{\ol{Y}_i} \tiol{Y}_i)$ 
such that, 
for any potentially semi-stable overconvergent isocrystal $\cE$ on 
$(X,\ol{X})/\cS_K$ and $q \in \N$, there exist uniquely 
the overconvergent isocrystals $\cF_i$ on $(Y_i,\ol{Y}_i)/\cS_K \,(0 \leq i 
\leq d)$ 
satisfying the following condition$:$ 
For any $(\ti{Y}_i, \tiol{Y}_i)$-triple $(Z,\ol{Z},\cZ)$ over 
$(S,S,\cS)$ such that $\ol{Z}$ is smooth over $S$, $(E_i \times_{\tiol{Y}_i} 
\ol{Z})_{\red}$ is contained in a simple normal crossing divisor on 
$\ol{Z}$ and that $\cZ$ is formally smooth over $\cS$, 
the restriction of $\cF_i$ to $I^{\d}((Z,\ol{Z})/\cS_K, \cZ)$ 
is given functorially by 
$(R^qf_{(X \times_{Y} Z, \ol{X} \times_{\ol{Y}} \ol{Z})/\cZ, \rig *}\cE, 
\epsilon)$, where $\epsilon$ is an isomorphism 
{\small{ $$ 
p_2^*R^qf_{(X \times_{Y} Z, 
\ol{X} \times_{\ol{Y}} \ol{Z})/\cZ, \rig *}\cE 
\os{\sim}{\ra} 
R^qf_{(X \times_{Y} Z, \ol{X} \times_{\ol{Y}} \ol{Z})/\cZ \times_{\cS} \cZ, 
\rig *}\cE \os{\sim}{\leftarrow}
p_1^*R^qf_{(X 
\times_{Y} Z, \ol{X} \times_{\ol{Y}} \ol{Z})/\cZ, \rig *}\cE
$$ }}
$(p_j$ denotes the $j$-th projection 
$]\ol{Z}[_{\cZ \times_{\cS} \cZ} \lra ]\ol{Z}[_{\cZ})$. 
\end{cor} 

\begin{pf} 
We may replace $\ol{Y}$ by the closure of $Y$ in $\ol{Y}$. 
As we remarked in the proof in Theorem \ref{mainthm}, we can take 
the open set $U_{\ol{Y}}$ in Theorem \ref{mainthm} to be 
dense in $\ol{Y}$ if $Y$ is dense in $\ol{Y}$. So we can take a 
dense open set $Y_0 \subseteq Y$ such that the conclusion of the corollary 
holds for $i=0$. 
Then put $Y'_1 := Y-Y_0$, $\ol{Y}_1$ its closure in $\ol{Y}$ and 
apply Theorem \ref{mainthm} to the pull-back of $f$ by 
$(Y'_1,\ol{Y}_1) \lra (Y,\ol{Y})$: Then we have an open dense subscheme 
$Y_1$ in $Y'_1$ such that the conclusion of the corollary is true for 
$i=1$. Repeating this process, we can prove the corollary. 
\end{pf} 

Combining with Theorem \ref{sconj}, we obtain the following corollary. 

\begin{cor}\label{maincor}
Let us assume given a morphism of pairs 
$f: (X,\ol{X}) \lra (Y,\ol{Y})$ such that $f:\ol{X} \lra \ol{Y}$ 
is proper. Then, 
\begin{enumerate}
\item 
There exists a dense open set 
$U_{\ol{Y}}$ of $Y$, a proper surjective map $\tiol{Y} \lra \ol{Y}$ 
and a simple normal crossing 
divisor $E$ on $\tiol{Y}$ with $E \cap U_{\tiol{Y}} = \emptyset$ 
$($where we put $U_{\tiol{Y}} := U_{\ol{Y}} \times_{\ol{Y}} 
\tiol{Y})$ such that, 
for any overconvergent $F$-isocrystal $\cE$ on 
$(X,\ol{X})/\cS_K$ and $q \in \N$, there exists 
the unique overconvergent isocrystal $\cF$ on $(U_{\ol{Y}},\ol{Y})/\cS_K$ 
satisfying the condition of Theorem \ref{mainthm}.
\item  
Moreover, there exists a stratification $Y = \coprod_{i=0}^d Y_i$ of $Y$ 
by locally closed subschemes $Y_i$ with $Y'_{i} = \coprod_{j\geq i}Y_{j}$ 
$($where $Y'_{i}$ denotes the closure of $Y_i$ in $Y)$
and proper surjective maps $\tiol{Y}_i \lra \ol{Y}_i \,(0 \leq i \leq d)$ 
$($where $\ol{Y}_i$ denotes the closure of $Y_i$ in $\ol{Y})$ 
and a simple normal crossing 
divisor $E_i$ on $\tiol{Y}_i$ with $E_i \cap \ti{Y}_i = \emptyset$ 
$($where we put $\ti{Y}_i := Y_i \times_{\ol{Y}_i} \tiol{Y}_i)$ 
such that, 
for any overconvergent $F$-isocrystal $\cE$ on 
$(X,\ol{X})/\cS_K$ and $q \in \N$, there exist uniquely 
the overconvergent isocrystals 
$\cF_i$ on $(Y_i,\ol{Y}_i)/\cS_K \,(0 \leq i 
\leq d)$ 
satisfying the conclusion of Corollary \ref{constr}. 
\end{enumerate}
Furthermore, 
the overconvergent isocrystals $\cF, \cF_i$ have Frobenius structures 
which are induced by the Frobenius structure of $\cE$. 
\end{cor} 

\begin{pf} 
The assertions except the last one 
follow immediately from Theorem \ref{mainthm}, Corollary \ref{constr} and 
Theorem \ref{sconj}. As for the last assertion, we can argue as in 
the proof of \cite[5.16]{shiho3} (see also Theorem \ref{thm2.4}) to 
reduce to the case $Y=\ol{Y}$ is equal to $S$, and in this case, 
the assertion follows from the bijectivity of the Frobenius endomorphism 
of the (absolute) rigid cohomology: In smooth case, it is shown in 
Theorem \ref{thm2.4} and we can reduce the general case to the smooth case 
by using proper descent of rigid cohomology. So we are done. 
\end{pf}



\begin{thebibliography}{[KKMS-]}

\bibitem[A]{andre}
  Y. ~Andr\'{e}, 
  {\it Filtrations de type Hasse-Arf et monodromie $p$-adique}, 
  Invent Math. {\bf 148}(2002), 285--317. 

\bibitem[Ba-Ch]{bc}
   F. ~Baldassarri and B. ~Chiarellotto, 
   {\it Formal and $p$-adic theory of differential 
   systems with logarithmic singularities depending 
   upon parameters}, 
   Duke Math. J. {\bf 72}(1993), 241--300.

\bibitem[Ba-Ch2]{bc2} 
   F. ~Baldassarri and B. ~Chiarellotto, 
   {\it Algebraic versus rigid Cohomology with logarithmic coefficients},
   in Barsotti Symposium in Algebraic Geometry, Academic Press.

\bibitem[Be1]{berthelot1} 
   P. ~Berthelot, 
   {\it Cohomologie cristalline des sch\'{e}mas de 
   caract\'{e}ristique $p > 0$}, 
   Lecture Note in Math. {\bf 407}, Springer-Verlag (1974).

\bibitem[Be2]{berthelot2}
   P. ~Berthelot,  
   {\it G\'{e}om\'{e}trie rigide et cohomologie des vari\'{e}t\'{e}s
    alg\'{e}briques de caract\'{e}ristique p},
   Bull. Soc. Math. de France, M\'{e}moire {\bf 23}(1986), 7--32. 

\bibitem[Be3]{berthelot3} 
  P. ~Berthelot, 
  {\it Cohomologie rigide et cohomologie rigide \`{a} supports propres
   \,\,\, premi\`{e}re partie}, pr\'{e}publication de 
   l'IRMAR 96-03.

\bibitem[Be4]{berthelot4}
  P. ~Berthelot, 
  {\it Finitude et puret\'{e} cohomologique en cohomologie rigide}, 
  Invent. Math., {\bf 128}(1997), 329--377.

\bibitem[Be-O]{berthelotogus}
   P. ~Berthelot and A. ~Ogus,
   {\it Notes on crystalline cohomology},
   Mathematical Notes, Princeton University Press, 1978.

\bibitem[Ch-T]{chts}
   B. ~Chiarellotto and N. Tsuzuki, 
   {\it Cohomological descent of rigid cohomology for etale 
   coverings}, Rend. ~Sem. ~Mat. ~Univ. Padova, {\bf 109}(2003), 
   63--215. 



\bibitem[I]{illusie} 
  L. ~Illusie, 
  {\it An overview of the work of K. ~Fujiwara, K. ~Kato and C. ~Nakayama 
  on logarithmic \'{e}tale cohomology}, 
  Ast\'{e}risque, {\bf 279}(2002), 271--332. 

\bibitem[dJ1]{dejong1}
   A. ~J. ~de Jong, 
   {\it Smoothness, semi-stability and alterations}, 
   Publ. Math. IHES, {\bf 83}(1996), 51--93. 

\bibitem[dJ2]{dejong2} 
   A. ~J. ~de Jong, 
   {\it Families of curves and alterations}, 
   Ann. Inst. Fourier, {\bf 47}(1997), 599--621. 

\bibitem[Ka1]{kkato1}
   K. ~Kato,
   {\it Logarithmic structures of Fontaine-Illusie},
   in Algebraic Analysis, Geometry, and Number Theory,
    J-I.Igusa ed., 1988, Johns Hopkins University, pp.\ 191--224.

\bibitem[Ka2]{kkato2} 
   K. ~Kato, 
   {\it Toric singularities}, 
   Amer. J. Math. {\bf 116}(1994), 1073--1099. 

\bibitem[Ke1]{kedlaya-p}
  K. ~S. ~Kedlaya, 
  {\it A $p$-adic local monodromy theorem}, 
  Annals of Math. {\bf 160}(2004), 93--184. 

\bibitem[Ke2]{kedlaya1}
   K. ~S. ~Kedlaya, 
   {\it Full faithfulness for overconvergent $F$-isocrystals}, 
   Geometric Aspects of Dwork Theory, Walter de Gruyter, 819--835 (2004). 

\bibitem[Ke3]{kedlaya2} 
   K. ~S. ~Kedlaya, 
   {\it Finiteness of rigid cohomology with coefficients}, 
   Duke Math. J. {\bf 134}(2006), 15--97. 

\bibitem[Ke4]{kedlayaI}
   K. ~S. ~Kedlaya, 
   {\it Semistable reduction for overconvergent $F$-isocrystals, I$:$ 
   Unipotence and logarithmic extensions}, 
   Compositio Math.. {\bf 143}(2007), 1164--1212. 

\bibitem[Ke5]{kedlayaII}
   K. ~S. ~Kedlaya, 
   {\it Semistable reduction for overconvergent $F$-isocrystals, II$:$ 
   A valuation-theoretic approach}, 
   to appear in Compositio Math.

\bibitem[Ke6]{kedlayaIII} 
   K. ~S. ~Kedlaya, 
   {\it Semistable reduction for overconvergent $F$-isocrystals, III$:$ 
   Local semistable reduction at monomial valuations}, preprint. 

\bibitem[Ke7]{kedlayaIV}
   K. ~S. ~Kedlaya, 
   {\it Semistable reduction for overconvergent $F$-isocrystals, IV$:$ 
   Local semistable reduction at nonmonomial valuations}, preprint. 

\bibitem[KKMS]{kkms} 
   G. ~Kempf, F. ~Knudsen, D. ~Mumford and B. ~Saint-Donat, 
   {\it Toroidal embeddings I}, Lecture Note in Math. {\bf 339}, 
   Springer Verlag (1973). 

\bibitem[M]{mebkhout} 
   Z. ~Mebkhout, 
   {\it Analogue $p$-adique du th\'{e}or\`{e}me de Turrittin et le 
   th\'{e}or\`{e}me de la monodromie $p$-adique}, 
   Invent Math. {\bf 148}(2002), 319--351. 

\bibitem[Na-Sh]{nakkshiho}
  Y. ~Nakkajima and A. ~Shiho, 
  {\it Weight filtrations on log crystalline cohomologies of 
   families of open smooth varieties in characteristic $p>0$}, 
   to appear in Lecture Note in Math. 

\bibitem[Ni]{niziol} 
  W. ~Niziol, 
  {\it Toric singularities: Log-blow-ups and global resolutions}, 
  J. Alg. Geom. {\bf 15}(2006), 1--29. 

\bibitem[O1]{ogus1}
   A. ~Ogus,
   {\it F-isocrystals and de Rham Cohomology II  --- Convergent
   Isocrystals}, Duke Math. J., {\bf 51}(1984), 765--850.

\bibitem[O2]{ogus2}
   A. ~Ogus,
   {\it The Convergent Topos in Characteristic p},
   in Grothendieck Festschrift, Progress in Math., Birkh\"{a}user.


\bibitem[SD]{saintdonat} 
   B. Saint-Donat, 
   {\it Techniques de descente cohomologique}, 
   Expos\'{e} V bis in 
   Th\'{e}orie des topos et cohomologie \'{e}tale des 
   sch\'{e}mas, Lecture Note in Math. {\bf 270}, 
   Springer-Verlag. 
   
\bibitem[Sa]{saito}
   T. ~Saito, 
   {\it Log smooth extension of a family of curves and semi-stable 
   reduction}, 
   J. Alg. Geom. {\bf 13}(2004), 287--321. 

\bibitem[Sh1]{shiho1}
   A. ~Shiho, 
   {\it Crystalline Fundamental Groups I ---
    Isocrystals on Log Crystalline Site and Log Convergent Site}, 
    J. ~Math. ~Sci. ~Univ. 
    ~Tokyo, {\bf 7}(2000), 509--656.

\bibitem[Sh2]{shiho2}
   A. ~Shiho, 
  {\it Crystalline fundamental groups II --- 
  Log convergent cohomology and rigid cohomology}.
  J. ~Math. ~Sci. ~Univ. ~Tokyo, {\bf 9}(2002),  1--163.

\bibitem[Sh3]{shiho3} 
  A. ~Shiho, 
  {\it Relative log convergent cohomology and relative 
  rigid cohomology I}, preprint. 

\bibitem[Sh4]{shiho4} 
  A. ~Shiho, 
  {\it Relative log convergent cohomology and relative 
  rigid cohomology II}, preprint. 


\bibitem[T1]{tsuzuki2} 
 N. ~Tsuzuki, 
 {\it On the Gysin isomorphism of rigid cohomology}, 
 Hiroshima Math. J., {\bf 29}(1999), 479--527. 

\bibitem[T2]{tsuzuki3}
 N. ~Tsuzuki, 
 {\it On base change theorem and coherence in rigid cohomology}, 
 Documenta Math., Extra Volume Kato (2003), 891--918. 

\bibitem[T3]{tsuzuki4}
 N. ~Tsuzuki, 
 {\it Cohomological descent of rigid cohomology for proper coverings}, 
Invent. Math. {\bf 151}(2003), 101-133. 

\bibitem[T4]{tsuzuki5} 
 N. ~Tsuzuki, 
 {\it Cohomological descent in rigid cohomology}, 
 in Geometric Aspects of Dwork Theory, Walter de Gruyter, 931-981 (2004). 

\end{thebibliography}
\end{document}